\RequirePackage{fix-cm}
\documentclass[11pt]{article}
\usepackage{graphicx,amsmath,amsfonts,amscd,amssymb,bm,cite,epsfig,epsf,url,color,algorithm,algorithmic,enumerate}
\usepackage{fullpage}
\usepackage[small,bf]{caption}
\usepackage{subfigure}
\usepackage{lscape}
\usepackage[space]{grffile}
\usepackage{epsfig,psfrag,amstext,subfigure,subeqnarray,nccfoots,color,multirow,bm}
\usepackage{appendix}

\newcommand\sss{\ensuremath{\bcl{s}}} 

\newcommand\KK{\ensuremath{\kappa}}

\newcommand\Ec{\ensuremath{\mathcal{E}}}
\newcommand\Gc{\ensuremath{\mathcal{G}}}
\newcommand\Nc{\ensuremath{\mathcal{N}}}

\newcommand\Yc{\ensuremath{\mathcal{Y}}}

\newcommand\Sc{\ensuremath{\mathcal{S}}}

\newcommand\Lc{\ensuremath{{\mathcal{L}}}}

\newcommand\xb{\ensuremath{{\bm x}}}
\newcommand\wb{\ensuremath{{\bm w}}}
\newcommand\yb{\ensuremath{{\bm y}}}

\newcommand\ub{\ensuremath{{\bm u}}}

\newcommand\Ab{\ensuremath{{\bm A}}}
\newcommand\ab{\ensuremath{{\bm a}}}
\newcommand\Bb{\ensuremath{{\bm B}}}
\newcommand\bb{\ensuremath{{\bm b}}}
\newcommand\Cb{\ensuremath{{\bm C}}}

\newcommand\db{\ensuremath{{\bm d}}}
\newcommand\Db{\ensuremath{{\bm D}}}

\newcommand\Eb{\ensuremath{{\bf E}}}

\newcommand\Ib{\ensuremath{{\bm I}}}
\newcommand\Lb{\ensuremath{{\bm L}}}
\newcommand\Pb{\ensuremath{{\bm P}}}
\newcommand\pb{\ensuremath{{\bm p}}}

\newcommand\Mb{\ensuremath{{\bf M}}}

\newcommand\qb{\ensuremath{{\bm q}}}

\newcommand\vb{\ensuremath{{\bm v}}}
\newcommand\Wb{\ensuremath{{\bm W}}}

\newcommand\zb{\ensuremath{{\bm z}}}

\newcommand\thetab{\ensuremath{{\bm \theta}}}

\newcommand\mub{\ensuremath{{\bm \mu}}}

\newcommand\zerob{\ensuremath{{\bm 0}}}

\newcommand\la{\ensuremath{{\langle}}}
\newcommand\ra{\ensuremath{{\rangle}}}
\newcommand\rpo{\ensuremath{{r+1}}}
\newcommand\rmo{\ensuremath{{r-1}}}

\newcommand\diag{\ensuremath{{\rm diag}}}
\newtheorem{Lemma}{Lemma}

\newtheorem{Theorem}{Theorem}
\newtheorem{Definition}{Definition}
\newtheorem{Corollary}{Corollary}

\newtheorem{Assumption}{Assumption}

\usepackage{algorithmic,algorithm}
\usepackage{array}
\usepackage{enumerate}

\def\MatrixFont{\bf}
\def\VectorFont{\bf}

\newcommand{\mA}{{\MatrixFont A}}

\newcommand{\mI}{{\MatrixFont I}}

\newcommand{\va}{{\VectorFont a}}
\allowdisplaybreaks[4]
\title{Decentralized Non-Convex Learning with Linearly Coupled Constraints: Algorithm Designs and Application to Vertical Learning Problem}
\title{\bf Decentralized Non-Convex Learning with Linearly Coupled Constraints: Algorithm Designs and Application to Vertical Learning Problem\thanks{ Jiawei Zhang and Songyang Ge contributed equally. Part of the work was presented in IEEE ICASSP 2020 \cite{ZhangICASSP20}. }}
\author{Jiawei Zhang \thanks{Jiawei Zhang is with Shenzhen Research Institute of Big Data, The Chinese University of Hong Kong, Shenzhen, China 518172. E-mail: 216019001@link.cuhk.edu.cn},
    Songyang Ge \thanks{Songyang Ge is with Shenzhen Research Institute of Big Data, The Chinese University of Hong Kong, Shenzhen, China 518172. E-mail: songyangge@link.cuhk.edu.cn},
    Tsung-Hui Chang \thanks{ Tsung-Hui Chang is the corresponding author.
    Address: Guangdong Provincial Key Laboratory of Big Data Computing and the School of Science and Engineering, The Chinese University of Hong Kong, Shenzhen, China 518172. E-mail: tsunghui.chang@ieee.org. },
    \ and \
    Zhi-Quan Luo
    \thanks{ Zhi-Quan Luo
    is with the Shenzhen Research Institute of Big Data, The Chinese University of Hong Kong, Shenzhen, China 518172.  E-mail: luozq@cuhk.edu.cn }
}
\begin{document}

\maketitle
\begin{abstract}
    Motivated by the need for decentralized learning, this paper aims at designing a distributed algorithm for solving nonconvex problems with general linear constraints over a multi-agent network. In the considered problem, each agent owns some local information and a local variable for jointly minimizing a cost function,
    but local variables are coupled by linear constraints.
    Most of the existing methods for such problems are only applicable for convex problems or problems with specific linear constraints.
    There still lacks a distributed algorithm for solving such problems with general linear constraints under the nonconvex setting. To tackle this problem, we propose a new algorithm, called \emph{proximal dual consensus} (PDC) algorithm, which combines a proximal technique and a dual consensus method.
   We show that under certain conditions the proposed PDC algorithm can generate an $\epsilon$-Karush-Kuhn-Tucker solution in $\mathcal{O}(1/\epsilon)$
   iterations, achieving the lower bound for distributed non-convex problems up to a constant. Numerical results are presented to demonstrate the good performance of the proposed algorithms for solving two vertical learning problems in machine learning over a multi-agent network.
\end{abstract}

\vspace{-0.0cm}
\section{Introduction}\vspace{-0.0cm}
\label{sec:intro}
Recently, distributed optimization has attracted significant attention in signal
processing, control and machine learning societies
\cite{giannakis2016decentralized,BK:Bekkerman12,chang2020distributed} due to the need of achieving
decentralized control and decision making in sensor networks, and large-scale data
processing and learning.
Since the distributed agents/nodes process only local data and messages exchanged
from direct neighbors, distributed optimization is well suited for applications
where
the data are acquired distributively  \cite{scutari2018parallel} or where data
privacy is of primary concern
\cite{konevcny2016federated}.

In this paper we are interested in distributed optimization methods for solving the
following
problem
\begin{align}
{\sf (P)}~~\min_{\substack{\xb_i \in \mathbb{R}^n, \\i=1,\cdots,N}}~ &\sum_{i=1}^N
f_i(\xb_i)  ~
~\text{s.t.}~ ~\sum_{i=1}^N\Bb_i\xb_i =\qb, ~~
\label{consensus problem
2}
\end{align}
where $N$ denotes the number of agents,
$f_i:\mathbb{R}^n \rightarrow \mathbb{R}$ are smooth and possibly non-convex
functions, and  $\Bb_i\in
\mathbb{R}^{M \times n}$ and $\qb\in \mathbb{R}^M$ are constants.
The target is to solve {\sf (P)} over a multi-agent network with
$N$ agents. It is assumed that each agent $i$ knows only $f_i$, $\Bb_i$ and $\qb$,
and can only communicate with its neighboring agents. The agent variables are
coupled due to the linear constraint $\sum_{i=1}^N\Bb_i\xb_i =\qb$.
Such problem arises, for example, in machine learning with distributed features
\cite{Chang14,Hu19},
power control in electrical power networks
\cite{molzahn2017survey,TsaiTSG2016},
interference management in wireless networks \cite{ShenTSP2012}, and the network
utility maximization problem \cite{palomar2007alternative}.

\vspace{-0.3cm}
\subsection{Related Works}

As a linearly constrained problem,  problem {\sf (P)}
can be easily solved \footnote{
Since finding a global minimum of problem {\sf (P)} is NP-hard in general,  what we mean ``solved" here is to reach a Karush-Kuhn-Tucker (KKT) solution of problem {\sf (P)}}. under the
centralized setting,
where all information are gathered in a central node.
Centralized methods for  solving linearly constrained problems are well studied in
the literature. For example, the penalty method \cite{Monteiro} penalizes the linear
constraints in the objective function and minimizes the
penalized objective function instead.
The augmented Lagrangian (AL) method \cite{BK:Bertsekas2003_NP} uses a min-max
approach to solve the resulting
AL function. Problem {\sf (P)} can also be efficiently handled by the alternating
direction method of multipliers (ADMM) \cite{hong2020block,wang2019global,Zhang18}.
However, solving  {\sf (P)} in a distributed fashion
over a
multi-agent network is challenging, due to the nonconvexity and coupled variables.
In the existing literature, most of the
distributed methods for solving {\sf (P)} assumed convexity of $f_i$'s.
For example, if {\sf (P)} is convex, it is well known that the dual decomposition
method \cite{palomar2006tutorial} is an efficient method to solve such problem over
a star network where all agents can communicate with a central node.
For a multi-agent network with general graph structures,  convexity and strong
duality are the keys to most of the existing distributed methods for solving {\sf
(P)}; see
\cite{Chang14,chang2016proximal,8472154,8695072,alghunaim2018dual,8709779,[8],Chang13}
and references therein.

Interestingly, when $\Bb =[\Bb_1,\cdots,\Bb_N]$ is the incidence matrix of the
network graph and $\qb$ is the zero vector, {\sf (P)} corresponds to a consensus
formulation of the finite sum minimization problem $\min_{\xb}\sum_{i=1}^N
f_i(\xb)$, and there exist rich results in the literature for both convex and
nonconvex cases,  see, for example,  the recent surveys
\cite{chang2020distributed,nedic2020distributed,
leinonen2019compressed}.
 Specifically,  it is worth pointing out that the recent work \cite{sun2019distributed} has shown that, for solving the non-convex finite sum problem in the distributed setting,  the first-order algorithms have an iteration complexity lower bound of $\Omega(1/\sqrt{\xi(\Gc)}\times 1/\epsilon)$, where $\xi(\Gc)$ is the spectral gap of the graph Laplacian matrix  and $\epsilon$ is the solution  accuracy to a stationary point.  The authors of \cite{sun2019distributed} also devised an optimal distributed algorithm that can attain the lower bound up to a factor.
However, since these decentralized methods are designed  for distributed finite sum
problems, they cannot be used to solve {\sf (P)} which has a general linear equality
constraint.

\vspace{-0.2cm}
\subsection{Contributions}
In this paper, we propose a new distributed algorithm for handling {\sf (P)} under
the general linearly coupled constraint and non-convex setting. The proposed method,
which we
call the \emph{proximal dual consensus} (PDC) method \cite{ZhangICASSP20}\footnote{
Compared with \cite{ZhangICASSP20}, this manuscript presents all the proof details and stronger theoretical results. Furthermore, the current manuscript presents an inexact update version of the proposed PDC algorithm as well as new experiment results.}, is
motivated by the proximal
technique recently proposed in \cite{Zhang18} and the dual consensus method in
\cite{Chang14, chang2016proximal}.
Specifically, a strongly convex counterpart of {\sf (P)} is obtained by introducing
a quadratic proximal term centered at an auxiliary primal variable.
Then, the dual consensus approach is applied to the proximal problem, resulting in a
distributed algorithm for {\sf (P)}.
In the proposed algorithm, the agents can choose to solve their local subproblem
exactly or inexactly by simply performing a gradient descent.

While the proposed PDC algorithms are seemingly a simple combination of \cite{Zhang18} and \cite{Chang14},  proving their convergence and convergence rate are far from easy as several new error bounds are needed.
Theoretically, we show that the proposed PDC algorithms can converge to a
Karush-Kuhn-Tucker (KKT) solution of {\sf (P)}, and generate an
$\epsilon$-optimal solution in $\mathcal{O}(1/\epsilon)$ iterations. This meets the iteration complexity lower bound  (up to a constant) of the first-order distributed algorithms for solving a nonconvex problem \cite{carmon2019lower,sun2019distributed}\footnote{
As mentioned,  since the finite sum problem is a special instance of problem {\sf (P)},  the iteration complexity lower bound
in \cite{sun2019distributed} also serves as a lower bound for {\sf (P)}.
}.
To evaluate the practical performance of the proposed algorithms,  we consider the less studied distributed vertical learning problem \cite{YingTSP2019,hu2019fdml}.
The presented numerical results
show that the proposed algorithms perform well
in a non-convex regularized logistic regression problem and a neural
network based distributed classification problem under the vertical learning setting. To the best of our knowledge, the proposed distributed algorithms are
the first to solve {\sf (P)} and the distributed vertical learning problem under the non-convex setting over the multi-agent network.

{\bf Synopsis:} Section \ref{section: application} presents two application examples
of {\sf (P)} about machine learning from  distributed features. The proposed PDC
algorithms are developed and analyzed in Section \ref{section: alg develop}. The
proof details of main theoretical results are shown in Section \ref{section: proof}.
Numerical results are given in Section \ref{section: simulation results} and
conclusions are drawn in Section \ref{section: conclusions}.

{\bf Notation: } $\mI_n$ is denoted as the $n$ by $n$ identity matrix.
$\langle \va, \bb \rangle$ represents the inner product of vectors $\va$ and $\bb$,
$\|\va\|$ is the Euclidean norm, and $\|\ab\|_{\mA}^2$ represents $\ab ^\top
{\mA}\ab$; for a matrix $\mA$,  both $\|\mA\|$ and $\sigma_{\max}(\mA)$ denote its matrix 2-norm (maximum singular value), and $\sigma_{\min}(\mA)$ denotes the minimum non-zero singular value of  $\mA$;
$\otimes$ denotes the Kronecker product;
 $i\in [N]$ denotes $i \in \{ 1, \ldots, N\}$.

\vspace{-0.0cm}
\section{Application Examples} \label{section: application}

\begin{figure}[!t]
  \centering
    \includegraphics[width=4.75in]{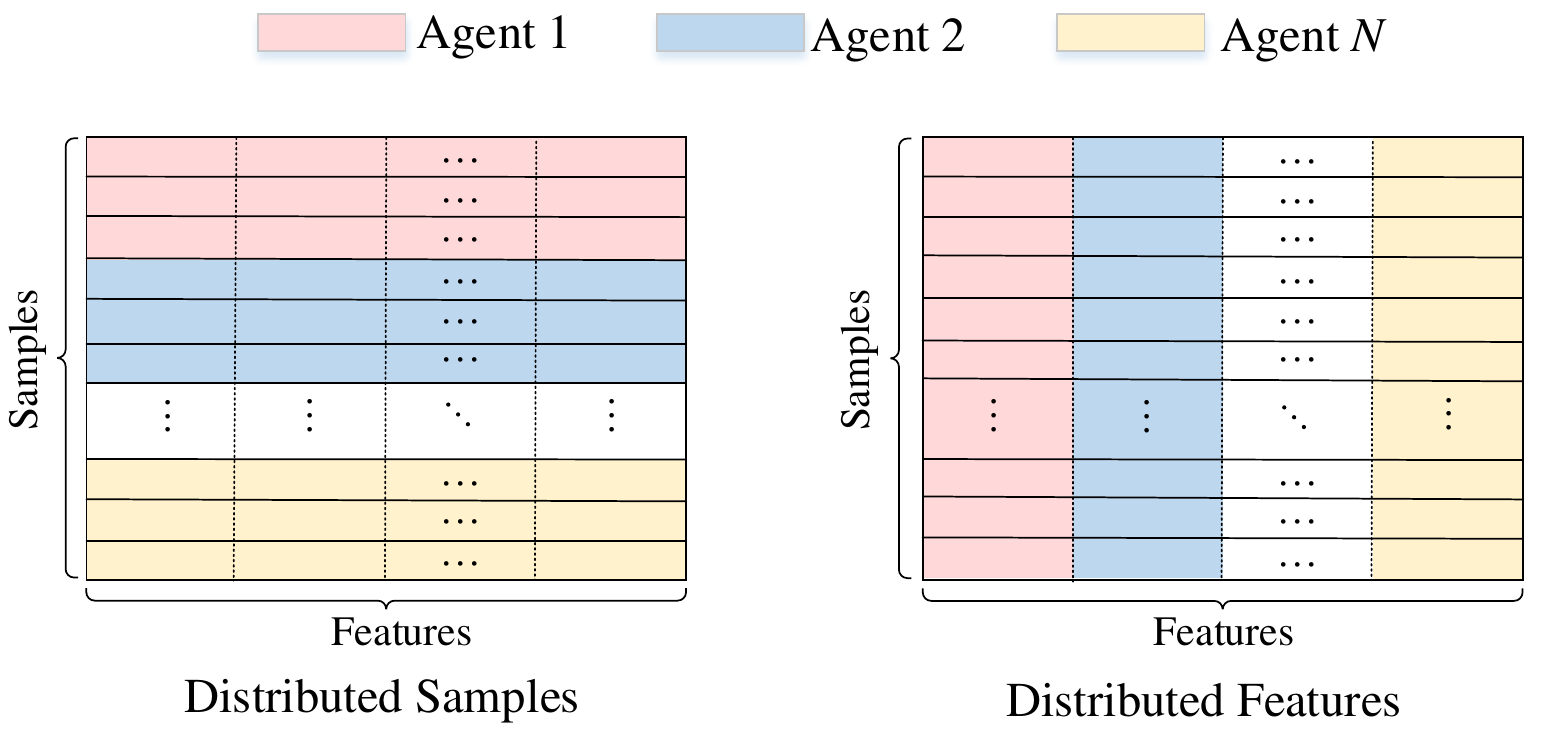}\\
  \caption{ Data partition for horizontal learning (Left) and
  vertical learning (Right). }\label{fig: Distributedfeatures}
  \vspace{-0.5cm}
\end{figure}

Problem {\sf (P)} appears in many engineering applications. As mentioned in Section \ref{sec:intro}, one classical example is the consensus formulation of the finite sum problem $\min_{\xb}\sum_{i=1}^N
f_i(\xb)$.  Specifically,  consider a multi-agent network with $N$ connected agents, and each agent $i$, $i\in [N]$, has a local copy of $\xb$,  denoted as $\xb_i \in \mathbb{R}^n$.  Suppose that $\tilde \Ab\in \mathbb{R}^{N\times N}$ is the incidence matrix of the network graph.
One can have the following consensus formulation of the finite sum problem
\begin{subequations}\label{eqn: finite sum}
\begin{align}
\min_{\xb =[\xb_1^\top, \cdots, \xb_N^\top]^\top}&~\sum_{i=1}^N f_i(\xb_i)\\
\text{s.t.}~&(\tilde \Ab \otimes \Ib_n) \xb = \mathbf{0},\label{eqn: consensus}
\end{align}
\end{subequations}
where \eqref{eqn: consensus} ensures all $\xb_i$'s are the same to each other. In the machine learning literature,  the finite sum problem corresponds to the conventional \emph{horizontal learning} problem where the data samples are partitioned and distributively owned
by the agents \cite{chang2020distributed} (left of Fig. \ref{fig: Distributedfeatures}),  and such problem has been extensively studied \cite{chang2020distributed,nedic2020distributed,
leinonen2019compressed}.
Unlike \eqref{eqn: finite sum},  in this section we discuss an important but less discussed learning problem: the \emph{vertical learning} problem,  where the agents own all data samples
but part of their features only (right of Fig. \ref{fig: Distributedfeatures}) \cite{YingTSP2019,hu2019fdml}. Such problems appear, for example, in collaborative
learning between multiple parties who hold different aspects of information for the
same set of samples.
Using two examples below, we show that the vertical learning problem can be expressed as the form
of {\sf (P)}.

%
%
%
%
%

{\bf Example 1} {(\bf  Empirical risk minimization)}
Consider a standard empirical risk minimization problem as follows
\begin{equation}\label{distributed-feature-general}
\min_{ \wb  }~ \sum\limits_{k=1}^M\psi (  \bb_k^\top  \wb, \vb_k  )+ R( \wb  ),
\end{equation}
where $M$ is the number of data samples, $\bb_k$ and $\vb_k$ are respectively the
feature vector and (one-hot vector) label
for the $k$-th sample,
$ \wb  \in \mathbb{R}^{nN}$ is the parameter to be optimized,
$\psi$ is the (possible non-convex) loss function, and $R(\cdot)$ is  a regularization function.

 Let us consider to solve problem \eqref{distributed-feature-general} over a multi-agent
network with $N$ agents under the vertical learning setting.  Each agent $i$ owns all $M$ data samples but has part of the data features only.  Specifically,
assume that each sample $\bb_k$ is partitioned as $\bb_k = [\bb_{1, k}^\top, \cdots,
\bb_{N, k}^\top]^\top$, where $\bb_{i,k}\in \mathbb{R}^{n}$, and each agent $i$ owns
samples with only partial features $\{\bb_{i, k}\}_{k=1}^M$, for all $i\in [N]$.
Let the parameter vector
$\wb$ be partitioned in the same fashion, i.e.,
$\wb = [\wb_1^\top, \cdots, \wb_N^\top]^\top$, and assume the regularization term
$R(\wb) = \sum_{i=1}^N R_i(\wb_i)$. So,  each agent $i$ under the vertical learning setting maintains only a subset of the model parameters $\wb_i$.

Then, problem \eqref{distributed-feature-general} can be reformulated as
\begin{subequations}\label{eq: partition distributed-feature-general}
\begin{align}
\min_{{\wb},\wb_0 } ~&\sum\limits_{k=1}^M\psi \big( w_{0,k}, \vb_k \big
)+\sum_{i=1}^N R_i(\wb_i) \\
{\rm s.t.}~& \sum_{i=1}^N \bb_{i, k}^\top \wb_i-w_{0,k}=0,~k\in[M],
\end{align}
\end{subequations}
where $\wb_0=[w_{0,1},\ldots,w_{0,M}]^\top$ is an introduced variable.
Further define ${\Bb}_i=[\bb_{i,1}, \cdots, \bb_{i,M}]^\top$, $i\in[N]$, and let
$\Psi(\wb_{0})=\sum_{k=1}^M\psi(w_{0, k},
\vb_k)$.
Thus  \eqref{eq: partition distributed-feature-general} can be written as
\begin{subequations}\label{consensus problem 22}
\begin{align}
~~ \min_{\substack{\wb_1, \ldots, \wb_N,\\ \wb_0}}~
&\Psi(\wb_{0})+\sum_{i=1}^NR_i(\wb_i) \\
~~\mathrm{ s.t. }~ &\sum_{i=1}^{N}\Bb_i\wb_i - \wb_0 =\zerob.
\end{align}
\end{subequations}
One can see that the above problem \eqref{consensus problem 22} is indeed  an instance of
{\sf (P)}.  It corresponds to the vertical learning problem where each agent $i$ knows partial feature matrix $\Bb_i$,  regularizer $R_i$ and partial model parameter $\wb_i$, for all $i\in [N]$; the variable $\wb_0$ and $\Psi$ can be handled by any one of the agents or an extra auxiliary agent.

\begin{figure}[!t]
  \centering
  \includegraphics[width=5in]{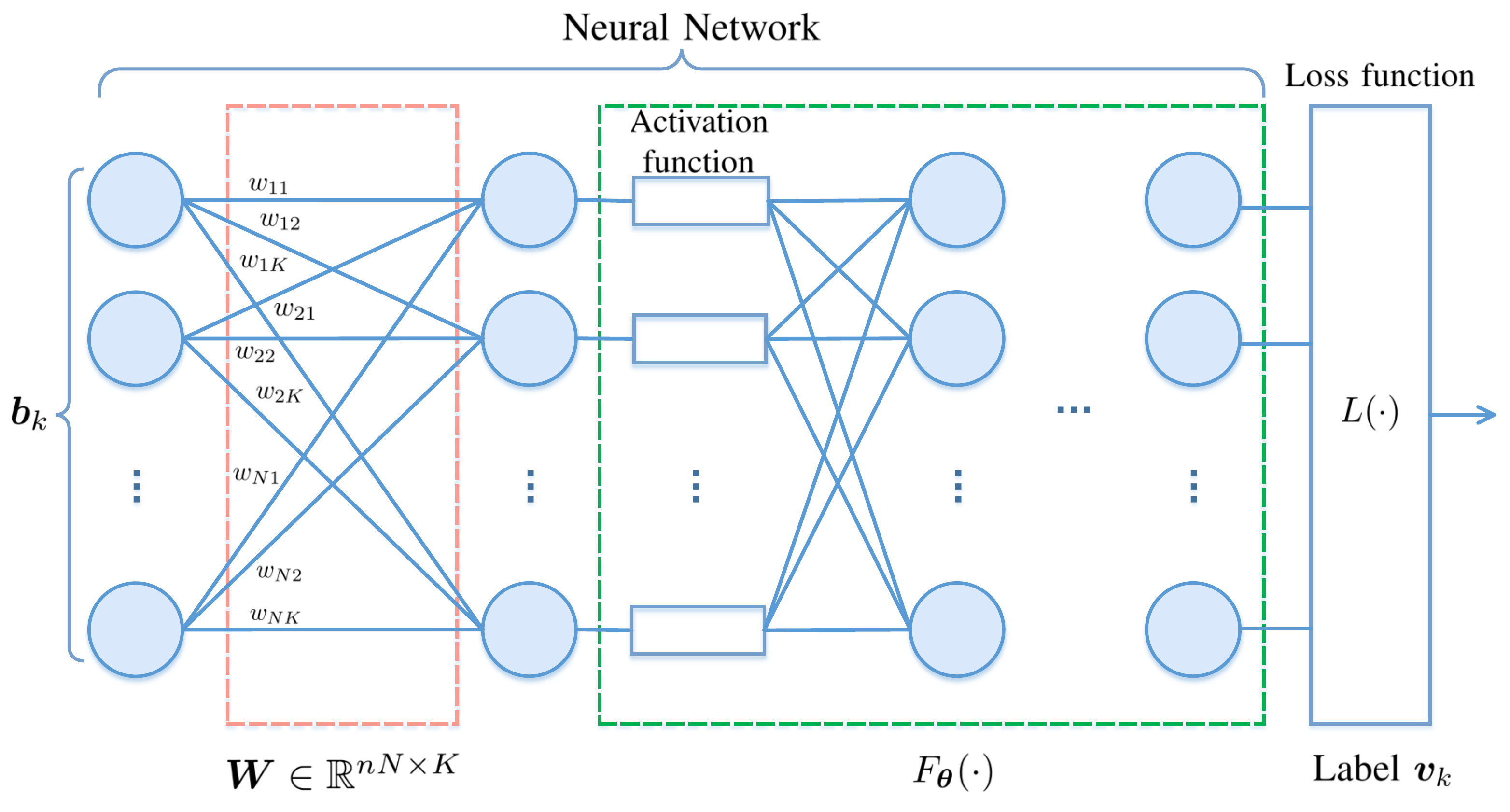}\\
  \caption{{A neural network structure used in problem \eqref{NN
  problem}.}}\label{DNN}\vspace{-0.2cm}
\end{figure}

{\bf Example 2} {(\bf  Classification by neural network)}
Consider to train a neural network (NN) for classification as shown in Fig.
\ref{DNN}.
The mapping of the
NN is expressed by the function $F_{\thetab}(\bb_k^\top\Wb)$, where $\Wb\in
\mathbb{R}^{nN \times K}$ denotes the coefficients of the first linear layer of the
NN, $K$ is the dimension of the 1st layer output, and $\thetab$ incorporates the
parameters from the 1st layer output to the output layer of the NN.
The learning problem is given by
\begin{align}\label{NN problem}
\min_{\substack{\Wb,\thetab}}~\sum\limits_{k=1}^M L\big(F_{\thetab}(\bb_k^\top\Wb),
\vb_k \big),
\end{align}
where $L$ is a loss function. By comparing \eqref{NN problem} with \eqref{distributed-feature-general},  one can see that the previous empirical risk minimization problem is equivalent to the NN problem with only a single linear layer,  $K=1$ and no $\thetab$.

Similarly,  let us consider the vertical learning setting.  Suppose that both the data samples and model parameters are partitioned as $\bb_k = [\bb_{1, k}^\top, \cdots, \bb_{N, k}^\top]^\top$ for all $k\in [M]$,  and $\Wb=[\Wb_1^\top, \ldots, \Wb_N^\top]^\top$.  By introducing
\begin{align}
\wb_{0,k} = \sum_{i=1}^N \Wb_i^\top \bb_{i,k}=\sum_{i=1}^N (\bb_{i,k}^\top \otimes
\Ib_{K})\wb_i,
\end{align}
where $\wb_i={\rm vec}(\Wb_i)$ is obtained by stacking the columns of
$\Wb_i $,  problem \eqref{NN problem} can be written as
\begin{subequations}\label{eq: distributed NN problem }
\begin{align}
\min_{\substack{\wb_1, \cdots, \wb_N, \\ (\wb_0,\thetab)}}~ &\sum\limits_{k=1}^M
L\big(F_\thetab( \wb_{0,k}^\top), \vb_k \big) \\
{\rm s.t. }~& \sum_{i=1}^N \Bb_i \wb_i - \wb_0 =\zerob,
\end{align}
\end{subequations}
where $\textstyle\Bb_i = \textstyle [(\bb_{i,1}^\top\otimes \Ib_K)^\top , \cdots,(\bb_{i,M}^\top\otimes
\Ib_K)^\top  ]^\top$ and $\wb_0=[\wb_{0,1}^\top,\ldots,\wb_{0,M}^\top]^\top$.

As a result, the distributed optimization problem {\sf (P)} includes both problems \eqref{consensus problem 22} and \eqref{eq: distributed NN
problem } as special cases, and thereby any distributed algorithm that can solve {\sf (P)} can be  applied for solving the above two vertical learning problems.  We should mention that the existing works for solving the vertical learning problem in the distributed networks consider only convex problems or over a simple network \cite{YingTSP2019,hu2019fdml}.  In the next section,  we present a new distributed algorithm for non-convex problem {\sf (P)} over general connected networks.

\vspace{-0.1cm}
\section{Proposed PDC Algorithm} \label{section: alg develop}
The network models and assumptions are given in Section \ref{sec: network models and
assumptions}. To solve {\sf (P)}, we develop the PDC algorithm in
Section \ref{sec: PDC alg} and further give its convergence analysis in Section
\ref{sec: PDC thm}. In Section \ref{sec: IPDC alg}, we further propose an inexact
version of PDC, referred as the IPDC.

\vspace{-0.2cm}
\subsection{Network Model and Assumptions }\label{sec: network models and
assumptions}
The multi-agent network is modeled as an undirected  graph $\Gc=(\Ec,\Nc)$, where
$\Ec$ is the set of edges and $\Nc$ is the set of agents.  Each agent can only
communicate and exchange information with its neighbors, i.e., agents $i, j$ can
communicate with each other if and only if $(i, j)\in \Ec$.

Let us define some useful notations regarding the network graph.
The neighbouring set of agent $i$ is defined as $\Nc_i=\{j\in \Nc~|~
    (i,j)\in \Ec\}$. Let $\tilde \Db ={\rm diag}\{|\Nc_1|,\cdots,|\Nc_N|\} \in
    \mathbb{R}^{N \times N}$ be the degree matrix of $\Gc$.
The agent-edge incidence matrix is defined as $\tilde \Ab  \in \mathbb{R}^{|\Ec|
\times |\Nc|}$ which has
    $\tilde \Ab _{e,i} =1$ and $\tilde \Ab _{e,j}=-1$
for $j>i$ and $e=(i,j)\in \Ec$, and zeros otherwise.
$\tilde \Lb^- = \tilde \Ab^\top \tilde \Ab$ is the (signed) Laplacian matrix
    of $\Gc$	. Then
$\tilde \Lb^+= 2\tilde \Db -\tilde \Lb^-$
is the signless Laplacian matrix of $\Gc$. Both $\tilde \Lb^-$ and $\tilde
\Lb^+$ are positive semi-definite and describe the connectivity of the network.
For ease of later use, we also define
\begin{align*}
&\Ab=\tilde \Ab \otimes \Ib_M=[\Ab_1,\ldots, \Ab_N],\\
&\Lb^+=\tilde \Lb^+ \otimes \Ib_M =\begin{bmatrix}
(\Lb_1^+)^\top,\cdots,
(\Lb_N^+)^\top
\end{bmatrix}^\top,\\
&\Lb^- =\tilde \Lb^- \otimes \Ib_M
=\begin{bmatrix}
(\Lb_1^-)^\top,\cdots,
(\Lb_N^-)^\top
\end{bmatrix}^\top,
\end{align*}
and let $\lambda_{\max}$ be the maximum eigenvalue
of $\Lb^+$.

We make the following
assumption for the network.
\vspace{-0.2cm}
\begin{Assumption}\label{assumption connected G}
The graph $\Gc$ is connected.
\end{Assumption}
It means that for any $i, j\in \Nc$, there exists at least one path connecting $i$
and $j$. Assumption \ref{assumption connected G} is commonly made in the literature
of distributed optimization.

For problem {\sf (P)}, we make the following assumption.
\vspace{-0.2cm}
\begin{Assumption}\label{assumption smooth}
Each $f_i$ is continuously differentiable and $\nabla f_i$ is Lipschitz continuous
with constant $\gamma^+>0$.
\end{Assumption}
Under Assumption \ref{assumption smooth}, there exists
 a (possibly negative) constant  $\gamma^-$  such that
\begin{equation}
\la \nabla f_i(\xb_i) - \nabla f_i(\xb_i'), \xb_i -\xb_i' \ra \geq \gamma^-
\|\xb_i -\xb_i'\|^2,
\end{equation}  for all $\xb_i,\xb_i' \in \mathbb{R}^n$.
Note that if $ \gamma^-$ is negative, then $f_i$ is not a convex function.

\subsection{Algorithm Development }\label{sec: PDC alg}
The proposed PDC method is motivated by the proximal approach recently presented in
\cite{Zhang18} for non-convex constrained optimization.
In particular, like  \cite{Zhang18}, we consider a ``proximal" counterpart of {\sf
(P)}
\begin{subequations}\label{consensus problem 2 PX}
\begin{align}
\min_{ \xb_i, i\in [N]}~ &\sum_{i=1}^Nf_i(\xb_i) +\frac{p }{2} \sum_{i=1}^N\|\xb_i-\zb_i\|^2
\\
\text{s.t.}~ & \sum_{i=1}^N\Bb_i\xb_i =\qb,\label{consensus problem 2 PX C1}
\end{align}
\end{subequations}
where $p>-\gamma^-$ is a penalty parameter and
$\zb=(\zb_1^\top,\ldots,\zb_N^\top)^\top$ is
an introduced proximal variable.
 Note that  \eqref{consensus problem 2 PX}
is still not a convex problem for the joint variable $(\xb^\top,\zb^\top)^\top$,  but it is convex with respect to $\xb$ when $\zb$ is fixed.  So,  the proximal problem \eqref{consensus problem 2 PX}  can be regarded as a convex approximation of problem  {\sf (P)} around the given point $\zb$.

Let us denote $(\xb(\zb),\yb_0(\zb))$ as a primal-dual optimal solution of \eqref{consensus
problem 2 PX} with fixed $\zb$. One can verify that $(\xb(\zb),\yb_0(\zb))$
is a KKT solution of  {\sf (P)}, i.e., satisfying
\begin{subequations}\label{eqn: KKT of P}
     \begin{align}
        &  \nabla f_i (\xb_i(\zb_i)) + \Bb_i^\top\yb_0(\zb)  = \mathbf{0} , \forall i\in [N] \\
        &\sum_{i=1}^N \Bb_i\xb_i(\zb_i) - \qb = \mathbf{0},
    \end{align}
    \end{subequations}
if and only if
\begin{align}\label{eqn: kkt}
   \xb(\zb) =\zb.
\end{align}
In \cite{Zhang18}, the authors proposed to handle \eqref{consensus problem 2 PX} by
using one step of primal-dual gradient update in each iteration $r$, followed by
updating $\{\zb^r\}$ by
\begin{align}
\zb^{r+1}=\zb^r+\beta(\xb^{r+1}-\zb^r),\label{eqn: alg original form z}
\end{align}
where $0<\beta<1$ is a positive stepsize. Notice that for $\beta = 1$, the algorithm may oscillate.
It is shown that \eqref{eqn: alg original form z} with sufficiently small $\beta$ can stabilize the primal-dual
sequence and the algorithm can achieve a KKT solution of {\sf (P)}.
However, the algorithm in \cite{Zhang18} is a centralized algorithm and cannot solve
{\sf (P)} over the multi-agent network.

To obtain a distributed algorithm,
we consider employing the dual consensus method in  \cite{Chang14}
to deal with the (strongly convex) proximal problem \eqref{consensus problem 2 PX} when $\zb$ is fixed.
First, by the Lagrange duality theory \cite{BK:BoydV04},  \eqref{consensus problem 2
PX} is equivalent to its dual problem  when $\zb$ is fixed
\begin{align}\label{eq: dual of consensus problem 2}
 \max_{\yb_0 \in \mathbb{R}^M} ~\sum_{i=1}^N \bigg(\phi_i(\yb_0, \zb_i) -
\yb_0^\top \frac{\qb }{N}\bigg),
\end{align}
where $\yb_0$ is the dual variable associated with the constraint \eqref{consensus
problem 2 PX C1}, and
\begin{equation}
\begin{aligned}
\phi_i(\yb_0,\zb_i) = \min_{\xb_i\in \mathbb{R}^n} f_i(\xb_i) +
\frac{p}{2}\|\xb_i-\zb_i\|^2 + \yb^\top_0\Bb_i\xb_i,
\end{aligned}
\end{equation} for $i\in [N],$
are concave dual functions of $\yb_0$.
In order to solve this dual problem in a decentralized manner, we let $\yb_i\in
\mathbb{R}^M$ be the local copy of $\yb_0$ at agent $i$. Then, $\phi_i$ can be
written as
\begin{equation}
\begin{aligned}\label{eq: psi}
\phi_i(\yb_i,\zb_i) = \min_{\xb_i\in \mathbb{R}^n} f_i(\xb_i) +
\frac{p}{2}\|\xb_i-\zb_i\|^2 + \yb^\top_i\Bb_i\xb_i,
\end{aligned}
\end{equation} for $i\in [N]$.
As the right hand side of \eqref{eq: psi} is a strongly convex problem with respect to $\xb_i$, $\phi_i(\yb_i,\zb_i)$ is smooth
and concave with respect to $\yb_i$ according to the Danskin's theorem \cite{bertsekas1997nonlinear}, and $\nabla \phi_i(\yb_i,\zb_i)$ is Lipschitz continuous
with constant $\frac{1}{p+\gamma^-}$ \cite{Nesterov05}.
Next,  similar to \eqref{eqn: finite sum}, we consider the following consensus
constrained counterpart
\begin{subequations}\label{eq: dual of consensus problem 2 eq}
	\begin{align}
	\max_{\yb=(\yb_1^\top, \cdots, \yb_N^\top)^\top} ~&\sum_{i=1}^N
\bigg(\phi_i(\yb_i,\zb_i) - \yb_i^\top\frac{\qb }{N}\bigg) \\
	{\rm s.t.}~& \Ab\yb=\zerob. \label{eq: dual of consensus problem 2 eq C1}
	\end{align}
\end{subequations}
 Note that under Assumption \ref{assumption connected G},
$\Ab\yb=\zerob$ is equivalent to $\yb_1=\yb_2 =\ldots=\yb_N$.

We further employ the inexact augmented Lagrangian (AL) method in \cite{Hong17}  to
problem \eqref{eq: dual of consensus problem 2 eq}. In particular, let
\begin{align} \label{eqn: AL of DP2}
&\Lc_\rho(\yb,\mub;\zb)
\notag \\
&=\sum_{i=1}^N \bigg(\phi_i(\yb_i,\zb_i) -
\yb_i^\top\frac{\qb}{N}\bigg)-\mub^\top\Ab \yb  -\frac{\rho}{2}\| \Ab \yb \|^2,
\end{align}
be the augmented Lagrangian function of \eqref{eq: dual of consensus problem 2 eq},
in which $\mub \in \mathbb{R}^{M|\Ec|}$ is the dual variable associated with
\eqref{eq: dual of consensus problem 2 eq C1}, and $\rho>0$ is a penalty parameter.
Then, according to  \cite{Hong17}, at each iteration $r$, we perform
\begin{subequations}\label{eqn: alg original form}
	\begin{align}
	&\mub^{r+1} = \mub^r + \alpha  \Ab \yb^r, \label{eqn: alg original form mu} \\
	& \yb^{r+1} = \arg\max_{\yb} \Lc_\rho(\yb,\mub^{r+1};\zb^r) -
\frac{\rho}{2}\|\yb - \yb^r\|_{\Lb^+}^2,\label{eqn: alg original form y}
	\end{align}
\end{subequations}
where $\alpha>0$ is a step size parameter.
The proximal term $\frac{\rho}{2}\|\yb - \yb^r\|_{\Lb^+}^2$ is used to make the
objective function of \eqref{eqn: alg original form y} decomposable across $\yb_1,
\cdots, \yb_N$ so as to obtain a fully distributed algorithm, as we delineate
below.

\vspace{0.2cm}
\noindent\underline{\bf  Updates of $\xb^r$ and $\yb^r$}:
According to \eqref{eq: psi} and \eqref{eqn: AL of DP2}, we define
\begin{align}\label{eq: potential function G}
 G(\xb, \yb, \mub; \zb)
&\triangleq
\sum_{i=1}^N \bigg( f_i(\xb_i) + \frac{p}{2}\|\xb_i-\zb_i\|^2 +
\yb_i^T\Bb_i\xb_i\bigg. \notag \\
\bigg.&~~~ - \yb_i^\top \frac{\qb}{N}\bigg)
-\mub^\top\Ab \yb -\frac{\rho}{2}\| \Ab \yb \|^2.
\end{align}
Then, by \eqref{eqn: alg original form y}, $(\yb^{r+1}, \xb^{r+1})$ is the solution
of the following max-min
problem:
\begin{align}\label{eqn: maxmin1}
\max_{\yb}\min_{\xb\in \mathbb{R}^{nN}}\bigg\{G(\xb, \yb, \mub^{r+1};
\zb^r)-\frac{\rho}{2}\|\yb-\yb^r\|^2_{\Lb^+}\bigg\}.
\end{align}

\noindent Similar to \cite{Chang14}, based on the
saddle point theory, the above max-min problem  is equivalent to the following
min-max problem
\begin{align}\label{eqn: minmax}
\min_{\xb\in\mathbb{R}^{nN}}\max_{\yb} \bigg\{G(\xb, \yb, \mub^{r+1};
\zb^r)-\frac{\rho}{2}\|\yb-\yb^r\|^2_{\Lb^+}\bigg\}.
\end{align}

When fixing $\xb$, the inner maximization problem for $\yb$ is a strongly concave
quadratic
problem. Then, one can show that the inner solution $\yb^{r+1}$
are decomposable and has a closed-form expression as
\begin{align}  \label{y: update}
\yb_i^{r+1} &=   \frac{1}{2\rho |\Nc_i|} \big( \Bb_i \xb_i -  \frac{\qb}{N}
-\Ab^\top_i\mub^{r+1}+ \rho  \Lb^+_i \yb^r \big),
\end{align}
for all $i \in [N]$.
By substituting \eqref{y: update} into \eqref{eqn: minmax},
the outer solution $\xb^{r+1}$ of \eqref{eqn: minmax} can be determined by
\begin{align} \label{eq: x update}
&\xb_i^{r+1} =\arg \min_{\xb_i\in \mathbb{R}^n}  \bigg\{ f_i(\xb_i) +
\frac{p}{2}\|\xb_i-\zb_i^r\|^2 \notag \\
&~~~+ \frac{1}{4\rho |\Nc_i|}
\|\Bb_i\xb_i - \qb/N - \Ab_i^\top \mub^{r+1} + \rho  \Lb_i^+  \yb^r\|^2\bigg\},
\end{align} for all $i \in [N]$, which are convex problems given $p>-\gamma^-$.
Notice that
\begin{align}
 \Lb_i^+  \yb^r\ =\sum_{j\in \mathcal{N}_i} (\yb_i^r+\yb_j^r),
\end{align} and therefore both \eqref{y: update} and \eqref{eq: x update} can be
implemented fully distributively using only local data and messages from the
neighbors.

\vspace{0.2cm}
\noindent\underline{{\bf Update of $\mub^r$}}:
Since in \eqref{y: update} and \eqref{eq: x update},  one needs $\Ab_i^\top \mub$
instead of $\mub$. We define $\pb_i^r=\Ab_i^\top\mub^r$, and according to
\eqref{eqn: alg original form mu}, obtain
\begin{align}\label{eqn: update of p99}
\pb^{r+1}_i&=\pb_i^{r}+ \alpha \Ab_i^\top \Ab \yb^r \notag\\
&=\pb_i^{r}+ \alpha \Lb_i^- \yb^r \notag\\
&=\pb_{i}^{r}+{\alpha}\textstyle \sum_{j \in \Nc_i}(\yb_i^{r}-\yb_j^{r}),	
\end{align} for $i \in [N]$.

By combining \eqref{y: update}, \eqref{eq: x update}, \eqref{eqn: update of p99},
and \eqref{eqn: alg original form z}, we obtain the PDC algorithm (see Algorithm
\ref{table: PDC}).
One can see that the algorithm is fully distributed since each agent $i$ uses only local information and
variables $\{\yb_j\}_{j\in \Nc_i}$ from its neighbors.

\subsection{Convergence Analysis}\label{sec: PDC thm}
We first define the approximate solution of problem {\sf (P)} as follows.
\begin{Definition}
 A vector $\xb \in\mathbb{R}^{Nn}$ is said to be an $\epsilon$-KKT solution of problem ${\sf {(P)}}$ if $\|\sum_{i=1}^N \Bb_i\xb_i -\qb\|^2 \le \epsilon$ and
 there exists $\yb\in
    \mathbb{R}^M$ such that $\sum_{i=1}^N\| \nabla f_i(\xb_i)+ \Bb_i^\top \yb \|^2 \leq \epsilon$.
\end{Definition}

\vspace{-0.0cm}
The main theoretical results for the proposed PDC method in Algorithm \ref{table:
PDC} are given in the following theorem.
\begin{Theorem}\label{thm: conv and iteration comp}
Let $p>-\gamma^-$, $\rho>0$.
Then
for sufficiently small $\alpha<\rho$ and $\beta<1$ (see \eqref{eqn: alpha condition}
and \eqref{eqn: beta condition}, the following results hold for Algorithm
\ref{table:
PDC}:
	\begin{enumerate}[(a)]
		\item Every limit point
		of the iteration sequence $\{(\zb^r, \yb^r)\}$ is a KKT solution of
problem
		 {\sf (P)};
		\item There exists a $t<r$ such that $\zb^t$ is an $\mathcal{O}(1/r)$-KKT solution of
problem {\sf (P)}.
	\end{enumerate}
\end{Theorem}

\begin{algorithm}[!t]
	\caption{Proposed PDC method for solving {\sf (P)}}
	\begin{algorithmic}[1]\label{table: PDC}
		\STATE {\bf Given} parameters $p, \beta, \alpha,
\rho$, and initial variables
		$\xb_i^{0}=\zb_i^0 \in \mathbb{R}^{n}$, $\yb_i^{0}\in \mathbb{R}^{M}$  and
$\pb_{i}^{0}=\zerob$ for $i\in [N]$. Set $r=1.$
		\REPEAT
		\STATE  For all  $i\in [N]$ (in parallel),
send $\yb_i^r $ to neighbour $j \in \Nc_i$ and receive $\yb_j^r$ from neighbour
$j\in \Nc_i$. Perform  \begin{align}
		&\pb_{i}^{r+1}= \textstyle ~ \pb_{i}^{r}+{\alpha}\sum_{j \in
\Nc_i}(\yb_i^{r}-\yb_j^{r}), 			\label{eq: dual ADMM p}\\
				\label{eq: dual ADMM x1}
		& \xb_i^{r+1}=\arg\min_{\xb_i }\textstyle
\bigg\{ f_i(\xb_i)+\frac{p}{2}\|\xb_i-\zb_i^r\|^2
\notag \\
		&~~~~~~~~+ \textstyle \frac{1}{4\rho|\mathcal{N}_i|}\!
		\big\|\Bb_i\xb_i\!
		  -\!\frac{\qb}{N}
		\textstyle \!-\! \pb_i^{r+1}\! +\! \rho \Lb^+_i \yb^r
		\big\|^2 \bigg\},
		\\
		\label{eq: dual ADMM lambda 21}
	& \textstyle \yb_i^{r+1} \!= \!  \frac{1}{2\rho |\Nc_i|}\big( \Bb_i
\xb_i^{r+1}\! -\! \frac{\qb}{N} \!-\!\pb_i^{r+1}\!  + \rho \Lb^+_i \yb^r\big),
	\\
	& \textstyle \zb_i^{r+1} = \zb_i^r +\beta (\xb_i^{r+1} - \zb_i^r),
	\label{eq: dual ADMM lambda 22}
		\end{align}
		\STATE {\bf Set} $r\leftarrow r+1.$
		\UNTIL {a predefined stopping criterion is satisfied.}
	\end{algorithmic}
\end{algorithm}

{\bf Proof:}
The proof relies on analyzing how the gap between $\zb^t$ and $\xb(\zb^t)$ progresses with the iterations,  and relating it with the gap between $\zb^t$ and a KKT solution of {\sf (P)} (Lemma \ref{lemma9}).  The details are given
in Section \ref{section: proof}. \hfill $\blacksquare$

Theorem \ref{thm: conv and iteration comp} implies that, given constants $\rho$ and
$p>-\gamma^-$, if step size parameters $\alpha$ and $\beta$ are sufficiently small,
$\{(\zb^r, \yb^r)\}$ can
asymptotically converge to the set of KKT solutions of problem {\sf (P)}. Moreover, Theorem \ref{thm: conv and iteration comp}(b) shows that an
$\epsilon$-KKT solution can be obtained in $\mathcal{O}( {1}/{\epsilon})$
iterations.
As mentioned in Section \ref{sec:intro},  since the finite sum problem is a special instance of problem {\sf (P)},  the iteration complexity lower bound
in \cite{sun2019distributed} also serves as a lower bound for {\sf (P)}.
Therefore,  Theorem \ref{thm: conv and iteration comp}(b) means that the proposed PDC algorithm can attain the  iteration complexity lower bound of first-order distributed algorithms up to a constant.

While Theorem \ref{thm: conv and iteration comp}(b) gives only the complexity order,  the following corollary presents the dependency
of the bound on the graph topology and data matrix.

 \begin{Corollary} \label{lem: convergence graph}
      Assume that $\Bb=[\Bb_1,\ldots,\Bb_N]$ is full row rank.    By choosing appropriate values for $\alpha, \beta$ and $\rho$, the convergence gap $\|\zb^t-\xb(\zb^t)\|^2$ in \eqref{eq: bound for z x(z) 1} can be further bounded as

        \vspace{-0.3cm}
        {\small
        \begin{align}\label{eqn: further bound}
           \textstyle & \|\zb^t-\xb(\zb^t)\|^2 \notag\\& \leq
            \textstyle \frac{\Phi^0-\underline{\Phi}}{r}\textstyle \left(
              D_1 B_{\max}^2 \lambda_{\max} \sqrt{D_2  \theta_1^4 +D_{3} \theta_1^2  } \right.\notag \\
            &~~~\left.
             +
             \sqrt{D_4 + \textstyle  B_{\max}^2 \lambda_{\max}^2 (D_5 \theta_1^4 + D_6\theta_1^2)
            \textstyle+ \textstyle \frac{D_7 + D_8 B_{\max}^2\theta_2^2}{B_{\max}^2\lambda_{\max}^2(D_9 \theta_1^4 + D_{10}\theta_1^2) }} \right.\notag \\
            &~~~\left.\textstyle
                +\textstyle D_{11} B_{\max}\theta_1 \sqrt{ D_{12} \theta_3^2 + D_{13}  }
            \times \max\bigg\{D_{14}, \sqrt{D_{15} \lambda_{\max}^2 \theta_1^2 + D_{16}}\bigg\}\right)^2,
        \end{align} }

                \vspace{-0.3cm}
\noindent        where $\theta_1$,  $\theta_2$ and $\theta_3$,  respectively defined in Lemma \ref{fact: perturb of y 2},  Lemma \ref{fact: perturb of y} and Lemma \ref{Fact y(z^r) - y(mu^{r+1, z^r})},  can be bounded as
        \begin{align}
            & \textstyle \theta_1=  \theta_2  \leq\textstyle
             (1+ \zeta_B)^2 \frac{\sigma_{\max}(\Lb^-)}{\sigma_{\min}^2(\Lb^-) }
             + \frac{4(1+ \zeta_B)^2}{ \sigma_{\max}(\Bb_{\mathrm{diag}})}  , \label{eq: theta_1}\\
    &\textstyle \theta _3 \leq
            \textstyle \frac{\left(2\sqrt{\sigma_{\max}(\Lb^-)}+\sigma_{\max}(\Lb^-)\right)(3+\sqrt{\sigma_{\max}(\Lb^-)} )^2 (1+ \zeta_B)^2 }
            {\sigma_{\min}(\Lb^-) }
            + \textstyle
                \frac{ 4(1+\zeta_B)^2}{\sigma_{\max}(\Bb_{\mathrm{diag}})}\left(3+\sqrt{\sigma_{\max}(\Lb^-)} \right)^2 , \label{eq: theta_3}\\
                & \textstyle B_{\max}  \textstyle \triangleq\textstyle \max_{i\in [N]} \|\Bb_i\| ,~~~\textstyle\zeta_B \triangleq \textstyle\frac{ 2\sqrt{N}\cdot \sigma_{\max}(\Bb_{\mathrm{diag}})}{\sigma_{\min}(\Bb)},
        \end{align}
$\lambda_{\max}$ is the largest eigenvalue of $\Lb^+$,  $\Bb_{\mathrm{diag}}$ is a block diagonal matrix with $\Bb_1,\ldots,\Bb_N$ being the diagonal blocks, and $D_1$ to $D_{16}$ are constants which depend only on $\gamma^-$ and $\gamma^+$.
 \end{Corollary}

 {\bf Proof:}  The bound \eqref{eqn: further bound} is obtained by further analyzing the dependency of the three parameters $\theta_1$,  $\theta_2,$ $\theta_3$ on the graph topology and data matrix. The details can be found in Appendix \ref{app: corollary 1}. \hfill $\blacksquare$

     As seen from Corollary \ref{lem: convergence graph},  both the graph topology ($\Lb^-$ and $\Lb^+$) and data matrix $\Bb$ can have impacts on the algorithm convergence.  However,  it is difficult to give a quantitative conclusion due to its intricacy,  except for some simple graphs.
The following corollary proved in \cite{zhang2021decentralized} shows that the upper bound of \eqref{eqn: further bound} can be reduced with the network connectivity for a cycle graph.
\begin{Corollary}\label{corollary cycle}
    Suppose that the network is a cycle graph.
    Then the upper bound in \eqref{eqn: further bound} decreases when the network size $N$ is decreased.
\end{Corollary}
Since for the cycle graph,  the smaller network size implies better network connectivity,  Corollary \ref{corollary cycle} also implies the algorithm convergence improves when the connectivity increases. The proof of Corollary 2 is shown in Appendix \ref{app: proof of corollary 2}.


\subsection{Inexact PDC Algorithm}\label{sec: IPDC alg}
Instead of globally solving the subproblem \eqref{eq: dual ADMM x1}, one may choose
to solve it inexactly by using one-step gradient descent only. Specifically, we
replace \eqref{eq: dual ADMM x1} by
\begin{align}
\xb_i^{r+1} &= \xb_i^{r}- \zeta\bigg[ \nabla f_i(\xb_i^r) +
p(\xb_i^r-\zb_i^r)
+ \frac{\Bb_i^\top}{2\rho |\Nc_i|}
\bigg(\Bb_i\xb_i^r - \frac{\qb}{N} - \pb_i^{r+1} + \rho  \Lb_i^+  \yb^r
\bigg) \bigg],
\end{align}
where $\zeta>0$ is a stepsize.
Then we obtain the inexact PDC (IPDC) as summarized in Algorithm \ref{table: IPDC}.

The convergence results of Algorithm \ref{table: IPDC}
are shown in the following theorem.
\begin{Theorem}\label{thm: conv and iteration comp inexact}
Let $\rho>1/2$,  $\beta$ satisfy \eqref{eqn: beta condition},  and $p$, $\zeta$ and $\alpha$ satisfy
\begin{align*}
&\textstyle p > \max\bigg\{ \gamma^+ - 2 \gamma^-, \textstyle \frac{32 ( 2\rho-1)B_{\max}^2}{5\lambda_{\max}^2 \rho^2} - \gamma^-, \textstyle32B_{\max}^2 - \gamma^-,
\textstyle\frac{2\gamma^+ - 7\gamma^-}{5}\bigg\},\\
& \textstyle\frac{1}{p+\gamma^-} < \zeta < \textstyle\min \bigg\{ \frac{5}{2(\gamma^+- \gamma^-)},\textstyle \frac{2}{p+\gamma^+}, \textstyle\frac{1}{32B_{\max}^2}, \textstyle\bigg.\frac{5a_2 \lambda_{\max}^2 \rho^2}{64(2\rho - 1)\theta_1^2 B_{\max}^2 } \bigg\} ,\\
&\textstyle  \alpha <\textstyle\min \bigg\{ \frac{\rho}{5},~ \textstyle\frac{2\rho-1}{16a_1
            \lambda_{\max}^2}, ~\textstyle \frac{\rho \zeta (p+\gamma^-)[1+\rho^2(p+\gamma^+)]}{8 B_{\max}^2[\zeta(p+\gamma^-) + 3]^2}\bigg\},
\end{align*}
where $a_1$ and $a_2$ are defined in \eqref{eq: a 1} and \eqref{eq: a 2}, respectively.
Then,  Algorithm \ref{table: IPDC} has the same convergence properties as
Algorithm \ref{table: PDC}.
  \end{Theorem}

Theorem 2 can be proved following similar ways as for Theorem 1; details are relegated to Appendix \ref{app: inexact}.
Theorem \ref{thm: conv and iteration comp inexact} shows that the IPDC algorithm
enjoys the same iteration complexity as the PDC algorithm. In practice, the inexact
algorithm may exhibit slower convergence speed than its exact counterpart in terms
of the iteration number, but has a substantially smaller computation time; see
\cite{Chang14, chang2016proximal}.

\begin{algorithm}[!t]
	\caption{Proposed IPDC method for solving {\sf (P)}}
	\begin{algorithmic}[1]\label{table: IPDC}
%
%

		\STATE {\bf Given} parameters $p, \beta, \alpha,
\rho$ and stepsize $\zeta$, and initial variables
		$\xb_i^{0}=\zb_i^0 \in \mathbb{R}^{n}$, $\yb_i^{0}\in \mathbb{R}^{M}$   and
$\pb_{i}^{0}=\zerob$ for $i\in [N]$. Set $r=1.$
		\REPEAT
		\STATE  For all  $i\in [N]$ (in parallel), send $\yb_i^r $ to neighbour $j
\in \Nc_i$ and receive $\yb_j^r$ from neighbour $j$. Perform
		\begin{align}
		&\pb_{i}^{r+1}= \textstyle ~ \pb_{i}^{r}+{\alpha}\sum_{j \in
\Nc_i}(\yb_i^{r}-\yb_j^{r}), 			\label{eq: dual ADMM p2}\\
&\xb_i^{r+1}  =  \textstyle \xb_i^{r}- \zeta\bigg[\nabla f_i(\xb_i^r) +
p(\xb_i^r-\zb_i^r)\bigg. \notag \\
&~~~  \bigg.  \textstyle + \frac{\Bb_i^\top}{2\rho |\Nc_i|}
\bigg(\Bb_i\xb_i^r - \frac{\qb}{N} - \pb_i^{r+1} + \rho  \Lb_i^+  \yb^r
\bigg) \bigg],\label{eq: dual ADMM x2}
		\\
	& \textstyle \yb_i^{r+1} \!= \!  \frac{1}{2\rho |\Nc_i|}\big( \Bb_i
\xb_i^{r+1}\! -\! \frac{\qb}{N} \!-\!\pb_i^{r+1}\!  + \rho \Lb^+_i \yb^r\big),\\
& \textstyle \zb_i^{r+1} = \textstyle \zb_i^r +\beta (\xb_i^{r+1} - \zb_i^r),
		\end{align}
		\STATE {\bf Set} $r\leftarrow r+1.$
		\UNTIL {a predefined stopping criterion is satisfied.}
	\end{algorithmic}
\end{algorithm}

\section{Proof of Theorem \ref{thm: conv and iteration comp}}\vspace{-0.0cm}
\label{section: proof}
This section presents the proof of Theorem \ref{thm: conv and iteration comp}.  The
potential function used in the proof is firstly developed, followed by showing some
key error bounds and perturbation bounds. Then,  we prove the
claims in the theorem.

\subsection{ The Potential Function }

To prove the convergence, we need to build a proper potential function that can
monotonically descend with the iteration number. Since our proposed algorithms are
primal-dual proximal algorithms, the potential function should involve both primal
and dual information of \eqref{eq: dual of consensus problem 2 eq}. For ease of
later use, we define
\begin{align}
\mathcal{Y}(\zb)= \arg\max_{\yb
=[\yb_1^\top,\ldots,\yb_N^\top]^\top}~&\sum_{i=1}^N\bigg(\phi_i(\yb_i,\zb_i) -
\yb_i^\top\frac{\qb  }{N}\bigg) \notag \\
	{\rm s.t.}~& \Ab\yb=\zerob\label{eq: y(z)}
\end{align}
as the set of optimal solutions of \eqref{eq: dual of consensus problem 2 eq}, which is non-empty since \eqref{eq: dual of consensus problem 2 eq} is a convex problem satisfying the Slater's condition \cite{boyd2004convex}. By
the definition of $\phi_i$ in \eqref{eq: psi}, we define for $i\in [N]$
\begin{align}
\!\!\!  \xb_i(\yb_i, \zb_i) =   \arg \min\limits_{\xb_i } f_i(\xb_i)
+\frac{p}{2}\|\xb_i - \zb_i\|^2 +  \yb_i^\top \Bb_i \xb_i,  \label{eq: x(y, z)}
\end{align}
and let $\xb(\yb, \zb)= [\xb_1^\top(\yb_1, \zb_1), \cdots, \xb_N^\top(\yb_N,
\zb_N)]^\top$. Besides, recalling $\Lc_\rho$ in \eqref{eqn: AL of
DP2}, the set of optimal primal variables of  \eqref{eq: dual of consensus problem 2
eq}
is given by
\begin{align}\label{eqn: dual function of DP2}
	\mathcal{Y}(\mub, \zb)=\arg\max_{\yb} \Lc_\rho(\yb,\mub;\zb),
	\end{align}
which can be shown to be nonempty given \eqref{eqn: update of p99} (see Appendix \ref{appendix: nonempty oof y(mu, z)}), and the associated objective value is denoted by
$d(\mub;\zb)$.

\vspace{0.0cm}
Next, we show some descent lemmas for $-\Lc_\rho$ and $d$.\vspace{-0.0cm}

\vspace{-0.2cm}
		\begin{Lemma}\label{fact: potential fun L} For the AL function in
\eqref{eqn: AL of DP2},
we have	
		\begin{align}\label{eq: potential function L}
		&-\Lc_\rho(\yb^{r+1},\mub^{r+1};\zb^{r+1}) + \Lc_\rho(\yb^r,\mub^r;\zb^r)
\notag \\
		& \leq  \underbrace{\frac{1}{\alpha}\|\mub^\rpo - \mub^r\|^2}_{\triangleq
\sf (A1)}
-\frac{1}{2} \|\yb^\rpo - \yb^r\|^2_{\rho\Lb^+ } -\frac{\rho}{2} \|\yb^\rpo -
\yb^r\|^2
		\notag \\
		& \underbrace{- \frac{p}{2} \sum_{i=1}^N \la \zb_i^\rpo -\zb_i^r,
\zb_i^\rpo +
\zb_i^r - 2 \xb_i(\yb_i^\rpo, \zb_i^\rpo)\ra}_{\triangleq \sf (A2)}.
		\end{align}	
	\end{Lemma}

\vspace{-0.2cm}
	\begin{Lemma}\label{fact: potential fun d}
		For the dual function  $d(\mub;\zb) $, we
have
		\begin{align}\label{eq: potential function d}
		&d(\mub^{r+1};\zb^{r+1}) - d(\mub^r;\zb^r)   \notag \\
		& \leq \underbrace{\frac{p}{2} \sum_{i=1}^N \la  \zb_i^\rpo -\zb^r,  - 2
\xb_i(\bar \yb_i^{r+1},\zb_i^r) + \zb_i^\rpo +\zb_i^r      \ra}_{\triangleq \sf
(B1)}
\underbrace{-\alpha \la \Ab \yb^r, \Ab \bar \yb^{r}\ra }_{\triangleq \sf
(B2)},
		\end{align}
		where  $\bar \yb^{r+1}=[(\bar \yb_1^{r+1})^\top, \ldots, (\bar
\yb_N^{r+1})^\top]^\top \in \mathcal{Y}(\mub^{r+1}, \zb^{r+1})$ and
	$\bar \yb^{r} \in \mathcal{Y}(\mub^{r+1}, \zb^{r})$	.
	\end{Lemma}
	
We also consider the descent property of $G$ in \eqref{eq: potential function G}.

\vspace{-0.2cm}
\begin{Lemma} \label{Fact G+c}	
Define function
\begin{align}
&\tilde{G}_r(\xb^r, \yb^r, \mub^r; \zb^r)
\triangleq G(\xb^r, \yb^r, \mub^r; \zb^r) +    \|\yb^r\|^2_{\Lb^-}  +
\frac{1}{2}\|\yb^r - \yb^\rmo\|^2_{\rho \Lb^+ }.
\end{align}
Then, we have
	\begin{align}\label{eq: potential function G 2}
	& \tilde{G}_{r+1}(\xb^{r+1}, \yb^{r+1}, \mub^{r+1}; \zb^{r+1})
-\tilde{G}_r(\xb^r, \yb^r, \mub^r; \zb^r)
	\notag \\
	&\leq \underbrace{-\frac{1}{\alpha}\|\mub^\rpo - \mub^r\|^2}_{\triangleq \sf
(C1)} -
\sum_{i=1}^N
\frac{p+\gamma^-}{2} \|\xb_i^\rpo-\xb_i^r\|^2
	\notag \\
	&~~~-\frac{1}{2}\| (\yb^r-\yb^\rmo) - (\yb^\rpo-\yb^r) \|^2_{\rho \Lb^+ }
	\notag \\
	&~~~-\bigg( \frac{\rho}{2}- \frac{\alpha}{2}  \bigg) \|\yb^\rpo -
\yb^r\|^2_{\Lb^-}
\notag \\
	&~~~+  \| \yb^\rpo - \yb^r\|^2_{\rho  \Lb^+ } +
\frac{p}{2}\bigg(1-\frac{2}{\beta}\bigg)\|\zb^\rpo
- \zb^r\|^2.
	\end{align}
\end{Lemma}

The above three lemmas can be proved by similar techniques in \cite{Zhang18,
Hong17}, and one may refer to Appendix \ref{appendix L bound}-\ref{appendix G bound}.
Based on the above three lemmas, let us define the potential function as
		\begin{align}\label{eq: potential fun}
\!\!\!\!	\Phi^r\!=\!  \tilde{G}_r(\xb^r, \yb^r, \mub^r; \zb^r)
	\!-\! 2 \Lc_\rho(\yb^r,\mub^r;\zb^r) \! +\! 2 d(\mub^r;\zb^r).
	\end{align}
By combining Lemmas \ref{fact: potential fun L},  \ref{fact:
potential fun d} and  \ref{Fact G+c}, we have
\begin{align} \label{eq: bound for Phi 1}
\Phi^{r+1} - \Phi^r
& \leq  2{\sf (A1)} + 2{\sf (A2)} + 2{\sf (B1)} + 2{\sf (B2)}+{\sf (C1)}
\notag \\
    &~~~ - \sum_{i=1}^N \frac{p+\gamma^-}{2} \|\xb_i^\rpo-\xb_i^r\|^2
\notag\\
	&~~~-\frac{1}{2}\| (\yb^r-\yb^\rmo) - (\yb^\rpo-\yb^r) \|^2_{\rho \Lb^+ }
\notag
\\
	&~~~-\bigg( \frac{\rho}{2}- \frac{\alpha}{2}  \bigg) \|\yb^\rpo -
\yb^r\|^2_{\Lb^-} + \frac{p}{2}\bigg(1-\frac{2}{\beta}\bigg)\|\zb^\rpo - \zb^r\|^2
\notag \\
	&~~~-{\rho} \|\yb^\rpo - \yb^r\|^2.
\end{align}

\subsection{Perturbation Bounds and Error Bounds}
To show that $\Phi^r$ is a descent function, we will need some perturbation bound
and error bounds.
The following lemma gives the upper bound on $\xb(\yb,\zb)$ when $\yb$ and $\zb$ are
perturbed.\vspace{-0.15cm}
	\begin{Lemma}\label{fact: pertubration of x_i}
Consider \eqref{eq: x(y, z)} for $i \in [N]$. With
$p> - \gamma^-$,
we have
		\begin{align}\label{eq: perturb of xi}
		& \|\xb_i(\yb_i;\zb_i) - \xb_i(\yb_i';\zb_i') \| \leq \sigma_1
\|\zb_i-\zb_i'\| + \sigma_2 \|\yb_i - \yb_i'\|,\\
\label{sigma1}
& \|\xb(\bar \yb^r, \zb^r)-\xb(\bar \yb^{r+1}, \zb^{r+1})\|
 \le {\sigma_3 }\|\zb^r-\zb^{r+1}\|,
		\end{align}
where $\bar \yb^{r} \in \mathcal{Y}(\mub^{r+1}, \zb^{r})$, $\bar \yb^{r+1}\in
\mathcal{Y}(\mub^{r+1}, \zb^{r+1})$, and
$\sigma_1  \textstyle\triangleq \textstyle\frac{p{(p + \gamma^+ + 1)}}
{p+\gamma^-},$
$\sigma_2   \textstyle \triangleq \textstyle\frac{B_{\max}{(p + \gamma^+ + 1)}}
{p+\gamma^-}$ with $B_{\max}  \textstyle \triangleq\textstyle \max_{i\in [N]}
\|\Bb_i\|$, and $\sigma_3 \triangleq \textstyle \frac{p}{p+\gamma^-}$.
	\end{Lemma}
	
	{\bf{Proof}:} See Appendix \ref{appendix: perturbation of x}. \hfill
$\blacksquare$

\vspace{0.2cm}
Next, we determine the error bound that characterizes the distance between
$\yb^{r+1}$ and the solution set $\Yc(\mub,\zb)$ of \eqref{eqn: dual function of
DP2}.
	\begin{Lemma}\label{fact: perturb of y 2}
		With $p>-\gamma^-$, there exists a
$\bar \yb^r\in \mathcal{Y}(\mub^{r+1}, \zb^r)$ of  \eqref{eqn: dual function of DP2}
such that
        \begin{align}\label{eq: y r+1 - y mu r+1 z r L-}
           & \| \yb^{r+1} - \bar \yb^r\|_{\Lb^-}^2
         \leq a_1
        \|\Lb^+(\yb^{r+1} - \yb^r)\|^2, \\
 \label{y^r+1 - y(mu^r+1, z^r) 2}
       & \|\yb^{r+1} - \bar \yb^r\|^2
        \leq a_2
        \| \Lb^+(\yb^{r+1} - \yb^r)\|^2 ,
    \end{align}
where

\vspace{-0.3cm}
{\small \begin{align}
a_1&\triangleq
        \frac{ \frac{\theta_1^2}{2\rho}  \bigg( \frac{1 }{p + \gamma^-} + \rho^2 (p
        + \gamma^+) \bigg)^2+\rho }
       { \frac{1 }{p + \gamma^-} + \rho^2 (p + \gamma^+) },\label{eq: a 1} \\
a_2 &\triangleq
       \frac{\theta_1^4}{\rho^2}  \bigg( \frac{1 }{p + \gamma^-} + \rho^2 (p +
       \gamma^+) \bigg)^2 + 2 \theta_1^2,\label{eq: a 2}
\end{align}    }

\vspace{-0.3cm}
\noindent and $\theta_1>0$ is a constant only depending on the problem instance.
	\end{Lemma}

     {\bf{ Proof}:} See Appendix \ref{appendix: perturbation of y 2}. \hfill
     $\blacksquare$


The following lemma describes the change of $\bar \yb\in \Yc(\mub,\zb)$ when $\zb$
is perturbed.
	\begin{Lemma}\label{fact: perturb of y}
		Consider \eqref{eqn: dual function of DP2} with $p>-\gamma^-$. For any $\bar
\yb^r \in \mathcal{Y}(\mub^{r+1}, \zb^{r})$, there exists a
$\bar \yb^{r+1}\in \mathcal{Y}(\mub^{r+1}, \zb^{r+1})$ such that
\begin{align}\label{eq: bound for y(mu r+1, zr) - y(mu r+1 , z r+1 )}
\|\bar \yb^r - \bar \yb^{r+1}\|^2 \leq {a_3} \|\zb^r-\zb^{r+1}\|^2.
\end{align}
        where
{$
        a_3
        \triangleq \theta_2^2\bigg( (p+\gamma^+)^2{\sigma_3 }^2
+\frac{B_{\max}^2{\sigma_3 }^2}{\rho^2}\bigg)
$}
       and $\theta_2>0$ only depends on the problem instance.
	\end{Lemma}

{\bf{ Proof}:} See Appendix \ref{appendix: perturbation of y}. \hfill
     $\blacksquare$

\vspace{0.2cm}
From above \eqref{y^r+1 - y(mu^r+1, z^r) 2} and \eqref{eq: bound for y(mu r+1, zr) -
y(mu r+1 , z r+1 )}, we can obtain that for $\yb^{r+1}$, there exists a
$\bar \yb^{r+1}\in \mathcal{Y}(\mub^{r+1}, \zb^{r+1})$ so that
\begin{align}\label{eqn: y r+1 to bar y r+1}
              \| \yb^{r+1} -  \bar \yb^{r+1}\|^2
            & \leq 2 a_2 \| \Lb^+ (\yb^{r+1} - \yb^r)\|^2
                +2 a_3\|\zb^{r+1} - \zb^r\|^2.
\end{align}

   Next, we consider the error bound  between \eqref{eq: y(z)} and \eqref{eqn: dual
   function of DP2} due to the violation of the consensus constraint.


\begin{Lemma} \label{Fact y(z^r) - y(mu^{r+1, z^r})}
For any $\bar \yb^r \in \mathcal{Y}(\mub^{r+1}, \zb^r)$, there exists a
$\yb(\zb^r)\in \mathcal{Y}(\zb^r)$ such that
\begin{align} \label{y(z^r) - y(mu^{r+1, z^r})}
   \| {\yb}(\zb^{r}) -\bar  \yb^{r}\|\leq \sqrt{a_4} \|\Ab\bar  \yb^{r}\|,
   \end{align}
   where $
   a_4\textstyle \triangleq {\textstyle \frac{1}{2} \bigg(\frac{\theta_3^2}{
   p+\gamma^- } + \rho^2 \theta_3^2 (p+\gamma^+) \bigg)^2 +
   \frac{\theta_3^2}{\rho^2} }$
for some $\theta_3>0$ only depending on the problem instance.
\end{Lemma}

{\bf{ Proof}:} See Appendix \ref{appendix:  y(z^r) - y(mu^{r+1, z^r})}. \hfill
$\blacksquare$

\vspace{0.2cm}

According to Lemma \ref{fact: pertubration of x_i}, Lemma \ref{fact: perturb of y
2}, and Lemma \ref{Fact y(z^r) - y(mu^{r+1, z^r})}, we can bound the difference
between $(\xb^{r+1},\yb^{r+1})$ generated by Algorithm \ref{table: PDC}
and the primal-dual solution $(\xb(\zb^r),\yb(\zb^r))$ of problem \eqref{consensus
    problem 2 PX} as follows.

\vspace{-0.15cm}
    \begin{Corollary}\label{eb}
    	For $\yb^{r+1}$, there exists a $\bar \yb^r \in \mathcal{Y}(\mub^{r+1},
    \zb^r)$ and $\yb(\zb^r) \in \mathcal{Y}(\zb^r)$ such that
    	\begin{align}
        	&\|\yb^{r+1}-\yb(\zb^r)\|\! \leq\! \sqrt{a_2}\lambda_{\max}
        \|\yb^r\!-\!\yb^{r+1}\| \! +\! \sqrt{a_4}
        \|\Ab \bar \yb^r\|, \label{eb2} \\
        	&\|\xb^{r+1}-\xb(\zb^r)\| \leq
        \sigma_2\|\yb^{r+1}-\yb(\zb^r)\|.\label{eb3}
    	\end{align}
    \end{Corollary}

{\bf{ Proof}:}
To show \eqref{eb2}, we have
    \begin{align} \label{eb2-1}
    \|\yb^{r+1}-\yb(\zb^r)\| & = \|\yb^{r+1}- \bar \yb^r + \bar \yb^r - \yb(\zb^r)\|
    \notag \\
          & \leq \|\yb^{r+1}-\bar \yb^r\| + \| \bar \yb^r -
          \yb(\zb^r)\|.
    \end{align}
By substituting \eqref{y^r+1 - y(mu^r+1, z^r) 2} in Lemma \ref{fact: perturb of y 2}
and  \eqref{y(z^r) - y(mu^{r+1, z^r})} in Lemma \ref{Fact y(z^r) - y(mu^{r+1, z^r})}
into \eqref{eb2-1}, one obtains \eqref{eb2}.
To prove \eqref{eb3}, we first note from
\eqref{eqn: maxmin1}, \eqref{eq: x update} and \eqref{eq: x(y, z)}
  that $\xb^{r+1}=\xb(\yb^{r+1},\zb^r)$. Thus,
by applying \eqref{eq: perturb of xi}, we obtain
   \begin{align}
    \|\xb^{r+1}-\xb(\zb^r)\|
        &
= \|\xb(\yb^{r+1},\zb^r) - \xb(\yb(\zb^r), \zb^r)\| \notag \\
&
\leq \sigma_2 \| \yb^{r+1} - \yb(\zb^r) \|,
   \end{align} which is \eqref{eb3}.
   \hfill $\blacksquare$

With the above results, we are ready to prove Theorem \ref{thm: conv and iteration
comp}.
	
\subsection{ Proof of Theorem \ref{thm: conv and iteration comp}(a)}
We first show that the potential function $\Phi^r$ in  \eqref{eq: potential fun} is
non-increasing and is lower bounded.

\begin{Lemma}\label{Fact descent potential}
	Let $p>-\gamma^-$, $\rho>0$ and

    \vspace{-0.3cm}
    {\small\begin{align}
    &\alpha  \leq \min \bigg\{\frac{\rho}{5}, \frac{\rho}{8a_1
    \lambda_{\max}^2}\bigg\}, \label{eqn: alpha condition}\\
    &\beta <\bigg(\frac{1}{2} + \frac{20 p\sigma_2^2 a_2
    \lambda_{\max}^2}{\rho} + \frac{\rho(\sigma_1^2 + 2\sigma_2^2a_3)}{10
    p\sigma_2^2 a_2 \lambda_{\max}^2} \bigg)^{-1}.\label{eqn: beta condition}
\end{align}}
 \vspace{-0.3cm}

\noindent (a) Then, we have
	\begin{align}
	&\Phi^\rpo -\Phi^r \notag \\
	& \leq - \sum_{i=1}^N \frac{p+\gamma^-}{2} \|\xb_i^\rpo-\xb_i^r\|^2
-C_1 \|\yb^\rpo - \yb^r\|^2_{\Lb^-}
\notag
\\
    &~~~-\frac{1}{2}\| (\yb^r-\yb^\rmo) - (\yb^\rpo-\yb^r) \|^2_{\Lb^+ } -\alpha \|
    \Ab \bar \yb^r\|^2
\notag \\
	&~~~
    -C_2  \|\yb^\rpo - \yb^r\|^2 - C_3  \|\zb^\rpo - \zb^r\|^2. \label{eq: Phi
    descent 2}
	\end{align}
    where

    \vspace{-0.3cm}
    {\small \begin{align}
    \textstyle C_1 & \textstyle \triangleq \frac{\rho - 5 \alpha}{2}\geq0 ,
    ~C_2  \triangleq\rho- \bigg(2\alpha a_1  +  \frac{  {\rho} }{{5
    \lambda_{\max}^2} }
    \bigg)\lambda_{\max}^2>0 ,\notag \\
    \textstyle C_3 & \textstyle \triangleq  p \bigg( -\frac{1}{2} + \frac{1}{\beta}- \frac{20 p\sigma_2^2
    a_2
    \lambda_{\max}^2}{\rho}- \frac{\rho(\sigma_1^2 + 2\sigma_2^2a_3)}{10 p\sigma_2^2
    a_2
    \lambda_{\max}^2}  \bigg)>0. \notag
    \end{align}}
(b) Moreover, $\Phi^r$ is lower bounded, i.e.,
$\sum_{r=0}^T \Phi^{r+1} > \underline{\Phi}$ for some constant
$\underline{\Phi}>-\infty$.
\end{Lemma}

{\bf{ Proof}:} Firstly, consider ${\sf (A1)}$ in \eqref{eq: potential function
L}, ${\sf (B2)}$ in \eqref{eq: potential function d}, $\sf{ (C1)}$ in \eqref{eq:
potential function G 2}, and their combination:
\begin{align}\label{eq: combined term 1}
	&2{\sf (A1)} + 2{\sf
(B2)} + {\sf (C1)} \notag \\
&=\frac{1}{\alpha}\|\mub^\rpo - \mub^r\|^2 - 2\alpha \la \Ab \yb^r, \Ab
\bar \yb^{r}\ra \notag \\
	&\overset{\mathrm{(i)}}{=}\alpha \|\Ab \yb^r\|^2 - 2\alpha \la \Ab \yb^r, \Ab
\bar \yb^{r}\ra \notag \\
	&\overset{\mathrm{(ii)}}{\leq} 2\alpha \|\yb^r - \yb^\rpo\|^2_{\Lb^-}  + 2\alpha
\|\yb^\rpo -
\bar \yb^{r}\|_{\Lb^-}^2
  -\alpha \| \Ab \bar \yb^{r}\|^2\notag  \\
	&\overset{\mathrm{(iii)}}{\leq} 2\alpha \|\yb^r - \yb^\rpo\|^2_{\Lb^-}  +
2\alpha  a_1 \|\Lb^+(\yb^{r+1} -
\yb^r)\|^2
     -\alpha \| \Ab \bar \yb^{r}\|^2,
	\end{align}
where $\mathrm{(i)}$ is due to \eqref{eqn: alg original form mu}, $\mathrm{(ii)}$ is
obtained by adding and subtracting a term $\yb^{r+1}$ followed by using inequality
$\|\ab+\bb\|^2 \leq 2\|\ab\|^2 + 2\|\bb\|^2$, and lastly $\mathrm{(iii)}$ is
owing to \eqref{eq: y r+1 - y mu r+1 z r L-} of Lemma \ref{fact: perturb of y 2}.
	
Secondly, consider 
	\begin{align}
	& {\sf 2(A2)} + 2{\sf (B1)} \notag \\
	& = 2 p \sum_{i=1}^N \la \zb_i^\rpo -\zb_i^r,  \xb_i(\yb_i^\rpo, \zb_i^\rpo) -
\xb_i(\bar \yb_i^\rpo,\zb_i^r) \ra
	\notag
	\\
	& \leq p\delta  \|\zb^\rpo -\zb^r \|^2
+ \frac{p}{\delta}\sum_{i=1}^N\|\xb_i(\yb_i^\rpo, \zb_i^\rpo)-
\xb_i(\bar \yb_i^\rpo,\zb_i^r)\|^2,
            \label{descent potential 3}
    \end{align}
    where the last inequality is obtained by Young's inequality and $\delta >0$.
By further applying \eqref{eq: perturb of xi} of Lemma \ref{fact: pertubration of x_i}
to \eqref{descent potential 3}, we obtain
\begin{align}
& {\sf 2(A2)} + 2{\sf (B1)} \notag \\
	& \leq  p\delta  \|\zb^\rpo -\zb^r \|^2
	       + \sum_{i=1}^N \frac{2p\sigma_1^2}{\delta}\|\zb_i^\rpo-\zb_i^r\|^2
        + \sum_{i=1}^N \frac{2 p\sigma_2^2}{\delta }\|\yb_i^\rpo -
   \bar \yb_i^\rpo\|^2 \notag \\
	& \overset{\mathrm{(i)}}\leq \bigg( p\delta  + \frac{2p\sigma_1^2}{\delta}
\bigg) \|\zb^\rpo-\zb^r\|^2
        +  \frac{2 p\sigma_2^2}{\delta }  \bigg[
	2 a_2 \| \Lb^+ (\yb^{r+1} - \yb^r)\|^2
                +2 a_3\|\zb^{r+1} - \zb^r\|^2   \bigg] \notag \\
	& = \bigg( p\delta  + \frac{2p\sigma_1^2}{\delta} + \frac{4 p\sigma_2^2
a_3}{\delta }\bigg) \|\zb^\rpo-\zb^r\|^2
        + \frac{4 p\sigma_2^2 a_2}{\delta }  \|\Lb^+ (\yb^\rpo - \yb^r)\|^2,
	\label{eq: combined term 2}
	\end{align}
   where $\mathrm{(i)}$ is obtained by \eqref{eqn: y r+1 to bar y r+1}.

By substituting \eqref{eq: combined term 1} and \eqref{eq: combined term 2} into
    \eqref{eq: bound for Phi 1}, we thereby obtain
    \begin{align} \label{eq: bound for Phi 2}
  \Phi^{r+1} - \Phi^r
& = - \sum_{i=1}^N \frac{p+\gamma^-}{2} \|\xb_i^\rpo-\xb_i^r\|^2
-\underbrace{\frac{ \rho - 5 \alpha}{2}}_{\triangleq \rm{C_1}} \|\yb^r -
\yb^{r+1}\|^2_{\Lb^-}
\notag \\
        &~~~ -\frac{1}{2}\| (\yb^r-\yb^\rmo) - (\yb^\rpo-\yb^r) \|^2_{\rho \Lb^+
        }-\alpha
\|\Ab \bar \yb^r\|^2
\notag
  \\
    &~~~- \underbrace{[ \rho - (2\alpha a_1  +  \frac{4 p\sigma_2^2 a_2}{\delta }
    )\lambda_{\max}^2
    ] }_{\triangleq \rm{C_2}}\|\yb^\rpo-\yb^r\|^2
    \notag
 \\
        &~~~ - \underbrace{p ( -\frac{1}{2} + \frac{1}{\beta}-\delta  -
    \frac{2\sigma_1^2}{\delta} - \frac{4 \sigma_2^2 a_3}{\delta }  )
    }_{\triangleq \rm{C_3}}\|\zb^\rpo-\zb^r\|^2.
\end{align}

We see that if $\alpha \leq \rho/5$, then $C_1$ is nonnegative.
In addition, choose $\delta = \frac{{20} p\sigma_2^2 a_2 \lambda_{\max}^2}{\rho} >
\frac{16 p\sigma_2^2 a_2 \lambda_{\max}^2}{\rho}$, i.e., $ \frac{4 p\sigma_2^2
a_2\lambda_{\max}^2 }{\delta }
    < \frac{\rho}{4}$, then as
long as  $ \alpha  \leq  \frac{\rho}{8a_1 \lambda_{\max}^2},$
we can have $C_2 \geq  {\rho/2} > 0$.
Lastly, if $\beta$ satisfies \eqref{eqn: beta condition},
then $C_3$ is positive.
Thus, Lemma \ref{Fact descent potential}(a) is proved.

The proof for showing $\Phi^r$ is lower bounded (Lemma \ref{Fact descent
potential}(b)) is referred to  Appendix \ref{appendix: lemma 8b}.
\hfill $\blacksquare$

	According to Lemma \ref{Fact descent potential}, we have
    \begin{align}
    &\|\xb^r-\xb^{r+1}\|\rightarrow 0,~~\|\yb^r-\yb^{r+1}\|\rightarrow 0,
    \label{conv 1} \\
    &\|{\Ab}\bar \yb^{r}\|\rightarrow
    0,~~\|\zb^{r}-\zb^{r+1}\|\rightarrow 0.
    \end{align}
    By \eqref{eqn: alg original form z} and {$\|\zb^{r}-\zb^{r+1}\|\rightarrow 0$},
    one can have
    $\|\xb^{r+1}-\zb^{r+1}\|\rightarrow 0$. Besides, combing \eqref{y^r+1 -
    y(mu^r+1, z^r) 2} and
    {$\|{\Ab}\bar \yb^{r}\|\rightarrow
    0$}, one then obtains $\|{\Ab}\yb^{r+1}\|\rightarrow 0$. Therefore,
    we obtain
    \begin{align}
        \|\xb^{r+1}-\zb^{r+1}\| \rightarrow 0,~~ \|{\Ab}\yb^{r+1}\|\rightarrow
        0.\label{conv 2}
    \end{align}
	By applying \eqref{conv 1} and \eqref{conv 2} to the KKT conditions of
\eqref{eqn: alg original form y}, one then concludes that every limit point of
$\{(\zb^r, \yb^r)\}$ is a KKT solution of {\sf (P)}
\hfill $\blacksquare$

\vspace{-0.2cm}
\subsection{ Proof of Theorem \ref{thm: conv and iteration comp}(b)}

Since $\Phi$ is non-increasing and lower bounded by
Lemma \ref{Fact descent potential},  one can show that for $r>0$, there exists a
$t\le r$ such that
\begin{align}
  \textstyle \Phi^t-\Phi^{t+1}\le \frac{\Phi^0-\underline{\Phi}}{r}.\label{conv 3}
\end{align}
Then, by \eqref{eqn: alg original form z}, \eqref{eb2} and \eqref{eb3}, we can bound
	\begin{align}\label{eq: bound for z x(z) 1}
	&\|\zb^t-\xb(\zb^t)\|\leq \|\zb^t-\xb^{t+1}\| +\| \xb^{t+1} - \xb(\zb^t)\|
\notag  \\
	&\leq \frac{1}{\beta}\|\zb^t-\zb^{t+1}\|+\sigma_2
\big(\sqrt{a_2}\lambda_{\max}\|\yb^t\!-\!\yb^{t+1}\| \!+
\sqrt{a_4} \|\Ab \bar \yb^t\|\big)
     \notag \\
	&\overset{\mathrm{(i)}}{\leq} {\textstyle \sqrt{\frac{\Phi^0-\underline{\Phi}}{r}}} \bigg(\frac{1}{\beta
\sqrt{C_3}} + \frac{\sigma_2\sqrt{a_2}\lambda_{\max}}{\sqrt{C_2}}  +
\frac{\sigma_2
\sqrt{a_4}}{\sqrt{\alpha}}\bigg)\notag \\
& \triangleq \tilde{C} \sqrt{\frac{1}{r}},
\end{align}where $\mathrm{(i)}$ is obtained by applying \eqref{eq: Phi descent
2}.

Next, we connect  \eqref{eq: bound for z x(z) 1} with an KKT solution of problem {\sf (P)} by the following lemma,  which is proved in Appendix \ref{appendix: proof of lemma 9}.
 \begin{Lemma}\label{lemma9}
           Suppose that $\|\zb-  \xb(  \zb)\|^2\leq \epsilon$ for some $\zb$. Then, $\zb$ is a $\kappa \epsilon$-KKT solution of problem {\sf (P)},
 where $  \kappa  = \max\{ 2(p^2+(\gamma^+)^2, NB_{\max}^2\} $ and $B_{\max} \triangleq \max_{i\in [N]} \|\Bb_i\|$.
        \end{Lemma}

By the above lemma,  we conclude that $\zb^t$ satisfying \eqref{eq: bound for z x(z) 1} is a $\big(\kappa \tilde{C}^2/r \big)$-KKT solution of problem {\sf (P)}.
\hfill $\blacksquare$

\vspace{0.0cm}
\section{Numerical Results}\label{section: simulation results}
In this section, we present two numerical examples to illustrate the empirical
performance of the proposed algorithms.

\vspace{-0.3cm}
\subsection{Distributed Logistic Regression} \label{Numerical-LR}
Following  \eqref{distributed-feature-general}, we formulate a non-convex
regularized
logistic regression problem \cite{antoniadis2011penalized} with
\begin{subequations}\label{LR problem}
\begin{align}
&\psi(\bb_{k}^\top\wb; v_k)=\log{\left( 1+e^{-
v_k(\bb_{k}^\top \wb)}\right) }, ~\forall k\in [M],\\
 &R_i(\wb_i)= \lambda
\sum_{s=1}^{n}\frac{ \xi w_{i, s}^2 }{1+\xi w_{i, s}^2},~\forall i\in [N],
\end{align}
\end{subequations}
where $v_k\in \{\pm 1\}$ are binary labels, and $\lambda, \xi >0$ are the
parameters.

\begin{figure}[!t]
\begin{minipage}[b]{1.0\linewidth}
  \centering
 \epsfig{figure=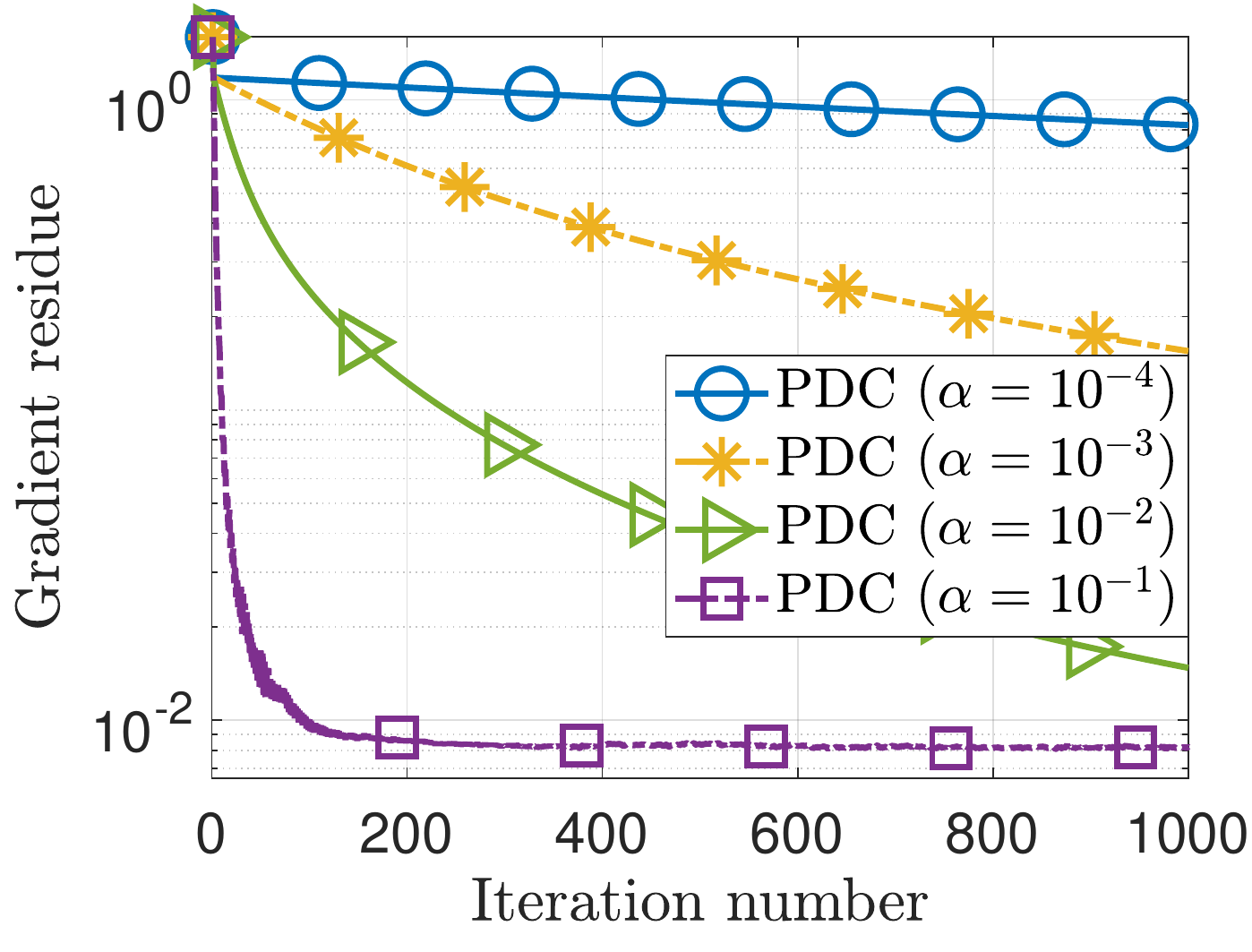,width=2.3in}
 \epsfig{figure=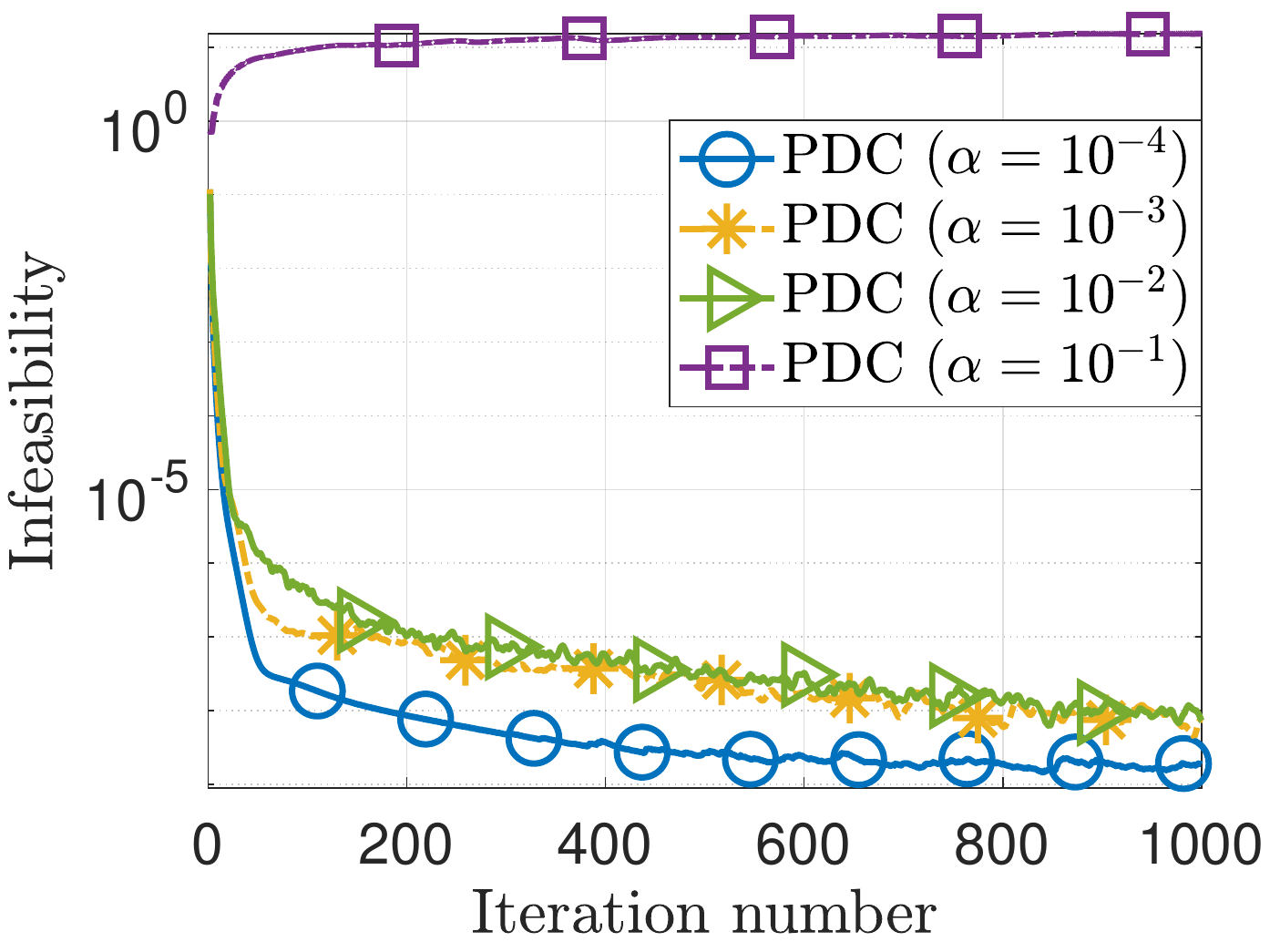,width=2.3in}
  \centerline{\small{(a) $\beta=0.1$, $ \rho = 0.01$, $p = 0.01$,
and various values of $\alpha$}}\medskip
\end{minipage}
\begin{minipage}[b]{1.0\linewidth}
  \centering
 \epsfig{figure=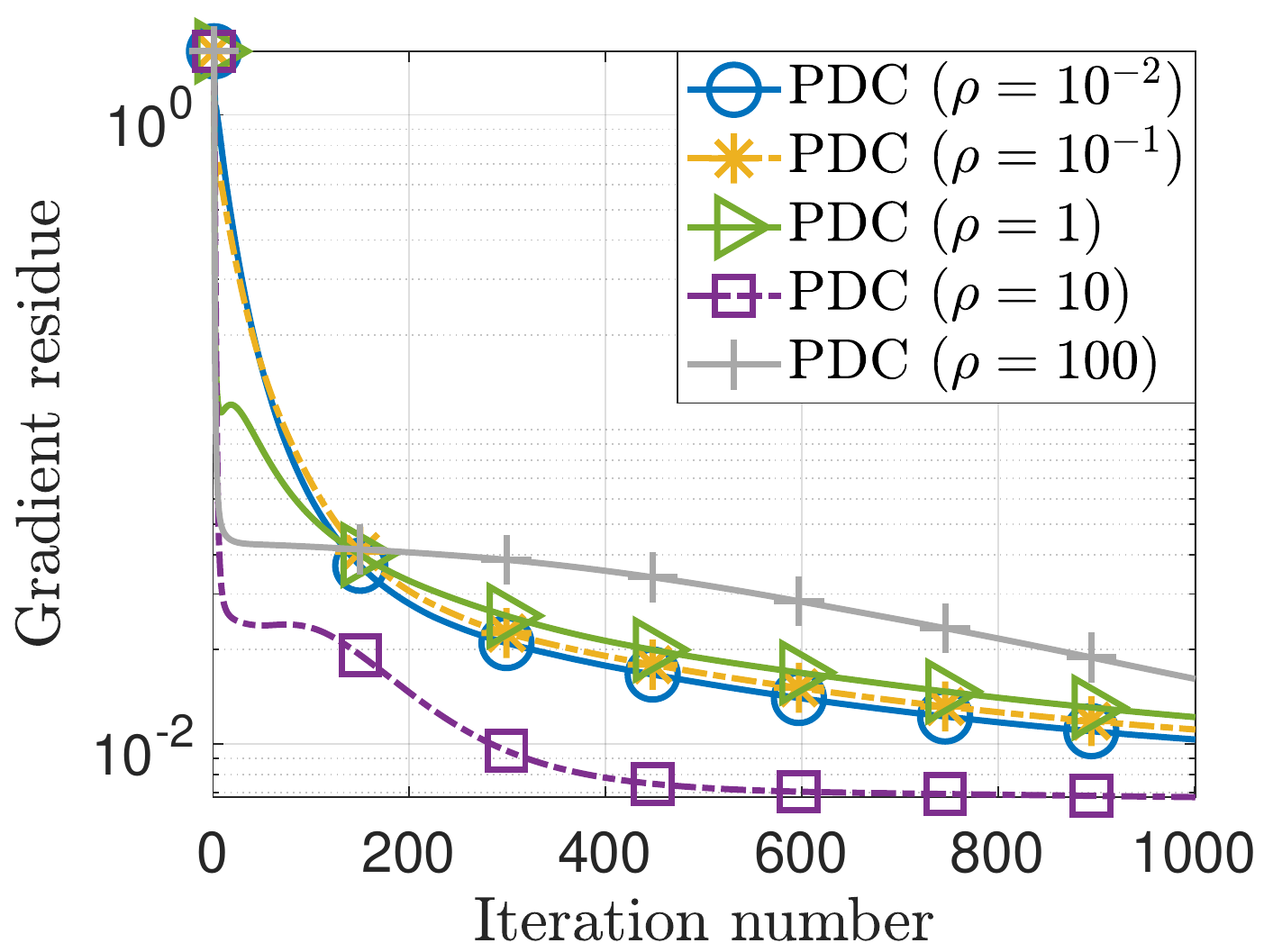,width=2.3in}
 \epsfig{figure=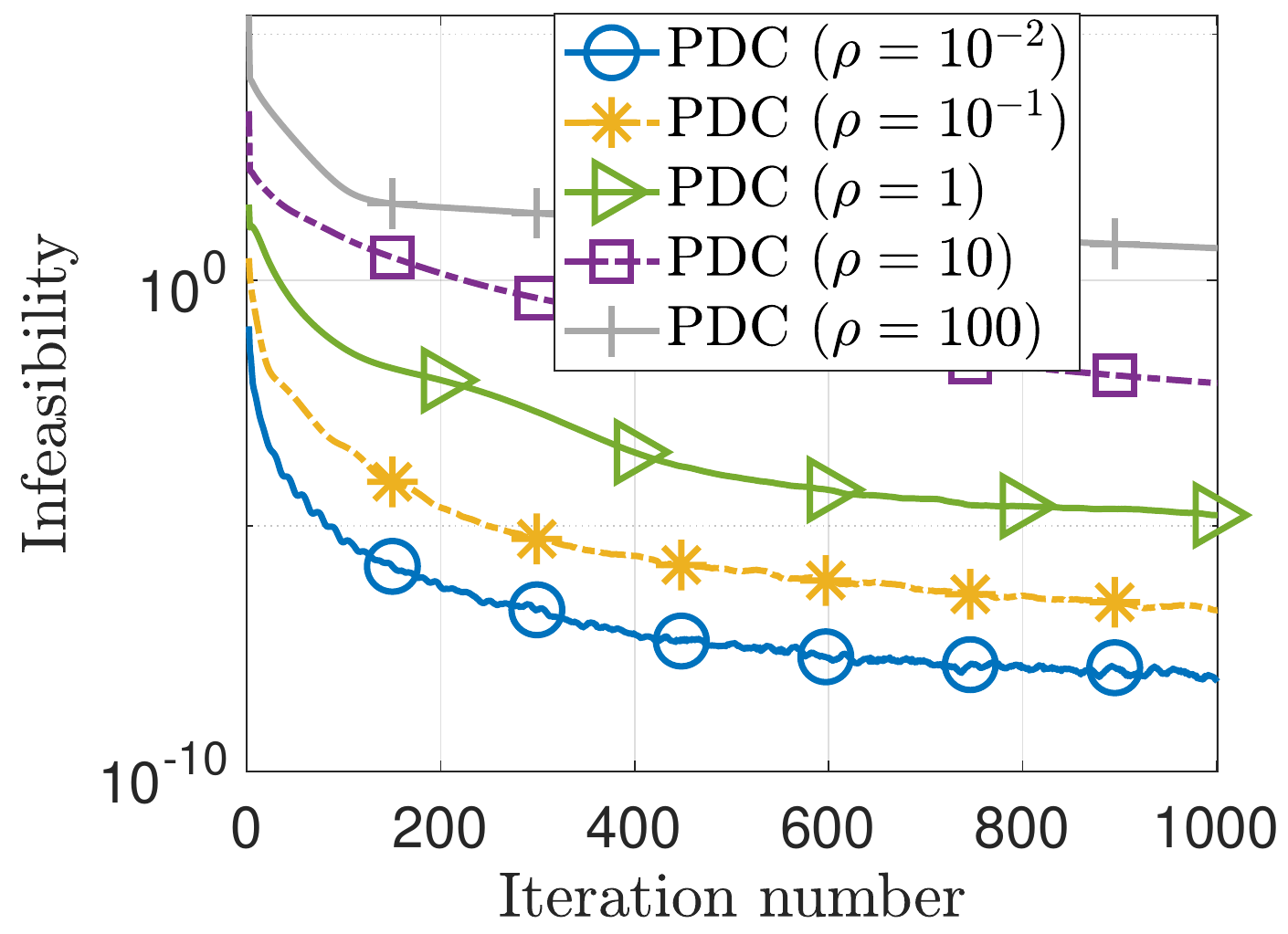,width=2.3in}
  \centerline{\small{(b) $\beta=0.1$, $p = 0.01$, $\alpha = 0.01$,
and various values of $\rho$}}\medskip
\end{minipage}
\caption{Convergence performance of the PDC algorithm with various values of $\alpha$ and $\rho$.} \label{fig: LR PDC}\vspace{-0.4cm}
\end{figure}

We consider the two images D24 and D68 of Brodatz data
set\footnote{http://www.ux.uis.no/tranden/broadatz.html}, and extract $M/2$
overlapping
patches with dimension $\sqrt{Nn} \times \sqrt{Nn}$ from these two images,
respectively. Each patch is vectorized into $1\times Nn$-vector, and they are
randomly shuffled and stacked into the data matrix $\Bb$.
Here we set $M = 100$ and $n = 100$, and assumed a network with $N = 25$ nodes. The
generation of the network graph follows the same method as in
\cite{YildizScag08}. For \eqref{LR problem}, we set $\lambda = 0.01$ and $\xi =
0.5$.

\begin{figure}[!t]
\begin{minipage}[b]{1.0\linewidth}
  \centering
 \epsfig{figure=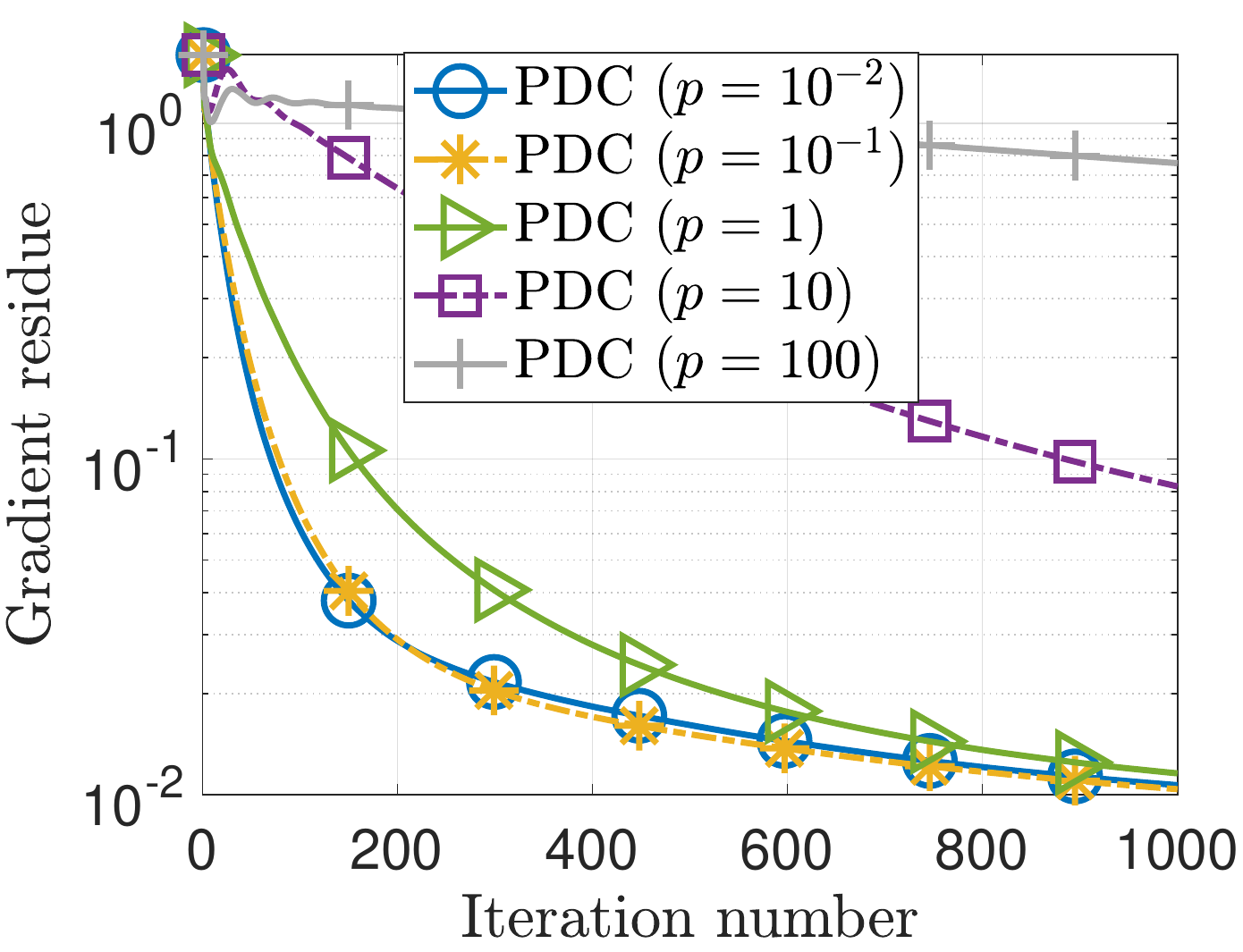,width=2.3in}
 \epsfig{figure=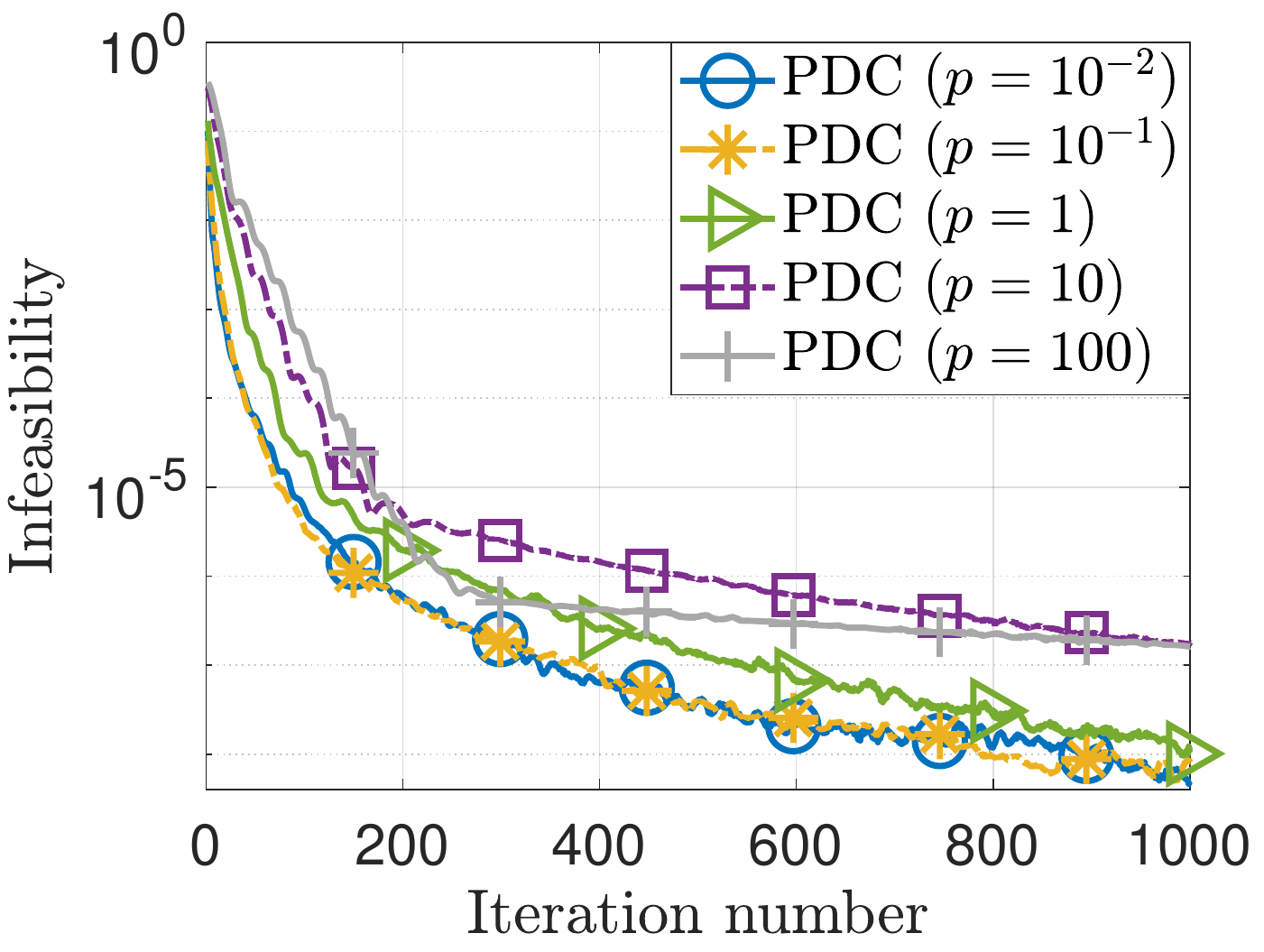,width=2.3in}
  \centerline{\small{(a) $\beta=0.1$, $ \rho = 0.01$, $\alpha = 0.01$,
and various values of $p$}}\medskip
\end{minipage}
\begin{minipage}[b]{1.0\linewidth}
  \centering
 \epsfig{figure=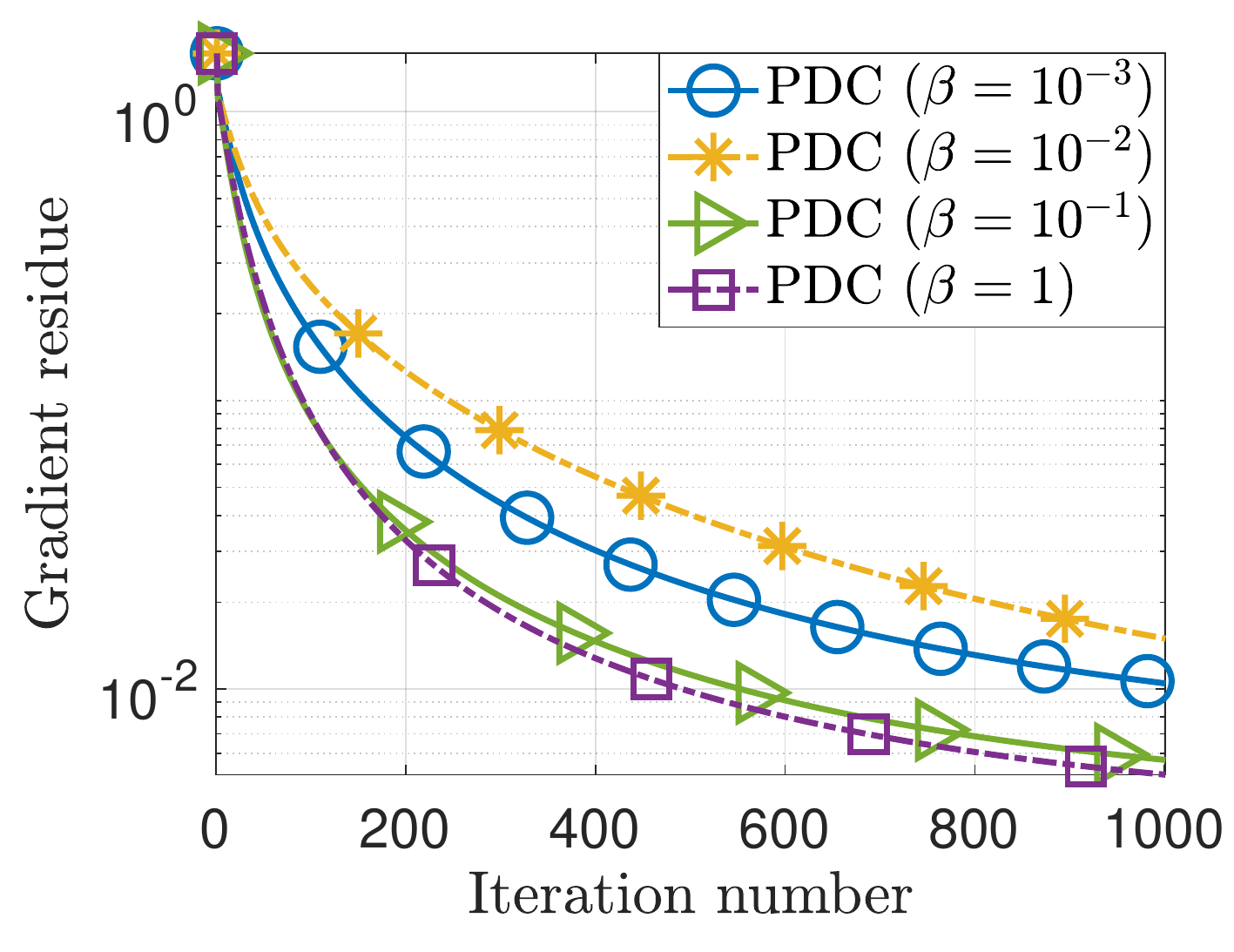,width=2.32in}
 \epsfig{figure=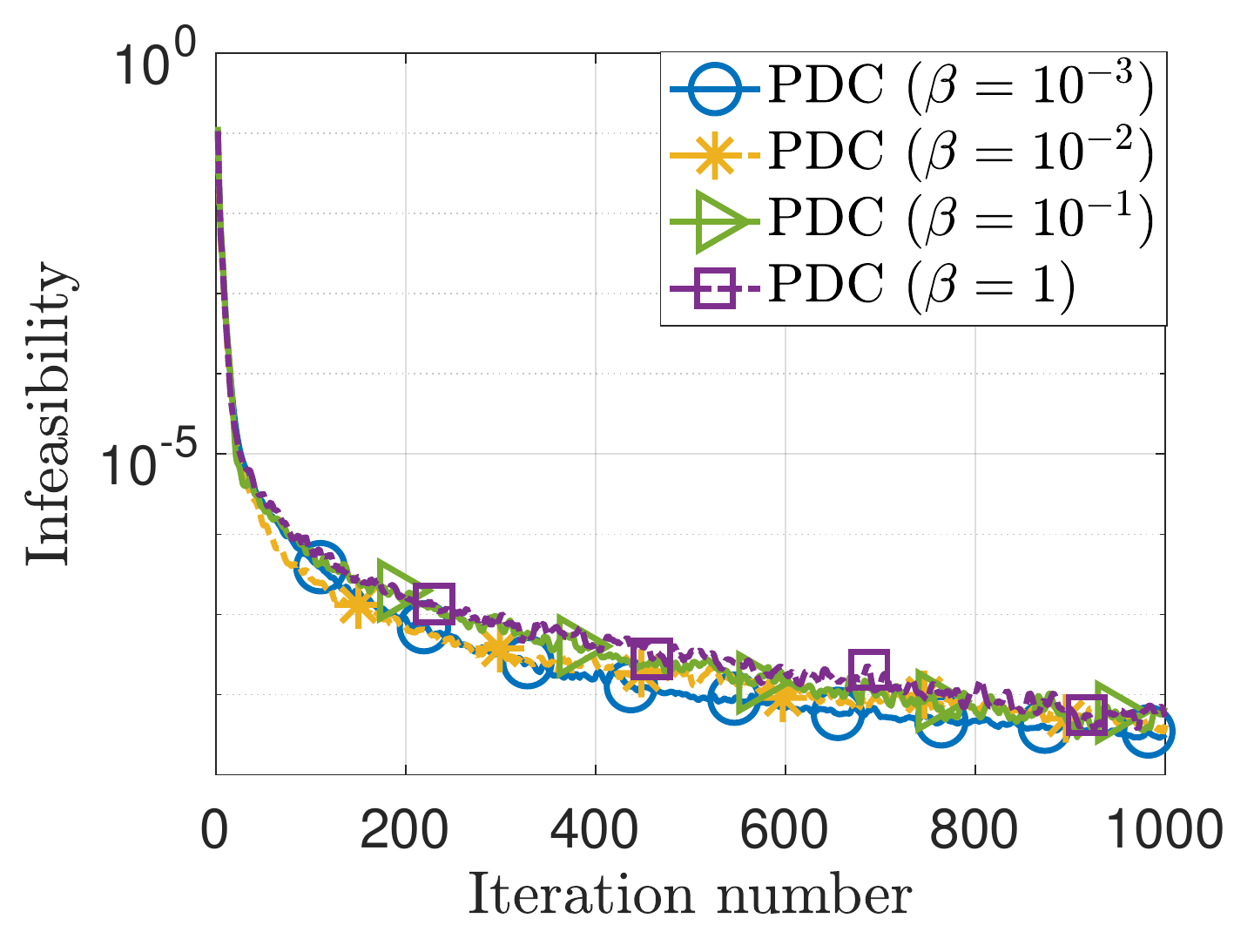,width=2.33in}
  \centerline{\small{(b) $p =0.01 $, $  \rho = 0.01$, $\alpha = 0.01$,
and various values of $\beta$}}\medskip
\end{minipage}
\caption{Convergence performance of the PDC algorithm with various values of $p$ and $\beta$.} \label{fig: LR PDC p rho}\vspace{-0.4cm}
\end{figure}

We employ the proposed PDC algorithm (Algorithm \ref{table: PDC}) to handle the
non-convex LR problem in the form of \eqref{consensus problem 22}.
Our aim is to use this example to examine how the algorithm parameters $\alpha$, $p$, $\beta$, and $\rho$ can affect the convergence of proposed PDC algorithms.
The fast iterative shrinkage thresholding algorithm (FISTA)
\cite{BeckFISTA2009} is used to solve subproblem \eqref{eq: dual ADMM x1}, and the
stopping condition of FISTA is that the normalized proximal gradient
\cite{BeckFISTA2009} is smaller than $10^{-5}$.
The entries of the initial variables $\xb^0$ and $\yb^0$ are randomly generated
following a uniform distribution over the interval $[-1, 1]$.
The experiments are performed 10 times, each with a different initial ($\xb^0$,
$\yb^0$), and the averaged results are presented.
In particular, we assess the algorithm performance by
calculating the following two terms
\begin{align}
 &{\rm Graident~ residue}= \textstyle \frac{1}{Nn} \sum_{i=1}^N\|\nabla
 f_i(\xb_i^{r}) +
 \Bb_i^\top\yb_i^{r}\|^2, \label{gradient residue}\\
 & {\rm Infeasibility}= \textstyle\frac{1}{M}\|\sum_{i=1}^N \Bb_i \xb_i^{r} -
 \qb\|^2,\label{infeasibility}
\end{align}
with respect to the iteration number (communication round) $r$.  The simulation results are displayed in Fig.~\ref{fig: LR PDC}, where we respectively examine the impacts of the four parameters $\alpha$, $\rho$, $p$ and $\beta$ on the algorithm convergence.

In Fig.~\ref{fig: LR PDC}(a),
we set $\beta=0.1$, $\rho = 0.01$, $p = 0.01$,  and $\alpha\in \{10^{-1},10^{-2},10^{-3},10^{-4}\}$.
One can see that
a smaller $\alpha$ may slow down
the proposed
PDC algorithm in terms of the gradient residue (left figure),  whereas slightly reducing the feasibility level (right figure).  It is observed that when $\alpha =
10^{-1}> \rho$, the algorithm cannot satisfy the linear constraint even though it achieves the smallest gradient residue.  Such observation is consistent with Theorem 1.

In Fig.~\ref{fig: LR PDC}(b),
we set $\beta=0.1$, $p = 0.01$, $\alpha = 0.01$,  and $\rho\in \{10^{-2},10^{-1},1,10,10^2\}$.  One can see from the two figures that parameter $\rho$ does not
significantly impact the convergence of gradient residue,  but can greatly affect the constraint feasibility.  In particular,  a larger $\rho$ can cause larger infeasibility values.

In Fig.~\ref{fig: LR PDC p rho}(a),
we set $\beta=0.1$, $\rho = 0.01$, $\alpha = 0.01$,  and $p\in \{10^{-2},10^{-1},1,10,10^2\}$.  One can see from the two figures that parameter $p$ does not
significantly impact the constraint feasibility,  whereas a larger $p$ can slow down the convergence of gradient residue.  This is expected since with a larger $p$,  problem \eqref{consensus problem 2 PX} is a more conservative convex approximation, and thus it slows down the algorithm convergence.

Lastly, in Fig.~\ref{fig: LR PDC p rho}(b),
we set $p=0.01$, $\rho = 0.01$, $\alpha = 0.01$,  and $\beta \in \{10^{-3},10^{-2},10^{-1},1\}$.  One can see
that the step size $\beta$ has a mild influence on the convergence of both gradient residue and infeasibility. Interestingly, while Theorem 1
suggests that $\beta$ should be less than one and a small number, it is found that
the algorithm behaves well with $\beta = 1$.

In Appendix \ref{Additional Numerical-LR}, we further present the convergence results for the IPDC algorithm with respect to the algorithm parameters $\alpha$,  $\rho$, $p$, $\beta$ and $\zeta$.

%
%
%
%

%
%

\vspace{-0.3cm}
\subsection{Distributed Neural Network}
\label{sec: simulation NN}

\begin{figure}[!t]
\begin{minipage}[b]{1.0\linewidth}
  \centering
 \epsfig{figure=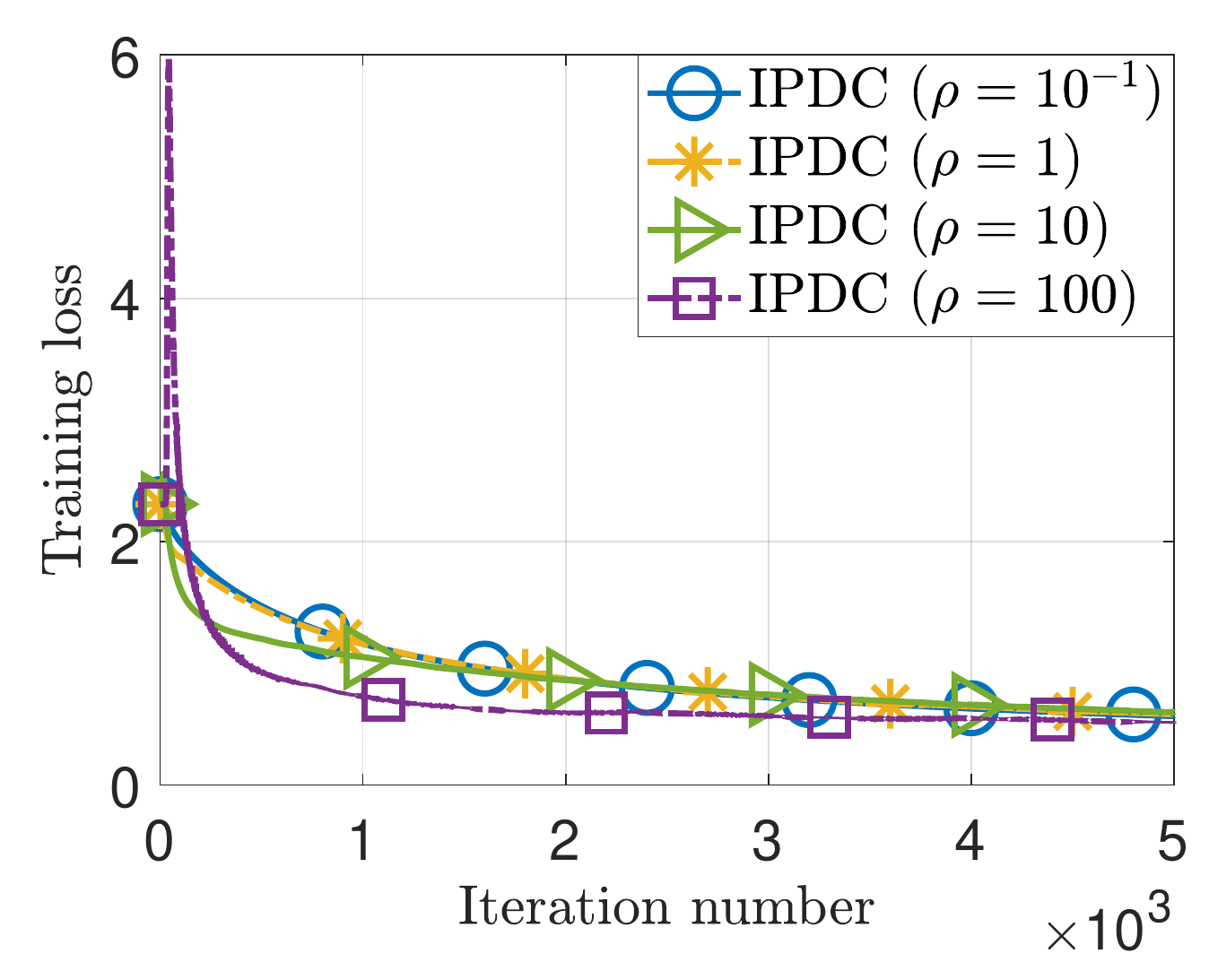,width=2.295in}
 \epsfig{figure=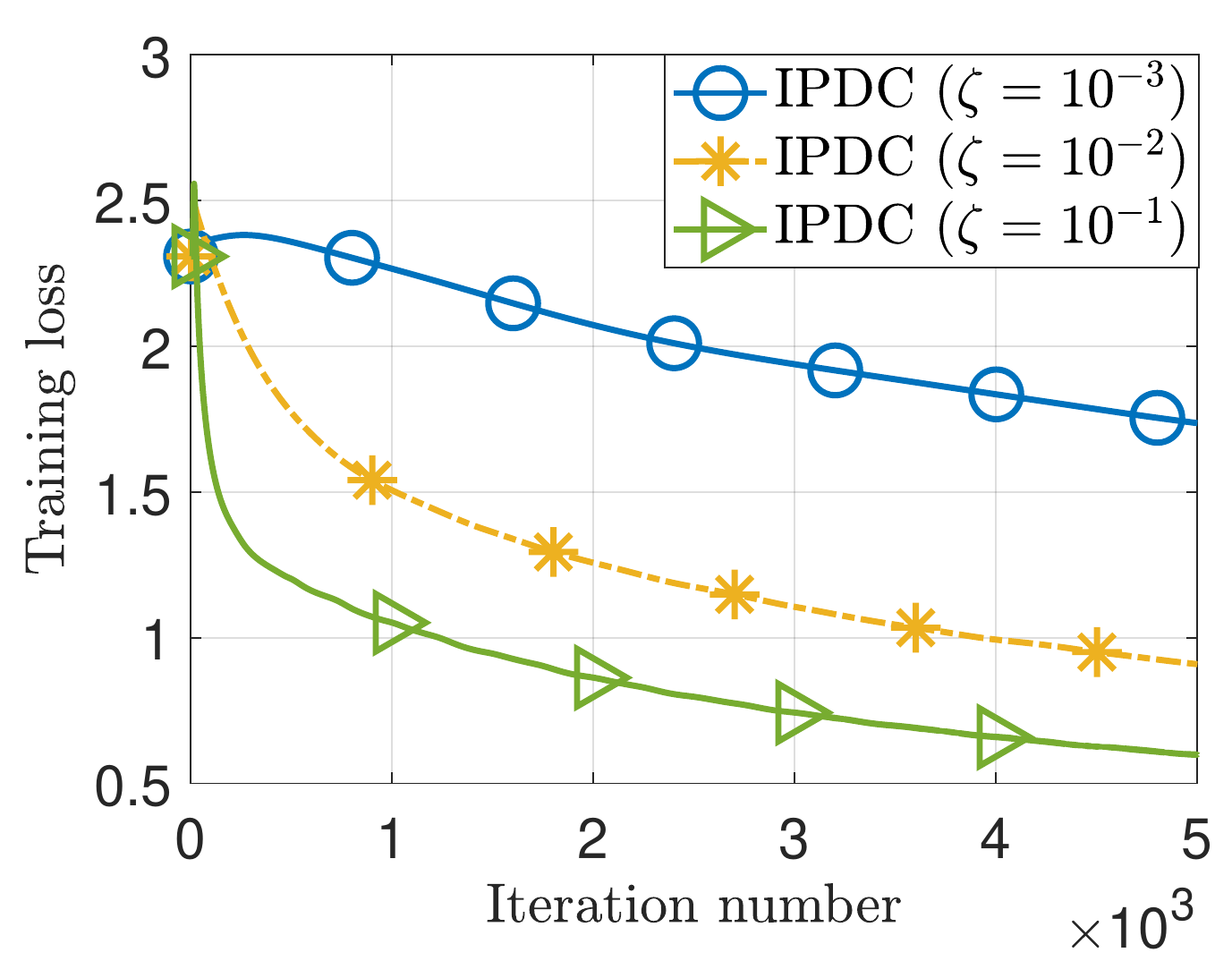,width=2.3in}
  \centerline{\small{(a) Training loss}}\medskip
\end{minipage}
\begin{minipage}[b]{1.0\linewidth}
  \centering
 \epsfig{figure= 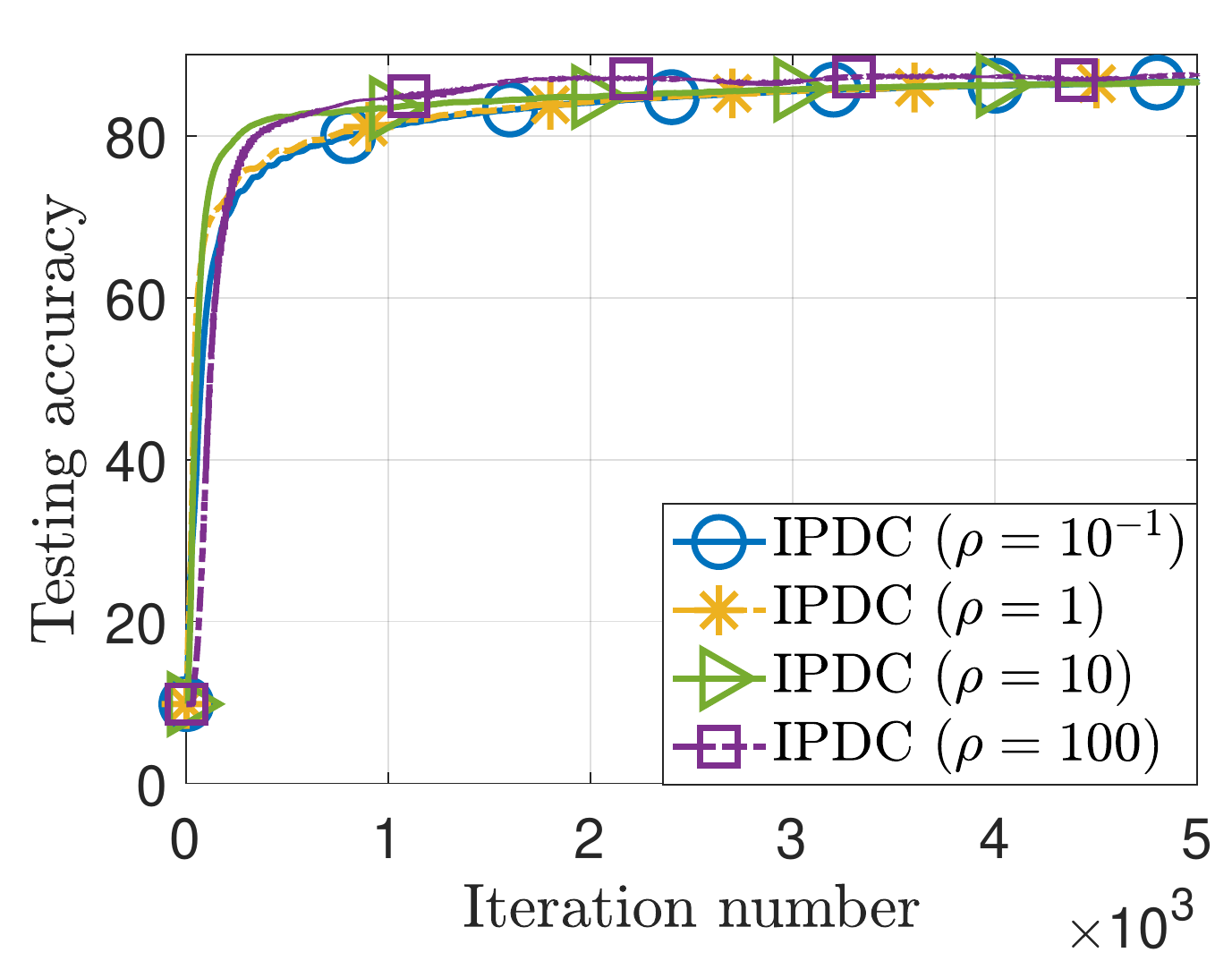,width=2.3in}
 \epsfig{figure= 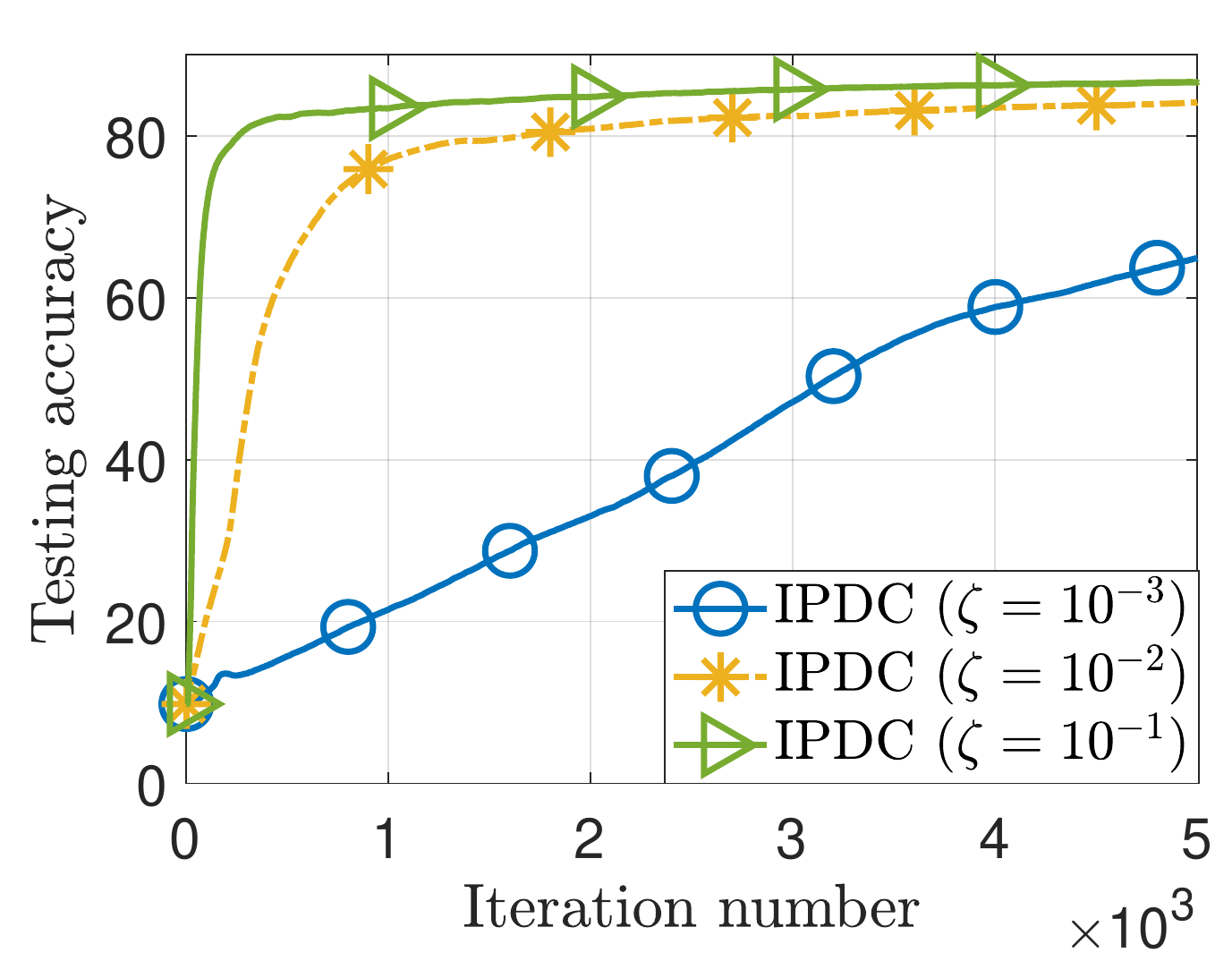,width=2.3in}
  \centerline{\small{(b) Testing accuracy}}\medskip
\end{minipage}
\caption{\small{Convergence performance of the IPDC algorithm with
$p = 10$, $\beta = 0.01$, $\alpha = 0.01$, and various values of $\rho$ and $\zeta$.
Left column: $\zeta=0.1$; Right column: $\rho=10$. }}\label{fig: NN
PDC IPDC alpha}
\vspace{-0.4cm}
\end{figure}
In this simulation, following \eqref{NN problem}, we consider the use of the NN to
classify handwritten digits in the MNIST dataset. The $60000$ training images are divided into $12$ batches. The dimension of each image is $28\times 28$ which is vectorized into a $1 \times
{784}$ vector ($nN=784$).  These vectors of the 5000 samples are stacked
as the data matrix $\Bb$.  Then this data matrix is column partitioned into $8$ equal parts to be distributed to $8$ agents ($N=8$). The network topology is generated in the same way as in \cite{YildizScag08}.

This NN has 2 layers, with $K = 30$ neurons in the hidden layer. The activation
function of the hidden layer and output layer are Rectified Linear Unit (ReLU) and
softmax function, respectively. The cross-entropy loss function
\cite{Goodfellow-et-al-2016} is applied in the last layer.
Then, the proposed IPDC algorithm are applied to train the
classification NN by solving problem \eqref{eq: distributed NN
problem }. The experiment is performed on the PyTorch platform.

For the IPDC algorithms, it is set that $p= 10$, $\alpha = 0.01$, $\beta = 0.01$, $\rho \in \{ 10^{-1}, 1, 10, 100\}$ and $\zeta   \in \{10^{-3}, 10^{-2},10^{-1} \}$.
Moreover, the simulations are preformed with $3$
randomly generated initial points, and their average results are presented in Fig.~\ref{fig: NN PDC IPDC alpha}.

Fig.~\ref{fig: NN PDC IPDC alpha} displays the training loss and testing accuracy of the NN versus the iteration number.
From the left column of Fig.~\ref{fig: NN PDC IPDC alpha}, one can observe that  the
IPDC algorithm with a larger $\rho$ has a sharper decrease in the training loss and increase in the testing accuracy in the first few iterations,  but the algorithm can eventually achieve similar training loss and testing accuracy for the tested values of $\rho$.
In the right column of Fig.~\ref{fig: NN PDC IPDC alpha},  one can observe that the step size $\zeta$ can influence the performance of the proposed IPDC algorithm significantly. Though it is shown in Theorem 2 that a sufficiently small $\zeta$ is required for the algorithm convergence of the IPDC algorithm,  it is found from the figure that  with a larger $\zeta$ the IPDC algorithm converges faster in both training loss and testing accuracy. Lastly,  it is interesting to see that,  while our theoretical claims in Theorem \ref{thm: conv and iteration comp inexact} are under the smoothness assumption (Assumption \ref{assumption smooth}),  the simulation results above demonstrate that the proposed algorithm can still work well for the non-smooth ReLU neural network.

\vspace{-0.1cm}
\section{Conclusions} \label{section: conclusions}
In this paper, two new distributed algorithms, i.e., the PDC and the IPDC
algorithms, for solving the
non-convex linearly constrained problem ${\sf (P)}$ over the multi-agent network
have been proposed. These algorithms can be regarded as extensions of the dual
consensus ADMM method in \cite{Chang14} to the challenging non-convex settings, by
leveraging the proximal technique in \cite{Zhang18}.
While the algorithm development is rather straightforward, it turns out that proving
the convergence and convergence rate of the proposed algorithms are not. By
developing some key error bounds and perturbation bounds, we have shown that, under
mild conditions on the algorithm parameters, the proposed algorithms both can
converge to a KKT solution of problem ${\sf (P)}$ with a complexity order of
$\mathcal{O}(1/\epsilon)$. Through experiments on a non-convex logistic regression
problem and an NN-based MNIST classification task,
we have demonstrated the good convergence behaviors of the proposed algorithms with
respect to various parameter settings.

To the best of our knowledge, the presented distributed algorithms and convergence analyses are the first for the linearly constrained non-convex problem ${\sf (P)}$. In the future, it is worthwhile to extend the current PDC framework to mini-batch
stochastic gradient descent (SGD) \cite{DPSGD_2017,xin2020near},  local SGD
{\cite{CE_DDNN_2017},  online learning \cite{xin2021improved},
and to handle more complex learning problems with heterogeneous data distribution \cite{tang2018d},  and hybrid partitioned data samples and features \cite{FML_CA_2019}.
Inspired by the simulation results in Section \ref{sec: simulation NN},  it is also meaningful to
analyze convergence properties of the proposed algorithm for ReLU neural networks \cite{zou2018stochastic, zhang2019learning} or devise distributed algorithms for general non-smooth problem ${\sf (P)}$.
Finally,  as we discuss after Corollary \ref{lem: convergence graph},  it needs more efforts to understand how graph topology and data matrix can affect the algorithm convergence.  Besides,  it is also important to devise optimal distributed schemes as in \cite{sun2019distributed} for the general problem {\sf (P)}.

\appendix
\appendixpage
\addappheadtotoc

\setcounter{subsubsection}{0}	
\setcounter{section}{0}
\setcounter{equation}{0}
\renewcommand{\theequation}{\thesection.\arabic{equation}}

\vspace{-0.3cm}
\vspace{0.3cm}
\section{Proof of nonempty property of $\Yc(\mub, \zb)$}\label{appendix: nonempty oof y(mu, z)}
To prove that the set $\Yc(\mub, \zb)$ is not empty, we first give a general lemma and its proof.
\begin{Lemma}\label{thm: existence of saddle points}
Let $F(x, y)=f(x)+y^\top Ax+v^\top y-\frac{1}{2}y^\top Qy$, where $f(x)$ is a strongly convex, smooth  function, $Q$ is a positive semi-definite matrix and $v\in \mathrm{Image}(Q)$. Then the saddle point of $F(x, y)$ exists, i.e., the solution of the problem $P1$ exists
\begin{align}
P1: \min_{x }\max_{y }F(x, y).
\end{align}
\end{Lemma}
{\bf Proof of Lemma \ref{thm: existence of saddle points}:}
First, by the basic property of quadratic optimization, $\arg\max_{y} F(x, y)$ exists if and only if
\begin{align}
Ax\in \mathrm{Image}(Q).
\end{align}

Let $\mathrm{rank}(Q)=r$.
Without loss of generality, we assume $Q=\begin{bmatrix}
Q_{11}&Q_{12}\\
Q_{12}&Q_{22}
\end{bmatrix}$ and $Q_{11}\in \mathbb{R}^{r\times r}$ is of full rank.
Let  $P=\begin{bmatrix}Q_{11}\\Q_{12} \end{bmatrix}$ and $Y=\{y\in \mathbb{R}^m\mid y_{r+1}=\cdots =y_m=0\}$.

Also let $\mathrm{Image}(Q)=\{w\mid Bw=0\}$, where the rows of $B\in \mathbb{R}^{(m-r)\times m}$ are the basis for $\mathrm{Null}(Q)$.
Define $X=\{x\in \mathbb{R}^n\mid BAx=0\}$. Then, $x \in X$ implies $x \in \mathrm{Image}(Q)$. For any $\tilde{x}\in X$, there exists some $\tilde{y}\in Y$ such that
\begin{align}
\tilde{y}=\arg\max_{y\in \mathbb{R}^m}F(\tilde{x}, y).
\end{align}
In fact, the matrix $P$ is of rank $r$. Therefore, $\mathrm{Image}(Q)=\mathrm{Image}(P)$.
Notice that $-A\tilde{x}-v\in \mathrm{Image}(Q)=\mathrm{Image}(P)$.
It implies that there exists some $u\in \mathbb{R}^r$ such that
$$-A\tilde{x}-v=Pu.$$
Hence, if $\tilde{y}=(u^T, 0, \cdots, 0)^T\in Y\subseteq \mathbb{R}^m$, we have
$$\nabla_yF(\tilde{x}, \tilde{y})=A\tilde{x}+v+Q\tilde{y}=0.$$
Consequently,
For any $\tilde{x}\in X$, there exists some $\tilde{y}\in Y$ such that
$$\tilde{y}=\arg\max_{y\in Y}F(\tilde{x}, y).$$
Then consider the following constrained strongly convex-strongly concave min-max problem:
\begin{align}
P2: \min_{x\in X}\max_{y\in Y}F(x, y).
\end{align}
Note that $F(x, y)$ is strongly concave for $y$ in $Y$, thus $P2$ is a strongly convex-strongly concave min-max problem. Since a strongly convex (concave) function has a compact sublevel set, by \cite[Prop. 3.6.9]{bertsekas2003convex}, there exists a saddle point for problem $P2$.
Suppose that $(\bar{x}, \bar{y})$ is a saddle point of P2.
Since  $\bar{x}\in X$, there exists $y\in Y$ such that,
\begin{equation}\label{gra0}
\bar{y}=\arg\max_{y\in Y}F(\bar{x}, y)=\arg\max_{y\in \mathbb{R}^m}F(\bar{x}, y).
\end{equation}
Hence, $\nabla_yF(\bar{x}, \bar{y})=0.$
We construct a saddle point of P1.
Since $\bar{x}=\arg\min_{x\in X}F(x, \bar{y})$, by the KKT conditions, there exists some multiplier $\lambda\in \mathbb{R}^{m-r}$ such that
\begin{eqnarray*}
&&\nabla f(\bar{x})+A^T\bar{y}+A^TB^T\lambda=0\\
&&BA\bar{x}=0.
\end{eqnarray*}
Let $x^*=\bar{x}$ and $y^*=\bar{y}+B^T\lambda$, we have
$$\nabla_xF(x^*, y^*)=0$$
and
\begin{eqnarray*}
\nabla_yF(x^*, y^*)&=&-Qy^*+v+Ax^*\\
&=&(-Q\bar{y}+v+A\bar{x})+Q(\bar{y}-y^*)\\
&=&\nabla_yF(\bar{x}, \bar{y})+Q(-B^T\lambda)\\
&=&0,
\end{eqnarray*}
where the last equality is because the  rows of $B$ spans $\mathrm{Null}(Q)$.
We finish the proof. \hfill $\blacksquare$

\vspace{0.3cm}
Note that this lemma implies the nonempty property of the $\Yc(\mub, \zb)$, where $f(x), y^\top Ax, Q, v$ correspond $\textstyle \sum_{i=1}^N f_i(\xb_i) + \frac{p}{2}\|\xb_i - \zb_i\|^2$, $\sum_{i=1}^N \yb_i^\top \Bb_i\xb_i - \yb^\top \qb/N$, $\rho \Lb^-$, $\Ab^\top \mub$, respectively.

\section{Proof of Lemma \ref{fact: potential fun L}} \label{appendix L bound}
Notice that the change of the $\Lc_{\rho}$ function after one iteration can be bounded as
\begin{align}
& -\Lc_\rho(\yb^{r+1},\mub^{r+1};\zb^{r+1}) + \Lc_\rho(\yb^r,\mub^r;\zb^r)\notag\\
&\leq \bigg( -\Lc_\rho(\yb^{r},\mub^{r+1};\zb^{r}) + \Lc_\rho(\yb^r,\mub^r;\zb^r) \bigg)
 + \bigg( -\Lc_\rho(\yb^{r+1},\mub^{r+1};\zb^{r}) + \Lc_\rho(\yb^r,\mub^{r+1};\zb^r) \bigg) \notag \\
&~~~ + \bigg(  -\Lc_\rho(\yb^{r+1},\mub^{r+1};\zb^{r+1}) + \Lc_\rho(\yb^{r+1},\mub^{r+1};\zb^r) \bigg).
\end{align}
Then, we will analyze this desired bound by following four steps below.
\begin{enumerate}[(a)]
	\item Bound of $-\Lc_\rho(\yb^{r},\mub^{r+1};\zb^{r}) +
\Lc_\rho(\yb^r,\mub^r;\zb^r)$:
	\begin{align}
	-\Lc_\rho(\yb^{r},\mub^{r+1};\zb^{r}) + \Lc_\rho(\yb^r,\mub^r;\zb^r)= \la \mub^\rpo
-\mub^r, \Ab \yb^r\ra  = \frac{1}{\alpha}\|\mub^\rpo - \mub^r\|^2, \label{proof of
fact: potential fun L 1}
	\end{align}
    where the last equality is obtained from \eqref{eqn: alg original form mu}.	

	\item Bound of $-\Lc_\rho(\yb^{r+1},\mub^{r+1};\zb^{r}) +
\Lc_\rho(\yb^r,\mub^{r+1};\zb^r)$:
	Since $\Lc_\rho(\yb,\mub^{r+1};\zb^r) - \frac{1}{2}\|\yb - \yb^r\|_{\rho \Lb^+
}^2$ is strongly concave with modulus $\rho $, and
	\begin{align}\label{eq: nabla y exact}
	\nabla_y \Lc_\rho(\yb^\rpo,\mub^{r+1};\zb^r)-  \rho \Lb^+  (\yb^\rpo - \yb^r)
=\zerob,
	\end{align} by the optimality condition of \eqref{eqn: alg original form y},
we have
	\begin{align}
	\Lc_\rho(\yb^r,\mub^{r+1};\zb^r) - \frac{1}{2}\|\yb^r - \yb^r\|_{\rho \Lb^+
}^2
	\leq \Lc_\rho(\yb^\rpo,\mub^{r+1};\zb^r)
- \frac{1}{2}\|\yb^\rpo - \yb^r\|_{\rho \Lb^+
}^2 -\frac{\rho}{2}\| \yb^\rpo - \yb^r \|^2.
	\label{proof of fact: potential fun L 2}
	\end{align}
	
	\item Bound of $-\Lc_\rho(\yb^{r+1},\mub^{r+1};\zb^{r+1}) +
\Lc_\rho(\yb^{r+1},\mub^{r+1};\zb^r)$:
	Note that
	\begin{align}
	&-\Lc_\rho(\yb^{r+1},\mub^{r+1};\zb^{r+1}) +
\Lc_\rho(\yb^{r+1},\mub^{r+1};\zb^r)
    \notag\\
	&=\sum_{i=1}^N \bigg( \phi_i(\yb_i^\rpo,\zb_i^r) -
\phi_i(\yb_i^\rpo,\zb_i^\rpo) \bigg) \notag \\
	&=\sum_{i=1}^N \bigg( f_i(\xb_i(\yb_i^\rpo,\zb_i^r))  +
\frac{p}{2}\|\xb_i(\yb_i^\rpo,\zb_i^r)-\zb_i^r\|^2
    + \la \yb_i^\rpo,
\Bb_i\xb_i(\yb_i^\rpo,\zb_i^r)\ra \bigg)
       \notag \\
       &~~~ - \sum_{i=1}^N \bigg( f_i(\xb_i(\yb_i^\rpo,\zb_i^\rpo))  +
\frac{p}{2}\|\xb_i(\yb_i^\rpo,\zb_i^\rpo)-\zb_i^\rpo\|^2
    + \la \yb_i^\rpo,
\Bb_i\xb_i(\yb_i^\rpo,\zb_i^\rpo)\ra \bigg) .
\end{align}
Since $\xb_i(\yb_i^{r+1}, \zb_i^r)$ is a minimizer of  \eqref{eq: dual ADMM x1}, we have the upper bound
\begin{align}
    &-\Lc_\rho(\yb^{r+1},\mub^{r+1};\zb^{r+1}) +
\Lc_\rho(\yb^{r+1},\mub^{r+1};\zb^r)
    \notag\\
	&\leq \sum_{i=1}^N \bigg( f_i(\xb_i(\yb_i^\rpo,\zb_i^\rpo))  +
\frac{p}{2}\|\xb_i(\yb_i^\rpo,\zb_i^\rpo)-\zb_i^r\|^2
         + \la \yb_i^\rpo,
\Bb_i\xb_i(\yb_i^\rpo,\zb_i^\rpo)\ra \bigg) \notag\\
	&~~~- \sum_{i=1}^N \bigg( f_i(\xb_i(\yb_i^\rpo,\zb_i^\rpo))  +
\frac{p}{2}\|\xb_i(\yb_i^\rpo,\zb_i^\rpo)-\zb_i^\rpo\|^2 + \la \yb_i^\rpo,
\Bb_i\xb_i(\yb_i^\rpo,\zb_i^\rpo)\ra \bigg) \notag \\
	& =\frac{p}{2}\sum_{i=1}^N \|\xb_i(\yb_i^\rpo,\zb_i^\rpo)-\zb_i^r\|^2 -  \frac{p}{2} \sum_{i=1}^N
\|\xb_i(\yb_i^\rpo,\zb_i^\rpo)-\zb_i^\rpo\|^2 \notag \\
	& =- \frac{p}{2}\sum_{i=1}^N \la \zb_i^\rpo - \zb_i^r,  \zb_i^\rpo + \zb_i^r -
2\xb_i(\yb_i^\rpo,\zb_i^\rpo) \ra.
	\label{proof of fact: potential fun L 3}
	\end{align}
    \end{enumerate}

	Summing \eqref{proof of fact: potential fun L 1}, \eqref{proof of fact:
potential fun L
2} and \eqref{proof of fact: potential fun L 3} leads to the desired result.
	\hfill $\blacksquare$

\section{Proof of Lemma \ref{fact: potential fun d}} \label{appendix d bound}
	For any $\bar \yb^r \in \mathcal{Y}(\mub^\rpo,\zb^r)$ and $\tilde \yb^r \in \mathcal{Y}(\mub^r,\zb^r)$,
	\begin{align}
	 d(\mub^{r+1};\zb^{r}) - d(\mub^r;\zb^r)
	&=
	\bigg [ \sum_{i=1}^N \bigg(\phi_i(\bar \yb_i^r,\zb_i^r) - \la
\bar \yb_i^r, \bb/N \ra \bigg)
        -\la \mub^\rpo, \Ab \bar \yb^r\ra  -\frac{\rho}{2}\|\Ab
\bar \yb^r
\|^2 \bigg] \notag \\
	&~~~- \bigg [ \sum_{i=1}^N \bigg(\phi_i(\tilde \yb_i^r,\zb_i^r) - \la
\tilde \yb_i^r, \bb/N \ra \bigg)
    -\la \mub^r, \Ab \tilde \yb^r\ra
\bigg.-\frac{\rho}{2}\|\Ab \tilde \yb^r \|^2
\bigg] \notag \\
	&\overset{\mathrm{(i)}}{\leq}
        \bigg [ \sum_{i=1}^N \bigg(\phi_i(\bar \yb_i^r,\zb_i^r) - \la \bar \yb_i^r, \bb/N \ra \bigg)
            -\la \mub^\rpo, \Ab \bar \yb^r\ra  -\frac{\rho}{2}\|\Ab \bar \yb^r \|^2 \bigg] \notag \\
	   &~~~- \bigg [ \sum_{i=1}^N \bigg(\phi_i(\bar \yb_i^r,\zb_i^r) - \la \bar \yb_i^r,
            \bb/N \ra \bigg)
            -\la \mub^r, \Ab \bar \yb^r\ra
        \bigg.-\frac{\rho}{2}\|\Ab \bar \yb^r \|^2
        \bigg] \notag \\
    & = - \la \mub^{r+1} -\mub^{r},  \Ab \bar \yb^r \ra \notag \\
	&=  - \alpha  \la \Ab
\yb^r,  \Ab \bar \yb^r \ra, \label{d2}
	\end{align}
    where $\mathrm{(i)}$ is due to the fact that $\tilde \yb^r$ is a maximizer of $\Lc_{\rho} (\yb, \mub^{r}; \zb^{r} )$, and the last equality is deduced by \eqref{eqn: alg original form mu}.

	On the other hand, we have
	\begin{align}
	&d(\mub^{r+1};\zb^{r+1}) - d(\mub^\rpo;\zb^r) \notag \\
	&=
	\bigg [ \sum_{i=1}^N \bigg(\phi_i(\bar \yb_i^{r+1},\zb_i^\rpo) - \la
\bar \yb_i^{r+1}, \bb/N \ra \bigg)
    -\la \mub^\rpo, \Ab \bar \yb^{r+1}\ra  -\frac{\rho}{2}\|\Ab
\bar \yb^{r+1} \|^2 \bigg] \notag \\
	&~~~- \bigg [ \sum_{i=1}^N \bigg(\phi_i(\bar \yb_i^r,\zb_i^r)
- \la
\bar \yb_i^r, \bb/N \ra \bigg)-\la \mub^\rpo, \Ab \bar \yb^r\ra
 -\frac{\rho}{2}\|\Ab \bar \yb^r
\|^2 \bigg] \notag\\
& \overset{\mathrm{(i)}}{\leq} \bigg [ \sum_{i=1}^N \bigg(\phi_i(\bar \yb_i^{r+1},\zb_i^\rpo) - \la
\bar \yb_i^{r+1}, \bb/N \ra \bigg)
    -\la \mub^\rpo, \Ab \bar \yb^{r+1}\ra  -\frac{\rho}{2}\|\Ab
\bar \yb^{r+1} \|^2 \bigg] \notag \\
	&~~~- \bigg [ \sum_{i=1}^N \bigg(\phi_i(\bar \yb_i^{r+1},\zb_i^r)
- \la
\bar \yb_i^{r+1}, \bb/N \ra \bigg)-\la \mub^\rpo, \Ab \bar \yb^{r+1}\ra
 -\frac{\rho}{2}\|\Ab \bar \yb^{r+1}
\|^2 \bigg] \notag\\
& = \sum_{i=1}^N \bigg(\phi_i(\bar \yb_i^\rpo,\zb_i^\rpo) -
\phi_i(\bar \yb_i^\rpo,\zb_i^r) \bigg) ,
\end{align}
where $\mathrm{(i)}$ is obtained from the fact that $\bar \yb^r$ is
a maximizer of $\Lc_{\rho} (\yb, \mub^{r+1}; \zb^{r} )$. Then according to \eqref{eq: psi}, we have
\begin{align}
    &d(\mub^{r+1};\zb^{r+1}) - d(\mub^\rpo;\zb^r) \notag \\
    &\leq \sum_{i=1}^N \bigg( f_i(\xb_i(\bar \yb_i^{r+1},\zb_i^\rpo))
            + \frac{p}{2}\|\xb_i(\bar \yb_i^{r+1},\zb_i^\rpo)-\zb_i^\rpo\|^2
            + \la \bar \yb_i^{r+1},\Bb_i\xb_i(\bar \yb_i^{r+1},\zb_i^\rpo)\ra \bigg)
            \notag\\
    	&~~~- \sum_{i=1}^N \bigg( f_i(\xb_i(\bar \yb_i^{r+1},\zb_i^r))
            + \frac{p}{2}\|\xb_i(\bar \yb_i^{r+1},\zb_i^r)-\zb_i^r\|^2
            + \la \bar \yb_i^{r+1}, \Bb_i\xb_i(\bar \yb_i^{r+1},\zb_i^r)\ra \bigg)
            \notag\\
    & \overset{\mathrm{(i)}}{\leq}
        \sum_{i=1}^N \bigg( f_i(\xb_i(\bar \yb_i^{r+1},\zb_i^r))
            + \frac{p}{2}\|\xb_i(\bar \yb_i^{r+1},\zb_i^r)-\zb_i^\rpo\|^2
                + \la \bar \yb_i^{r+1}, \Bb_i\xb_i(\bar \yb_i^{r+1},\zb_i^r)\ra \bigg) \notag\\
    	&~~~
        - \sum_{i=1}^N \bigg( f_i(\xb_i(\bar \yb_i^{r+1},\zb_i^r))
            + \frac{p}{2}\|\xb_i(\bar \yb_i^{r+1},\zb_i^r)-\zb_i^r\|^2
            + \la \bar \yb_i^{r+1}, \Bb_i\xb_i(\bar \yb_i^{r+1},\zb_i^r)\ra \bigg)
                 \notag\\
    & = \frac{p}{2} \sum_{i=1}^N \bigg(
            \|\xb_i(\bar \yb_i^{r+1},\zb_i^r)-\zb_i^\rpo\|^2
    	      - \|\xb_i(\bar \yb_i^{r+1},\zb_i^r)-\zb_i^r\|^2 \bigg )  \notag\\
    & =\frac{p}{2} \sum_{i=1}^N \la  \zb_i^\rpo -\zb^r,
                - 2\xb_i(\bar \yb_i^{r+1},\zb_i^r) + \zb_i^\rpo +\zb_i^r      \ra,
	\end{align}
where $\mathrm{(i)}$ is due to the fact that $\xb_i(\bar \yb_i^{r+1}, \zb_i^{r+1})$ is a minimizer of $\phi_i(\bar \yb_i^{r+1}, \zb_i^{r+1})$.
	
	\hfill $\blacksquare$
\section{Proof of Lemma \ref{Fact G+c}} \label{appendix G bound}
The proof of Lemma \ref{Fact G+c} is a direct corollary of the following two lemmas.

\begin{Lemma}\label{fact: potential func G}
We have
		\begin{align} \label{eq: Bound for G}
		&G(\xb^{r+1}, \yb^{r+1}, \mub^{r+1}, \zb^{r+1}) -G(\xb^r, \yb^r, \mub^r,
\zb^r)\notag \\
		&\leq -\frac{1}{\alpha}\|\mub^\rpo - \mub^r\|^2 - \sum_{i=1}^N
\frac{p+\gamma^-}{2} \|\xb_i^\rpo-\xb_i^r\|^2
         - \frac{\rho}{2}\|\yb^\rpo -
\yb^r\|^2_{\Lb^-}+ \sum_{i=1}^N \la \yb_i^\rpo - \yb_i^r, \Bb_i(\xb_i^r -
\xb_i^\rpo) \ra \notag
\\
		&~~~+ \sum_{i=1}^N (2\rho |\Nc_i|)\| \yb_i^\rpo - \yb_i^r\|^2 +
\frac{p}{2}\bigg(1-\frac{2}{\beta}\bigg)\|\zb^\rpo - \zb^r\|^2.
		\end{align}
	\end{Lemma}

\begin{Lemma}\label{fact: potential fun cross}
		 \begin{align}\label{eq: potential fun cross}
		& \bigg(\frac{\alpha}{2}\|\yb^\rpo\|^2_{\Lb^-} +\frac{1}{2}\|\yb^\rpo- \yb^r\|^2_{\rho \Lb^+}\bigg)
             - \bigg(\frac{\alpha}{2}\|\yb^r \|^2_{\Lb^-} +\frac{1}{2}\|\yb^r- \yb^{r-1}\|^2_{\rho \Lb^+}\bigg)
                \notag \\
        & \leq -\sum_{i=1}^N \la \yb_i^\rpo - \yb_i^r, \Bb_i(\xb_i^r - \xb_i^\rpo) \ra
                  -\frac{1}{2}\|(\yb^r-\yb^\rmo) - (\yb^\rpo-\yb^r) \|^2_{\rho \Lb^+ }
           -\bigg( \rho-\frac{\alpha}{2} \bigg) \|\yb^\rpo - \yb^r\|^2_{\Lb^-}.
		\end{align}
	\end{Lemma}

We then give the proofs  of the above two lemmas.

{\bf Proof of Lemma \ref{fact: potential func G}:}
Notice that the change of the function $G$ after one iteration can be bounded as
\begin{align}
& G(\xb^{r+1}, \yb^{r+1}, \mub^{r+1}, \zb^{r+1}) - G(\xb^{r}, \yb^{r}, \mub^{r}, \zb^{r}) \notag \\
 & \leq
 \bigg( G(\xb^{r}, \yb^{r}, \mub^{r+1}, \zb^{r}) - G(\xb^{r},
\yb^{r}, \mub^{r}, \zb^{r})  \bigg) +  \bigg(  G(\xb^{r+1}, \yb^{r}, \mub^{r+1}, \zb^{r}) - G(\xb^{r}, \yb^{r},
\mub^{r+1}, \zb^{r})  \bigg) \notag \\
&~~~+ \bigg( G(\xb^{r+1}, \yb^{r+1}, \mub^{r+1}, \zb^{r}) - G(\xb^{r+1}, \yb^{r},
\mub^{r+1}, \zb^{r}) \bigg) \notag \\
&~~~+ \bigg( G(\xb^{r+1}, \yb^{r+1}, \mub^{r+1}, \zb^{r+1}) - G(\xb^{r+1}, \yb^{r+1},
\mub^{r+1}, \zb^{r}) \bigg).
\end{align}

Next, we will bound four terms in the right hand side (RHS) of the above inequality respectively.
\begin{enumerate}[(a)]
	\item Bound of $G(\xb^{r}, \yb^{r}, \mub^{r+1}, \zb^{r}) - G(\xb^{r},
\yb^{r}, \mub^{r},
\zb^{r}) $:
	\begin{align}
	G(\xb^{r}, \yb^{r}, \mub^{r+1}, \zb^{r}) - G(\xb^{r}, \yb^{r}, \mub^{r},
\zb^{r})=
-\la \mub^\rpo -\mub^r, \Ab \yb^r\ra = -\frac{1}{\alpha}\|\mub^\rpo - \mub^r\|^2, \label{proof of fact: potential
func G 1}
	\end{align}
    where the last equality is due to \eqref{eqn: alg original form mu}.
	
    \item

Bound of  $G(\xb^{r+1}, \yb^{r}, \mub^{r+1}, \zb^{r}) - G(\xb^{r}, \yb^{r},
\mub^{r+1}, \zb^{r}) $:

	Firstly, by the optimality condition of \eqref{eq: psi} , we
have
	\begin{align}\label{eq: optimality of udpate of x}
	&\la  \nabla f_i(\xb_i^\rpo)  +  p(\xb_i^\rpo-\zb_i^r) + \Bb_i^\top\yb_i^\rpo,
\xb_i
-\xb_i^\rpo \ra \geq 0 ,~\forall~ \xb_i.
	\end{align}
	Secondly, since $f_i(\xb_i)  +
\frac{p}{2}\|\xb_i-\zb_i^r\|^2$
is strongly convex with modulus $p+\gamma^-$, we  obtain
	\begin{align}
	&f_i(\xb_i^r)  +  \frac{p}{2}\|\xb_i^r-\zb_i^r\|^2\notag \\
    &\geq f_i(\xb_i^\rpo)  +
\frac{p}{2}\|\xb_i^\rpo-\zb_i^r\|^2 + \la \nabla f_i(\xb_i^\rpo)  +  p(\xb_i^\rpo-\zb_i^r), \xb_i^r -
\xb_i^\rpo \ra + \frac{p+\gamma^-}{2} \|\xb_i^r - \xb_i^\rpo\|^2 \notag \\
	& \geq f_i(\xb_i^\rpo)  +  \frac{p}{2}\|\xb_i^\rpo-\zb_i^r\|^2
+ \la - \Bb_i^\top\yb_i^\rpo, \xb_i^r - \xb_i^\rpo \ra
+ \frac{p+\gamma^-}{2}
\|\xb_i^r - \xb_i^\rpo\|^2, \label{proof of fact: potential func G 2}
	\end{align}
    where the last inequality is due to \eqref{eq: optimality of udpate of x}.

	Thus, we have the upper bound for the difference $G(\xb^{r+1}, \yb^{r},
\mub^{r+1}, \zb^{r}) - G(\xb^{r}, \yb^{r},
\mub^{r+1}, \zb^{r}) $ as
	\begin{align}
	&G(\xb^{r+1}, \yb^{r}, \mub^{r+1}, \zb^{r}) - G(\xb^{r}, \yb^{r}, \mub^{r+1},
\zb^{r})
\notag\\
& = \sum_{i=1}^N \bigg( f_i(\xb_i^\rpo)  +
\frac{p}{2}\|\xb_i^\rpo-\zb_i^r\|^2  -
f_i(\xb_i^r) -  \frac{p}{2}\|\xb_i^r-\zb_i^r\|^2 \bigg)
+  \sum_{i=1}^N  \la  \yb_i^r, \Bb_i ( \xb_i^\rpo -\xb_i^r)  \ra   \notag \\
	&  \leq \sum_{i=1}^N\la  \yb_i^\rpo - \yb_i^r, \Bb_i(\xb_i^r - \xb_i^\rpo) \ra- \frac{p+\gamma^-}{2} \|\xb_i^r - \xb_i^\rpo\|^2,
\label{proof of fact: potential func G 3}
	\end{align}
	where the last equality is obtained by \eqref{proof of fact: potential func G
2}.

	\item  Bound of $G(\xb^{r+1},
\yb^{r+1}, \mub^{r+1}, \zb^{r}) - G(\xb^{r+1}, \yb^{r},
\mub^{r+1}, \zb^{r}) $:
	
Since
	\begin{align}
	 \frac{\rho}{2}\|\Ab \yb^\rpo \|^2
    & \geq \frac{\rho}{2}\|\Ab \yb^r \|^2 +\rho\la \Lb^-
\yb^r, \yb^\rpo - \yb^r \ra
    +\frac{\rho}{2}\|\yb^\rpo - \yb^r\|^2_{\Lb^-},
	\end{align}
    we can have
	\begin{align}
	&G(\xb^{r+1}, \yb^{r+1}, \mub^{r+1}, \zb^{r}) - G(\xb^{r+1}, \yb^{r},
\mub^{r+1},
\zb^{r}) \notag \\
	&= \sum_{i=1}^N \la \yb_i^\rpo - \yb_i^r, \Bb_i\xb_i - \frac{\qb}{N} \ra
	-\la \mub^\rpo, \Ab ( \yb^\rpo  - \yb^r ) \ra-\frac{\rho}{2}\|\Ab \yb^\rpo \|^2 +
\frac{\rho}{2}\|\Ab \yb^r \|^2 \notag \\
	& = \sum_{i=1}^N \la \yb_i^\rpo - \yb_i^r, \Bb_i\xb_i - \frac{\qb}{N}   - \Ab_i^\top
\mub^\rpo \ra -\frac{\rho}{2}\|\Ab \yb^\rpo \|^2 + \frac{\rho}{2}\|\Ab \yb^r \|^2
	\notag \\
&\leq \sum_{i=1}^N \la \yb_i^\rpo - \yb_i^r, \Bb_i\xb_i - \frac{\qb}{N}   - \Ab_i^\top
\mub^\rpo \ra  - \rho\la \Lb^- \yb^r, \yb^\rpo - \yb^r \ra
-\frac{\rho}{2}\|\yb^\rpo - \yb^r\|^2_{\Lb^-}.
	\notag
	\end{align}
Using the fact that $\Lb^+ = 2\Db - \Lb^-$, we further have
	\begin{align}
&G(\xb^{r+1}, \yb^{r+1}, \mub^{r+1}, \zb^{r}) - G(\xb^{r+1}, \yb^{r},
\mub^{r+1},
\zb^{r}) \notag \\
	& \leq \sum_{i=1}^N \la \yb_i^\rpo - \yb_i^r, \Bb_i\xb_i - \frac{\qb}{N}  - \Ab_i^\top
\mub^\rpo  +
\rho \Lb_i^+ \yb^r\ra- 2\rho\la \Db   \yb^r, \yb^\rpo - \yb^r \ra -
	\frac{\rho}{2}\|\yb^\rpo - \yb^r\|^2_{\Lb^-}     \notag \\
	&\overset{\mathrm{(i)}}{=} \sum_{i=1}^N \la \yb_i^\rpo - \yb_i^r, 2\rho |\Nc_i| \yb_i^\rpo\ra
	- 2\rho\la \Db   \yb^r, \yb^\rpo - \yb^r \ra-
	\frac{\rho}{2}\|\yb^\rpo - \yb^r\|^2_{\Lb^-}\notag \\
	&= \sum_{i=1}^N (2\rho |\Nc_i|)\| \yb_i^\rpo - \yb_i^r\|^2 -
\frac{\rho}{2}\|\yb^\rpo -
\yb^r\|^2_{\Lb^-}, \label{proof of fact: potential func G 4}
	\end{align}
	where $\mathrm{(i)}$ is obtained by \eqref{y: update}.
	
	\item  Bound of $G(\xb^{r+1},
\yb^{r+1}, \mub^{r+1}, \zb^{r+1}) - G(\xb^{r+1}, \yb^{r+1},
\mub^{r+1}, \zb^{r}) $:
	It can be obtained as
	\begin{align}
	&G(\xb^{r+1}, \yb^{r+1}, \mub^{r+1}, \zb^{r+1}) - G(\xb^{r+1}, \yb^{r+1},
\mub^{r+1},
\zb^{r})\notag\\& = \frac{p}{2}\|\xb^\rpo-\zb^r\|^2 -\frac{p}{2}\|\xb^\rpo-\zb^r\|^2
\notag \\
	& = \frac{p}{2} \la -\zb^\rpo + \zb^r, 2(\xb^\rpo - \zb^r) - \zb^\rpo + \zb^r
\ra \notag
\\
	& \overset{\mathrm{(i)}}{=}   \frac{p}{2} \la -\zb^\rpo + \zb^r, \frac{2}{\beta}(\zb^\rpo -\zb^r) -
\zb^\rpo +
\zb^r \ra \notag \\
	& = \frac{p}{2}\bigg( 1 - \frac{2}{\beta} \bigg) \|\zb^\rpo - \zb^r\|^2,
	\label{proof of fact: potential func G 5}
	\end{align}
where $\mathrm{(i)}$ is due to the updating step \eqref{eq: dual ADMM lambda 22} for $\zb$ of algorithm 1.

	Summing \eqref{proof of fact: potential func G
1}, \eqref{proof of fact: potential func G 3}, \eqref{proof of fact: potential
func
G 4}, and \eqref{proof of fact: potential func G 5} gives the desired bound in \eqref{eq: Bound for G}. \hfill $\blacksquare$
    \end{enumerate}

{\bf{Proof of Lemma \ref{fact: potential fun cross}:}}
Firstly, the optimality condition of \eqref{eqn: alg original form y} for
iteration $r$ and $r-1$ are given as
    \begin{align}\label{opt-L}
    &\bigg  \langle  \nabla_{\yb}\left( \Lc_{\rho} (\yb^{r+1}, \mub^{r+1}; \zb^r) -
    \frac{\rho}{2} \|\yb^{r+1} - \yb^r\|_{\Lb^+}^2\right), \yb - \yb^{r+1} \bigg  \rangle  \leq 0, \forall \yb,
    \\
    & \bigg  \langle\nabla_{\yb}\left( \Lc_{\rho} (\yb^{r}, \mub^{r+1}; \zb^r) -
    \frac{\rho}{2} \|\yb^{r} - \yb^{r-1}\|_{\Lb^+}^2\right), \yb - \yb^{r} \bigg  \rangle\leq 0, \forall \yb.
    \end{align}
    According to Danskin's theorem \cite[Proposition B.22]{BK:Bertsekas2003_NP}, we have
    \begin{align} \label{gradient of phi}
    \nabla_{\yb_i} \phi_i(\yb_i^{r+1}, \zb_i^r) = \Bb_i \xb_i^{r+1}.
    \end{align}
Hence, the gradient of $\Lc_{\rho}(\yb, \mub; \zb)$ is given by
$$
\nabla_{\yb_i}\Lc_{\rho} (\yb^{r+1}, \mub^{r+1}; \zb^{r})=\Bb_i \xb_i^{r+1} - \frac{\qb}{N}   - \Ab_i^\top \mub^{r+1} -\rho \Lb_i^-
\yb_i^{r+1}.
$$
Substituting the above equality into \eqref{opt-L}, we have
	\begin{align}
	&\la \Bb_i \xb_i^\rpo - \bb/N  - \Ab_i^\top \mub^\rpo  -\rho \Lb_i^- \yb_i^\rpo  -
\rho
\Lb_i^+ (\yb_i^\rpo -\yb_i^r), \yb_i -\yb_i^\rpo\ra
\leq 0,~\forall ~\yb_i,     \\
	&\la \Bb_i \xb_i^r - \bb/N  - \Ab_i^\top \mub^r  -\rho \Lb_i^- \yb_i^r  - \rho \Lb_i^+(\yb_i^r -\yb_i^\rmo),\yb_i -\yb_i^r\ra
\leq 0,~\forall ~\yb_i.  \label{eq: proof of fact: potential fun cross 1}
	\end{align}
	Letting $\yb_i=\yb_i^r$ and $\yb_i=\yb_i^\rpo$ in the first and second equations above
respectively,
and summing up these two inequalities, we obtain
	\begin{align}
	& \la \Bb_i (\xb_i^r - \xb_i^\rpo), \yb_i^\rpo -\yb_i^r\ra - \la \Ab_i^\top
(\mub^r
-\mub^\rpo ), \yb_i^\rpo -\yb_i^r\ra
	- \la \rho \Lb_i^- (\yb_i^r - \yb_i^\rpo), \yb_i^\rpo -\yb_i^r\ra
    \notag \\
        & - \la  (\rho \Lb_i^+ )(
(\yb_i^r -\yb_i^\rmo) - (\yb_i^\rpo-\yb_i^r) ),~ \yb_i^\rpo -\yb_i^r\ra
	\leq 0, \label{eq: proof of fact: potential fun cross 2}
	\end{align} for $\forall i \in [N]$. By summing the $N$ equations, one obtains
	\begin{align}
	\sum_{i=1}^N\la \Bb_i (\xb_i^r - \xb_i^\rpo), \yb_i^\rpo -\yb_i^r\ra
& \leq \la  \mub^r -\mub^\rpo , \Ab (\yb^\rpo -\yb^r)\ra - \rho \|\yb^\rpo -
\yb^r\|_{\Lb^-}^2 \notag \\
	&~~~+\la   (\yb^r -\yb^\rmo) - (\yb^\rpo-\yb^r) ,  \rho \Lb^+  (\yb^\rpo
-\yb^r)\ra.
	\label{eq: proof of fact: potential fun cross 3}
	\end{align}
	Notice that
	\begin{align}
	 \la  \mub^r -\mub^\rpo , \Ab (\yb^\rpo -\yb^r)\ra
&\overset{\mathrm{(i)}}{=} \la  \alpha \Ab \yb^r , \Ab (\yb^r
-\yb^\rpo)\ra \notag  \\
	&=\frac{\alpha}{2}\|\yb^\rpo - \yb^r\|_{\Lb^-}^2 + \frac{\alpha}{2}\|
\yb^r\|_{\Lb^-}^2
- \frac{\alpha}{2}\| \yb^\rpo \|_{\Lb^-}^2,
	\label{eq: proof of fact: potential fun cross 4}
	\end{align}
in which $\mathrm{(i)}$ is obtained from \eqref{eqn: alg original form mu}. Moreover, we have
	\begin{align}
	&\la   (\yb^r -\yb^\rmo) - (\yb^\rpo-\yb^r) ,  \rho \Lb^+   (\yb^\rpo -\yb^r)\ra
\notag\\
	&=  \frac{1}{2}\|\yb^r -\yb^\rmo \|^2_{\rho \Lb^+  } -  \frac{1}{2}\|(\yb^r
-\yb^\rmo)-(\yb^\rpo -\yb^r) \|^2_{\rho \Lb^+ }- \frac{1}{2}\|\yb^\rpo -\yb^r \|^2_{\rho \Lb^+ }.
	\label{eq: proof of fact: potential fun cross 5}
	\end{align}
	By substituting \eqref{eq: proof of fact: potential fun cross 4} and \eqref{eq:
proof of
fact: potential fun cross 5} into \eqref{eq: proof of fact: potential fun cross 3},
we
obtain
    \begin{align}
		& \sum_{i=1}^N \la \yb_i^\rpo - \yb_i^r, \Bb_i(\xb_i^r - \xb_i^\rpo) \ra
\notag \\
		& \leq  \frac{\alpha}{2}\|\yb^r\|^2_{\Lb^-} -\frac{\alpha}{2} \|\yb^\rpo\|^2_{\Lb^-} + \frac{1}{2}\|\yb^r -
\yb^\rmo\|^2_{\rho \Lb^+}-\frac{1}{2}\|\yb^\rpo - \yb^r\|^2_{\rho \Lb^+}
                -\frac{1}{2}\|(\yb^r-\yb^\rmo) - (\yb^\rpo-\yb^r) \|^2_{\rho \Lb^+ }
		\notag \\
		&~~~-\bigg( \rho-\frac{\alpha}{2} \bigg) \|\yb^\rpo - \yb^r\|^2_{\Lb^-}.
		\end{align}

    After rearranging the inequality above, we have
    \begin{align}
		& \bigg(\frac{\alpha}{2}\|\yb^\rpo\|^2_{\Lb^-} +\frac{1}{2}\|\yb^\rpo- \yb^r\|^2_{\rho \Lb^+}\bigg)
            - \bigg(\frac{\alpha}{2}\|\yb^r \|^2_{\Lb^-} +\frac{1}{2}\|\yb^r- \yb^{r-1}\|^2_{\rho \Lb^+}\bigg)
                \notag \\
        & \leq -\sum_{i=1}^N \la \yb_i^\rpo - \yb_i^r, \Bb_i(\xb_i^r - \xb_i^\rpo) \ra
            -\frac{1}{2}\|(\yb^r-\yb^\rmo) - (\yb^\rpo-\yb^r) \|^2_{\rho \Lb^+ }
            -\bigg( \rho-\frac{\alpha}{2} \bigg) \|\yb^\rpo - \yb^r\|^2_{\Lb^-}.
		\end{align}
	\hfill $\blacksquare$

\section{Proof of Lemma \ref{fact: pertubration of x_i}} \label{appendix:
perturbation of x}
\subsection{Proof of \eqref{eq: perturb of xi} in Lemma \ref{fact: pertubration of x_i}}
The proof of \eqref{eq: perturb of xi} is similar to \cite[Lemma 3.6]{Zhang18}.
Recall from \eqref{eq: x(y, z)} that
	\begin{align}\label{eq: psi2}
	\xb_i(\yb_i,\zb_i) = \arg \min_{\xb_i} \bigg\{f_i(\xb_i)  +
\frac{p}{2}\|\xb_i-\zb_i\|^2 + \yb_i^\top\Bb_i\xb_i\bigg\},
	\end{align}
	which is a strongly convex problem with modulus $p+\gamma^-$.
Thus,
according to the error bound for strongly convex problems in \cite[Theorem 3.1]{pang1987posteriori},
we have
	\begin{align}\label{proof of fact: pertubration of x_i 1}
	     \| \xb_i -\xb_i(\yb_i,\zb_i)  \| \leq \frac{ (p + \gamma^+ + 1) }{p+\gamma^-} \| \nabla f_i(\xb_i)+ p(\xb_i - \zb_i) + \Bb_i^\top \yb_i \|,
	\end{align} for all $ \xb_i  \in \mathbb{R}^{n}$. Since
	\begin{align}\label{proof of fact: pertubration of x_i 2}
	     \nabla f_i(\xb_i(\yb_i',\zb_i')) +p( \xb_i(\yb_i',\zb_i')- \zb_i') + \Bb_i^\top \yb_i'=0,
	\end{align} by letting $\xb_i = \xb_i(\yb_i',\zb_i')$ in \eqref{proof of fact:
pertubration of x_i 1}, we obtain
		\begin{align}\label{proof of fact: pertubration of x_i 3}
	   \| \xb_i(\yb_i',\zb_i')-\xb_i(\yb_i,\zb_i)  \|
 \leq \frac{ (p + \gamma^+ + 1) }{p+\gamma^-} \|\nabla f_i(\xb_i(\yb_i',\zb_i')) +
p(\xb_i(\yb_i',\zb_i')- \zb_i) + \Bb_i^\top \yb_i \|.
\end{align}

By substituting \eqref{proof of fact: pertubration of x_i 2} into \eqref{proof of fact: pertubration of x_i 3}, we have
\begin{align}
	\| \xb_i(\yb_i',\zb_i')-\xb_i(\yb_i,\zb_i)  \|& \leq \frac{ (p + \gamma^+ + 1) }{p+\gamma^-} \|  p(\zb_i -\zb_i') +  \Bb_i^\top (\yb_i -\yb_i')\| \notag \\
	& \overset{\mathrm{(i)}}\leq \frac{p (p + \gamma^+ + 1) }{p+\gamma^-} \| \zb_i -\zb_i' \| +   \frac{\|\Bb_i\| (p + \gamma^+ + 1) }{p+\gamma^-} \|
\yb_i -\yb_i'\|\notag\\
& \leq \frac{p (p + \gamma^+ + 1) }{p+\gamma^-} \| \zb_i -\zb_i' \| +   \frac{B_{\max} (p + \gamma^+ + 1) }{p+\gamma^-} \|
\yb_i -\yb_i'\|,
	\end{align}
where $\mathrm{(i)}$ is obtained by the triangle inequality and Cauchy-Schwartz inequality. \hfill $\blacksquare$

\vspace{0.2cm}

\subsection{Proof of \eqref{sigma1} in Lemma \ref{fact: pertubration of x_i}}
To show \eqref{sigma1}, let $$ \bar{h}(\xb; \mub, \zb)=  \max_{\yb}\bigg\{G(\xb,
\yb, \mub, \zb) - \frac{p}{2}\|\zb\|^2 \bigg\}.$$
By \eqref{eq: x(y, z)}, $\xb(\bar \yb  , \zb)$ is the  minimizer of the function
$\bar{h}(\xb;\mub, \zb)$, where $\bar \yb \in \mathcal{Y}(\mub, \zb)$. Then for any
$\bar \yb^r \in \mathcal{Y}(\mub^{r+1}, \zb^r)$ and $\bar \yb^{r+1} \in
\mathcal{Y}(\mub^{r+1}, \zb^{r+1})$, we can have $\nabla \bar{h}(\xb(\bar\yb^r ,
\zb^r); \mub^{r+1}, \zb^{r}) = \zerob $ and $\nabla  \bar{h}(\xb(\bar \yb^{r+1} ,
\zb^{r+1});\mub^{r+1}, \zb^{r+1}) = \zerob$.
Due to the fact that $\bar{h}$ is strongly convex with modulus $(p+\gamma^-)$, we
have
\begin{align}\label{eq: pertub of x strongly convex 1}
& \bar{h}(\xb(\bar \yb^{r+1}, \zb^{r+1});\mub^{r+1}, \zb^r)-\bar{h}(\xb(\bar\yb^r,
\zb^r); \mub^{r+1}, \zb^r)  \ge \frac{p+\gamma^-}{2}\|\xb(\bar\yb^r, \zb^r)-\xb(\bar \yb^{r+1} ,
\zb^{r+1})\|^2,
\end{align}
and
\begin{align}\label{eq: pertub of x strongly convex 2}
&\bar{h}(\xb(\bar\yb^r, \zb^r); \mub^{r+1}, \zb^{r+1})-\bar{h}(\xb(\bar \yb^{r+1},
\zb^{r+1}); \mub^{r+1}, \zb^{r+1}) \ge \frac{p+\gamma^-}{2}\|\xb(\bar\yb^r, \zb^r)-\xb(\bar \yb^{r+1},
\zb^{r+1})\|^2.
\end{align}
Then, by summing \eqref{eq: pertub of x strongly convex 1} and \eqref{eq: pertub of
x strongly convex 2}, we obtain
\begin{align}\label{eq: pertub of x h1}
& \bar{h}(\xb(\bar \yb^{r+1}, \zb^{r+1}); \mub^{r+1}, \zb^r)-\bar{h}(\xb(\bar\yb^r,
\zb^r);  \mub^{r+1}, \zb^r) \notag \\
&  + \bar{h}(\xb(\bar\yb^r, \zb^r);  \mub^{r+1}, \zb^{r+1})-\bar{h}(\xb(\bar
\yb^{r+1}, \zb^{r+1});  \mub^{r+1}, \zb^{r+1})\notag \\
& \ge (p+\gamma^- )\|\xb(\bar\yb^r, \zb^r)-\xb(\bar \yb^{r+1}, \zb^{r+1})\|^2.
\end{align}

On the other hand, according to the definition of $\bar{h}$, we have
\begin{align}\label{eq: pertub of x h11}
&\bar{h}(\xb(\bar\yb^r, \zb^r); \mub^{r+1}, \zb^r)-\bar{h}(\xb(\bar\yb^r, \zb^r);
\mub^{r+1}, \zb^{r+1})
\notag \\
& = \frac{p}{2}\|\xb(\bar\yb^r, \zb^r) - \zb^r\|^2 -  \frac{p}{2}\|\zb^r\|^2
     - \bigg(\frac{p}{2}\|\xb(\bar\yb^r, \zb^r)- \zb^{r+1}\|^2 -
    \frac{p}{2}\|\zb^{r+1}\|^2 \bigg)
\notag \\
& =p\xb(\bar\yb^r, \zb^r)^\top(\zb^{r+1}-\zb^r),
\end{align}
and
\begin{align}\label{eq: pertub of x h21}
&\bar{h}(\xb(\bar \yb^{r+1}, \zb^{r+1}); \mub^{r+1}, \zb^r)-\bar{h}(\xb(\bar
\yb^{r+1}, \zb^{r+1}); \mub^{r+1}, \zb^{r+1})
  =p\xb(\bar \yb^{r+1}, \zb^{r+1})^\top(\zb^{r+1}-\zb^r).
\end{align}
By taking the difference between \eqref{eq: pertub of x h11} and \eqref{eq: pertub
of x h21}, and substituting it into \eqref{eq: pertub of x h1}, we obtain
\begin{align}
&(p+\gamma^-)\|\xb(\bar\yb^r, \zb^r)-\xb(\bar \yb^{r+1}, \zb^{r+1})\|^2
\notag \\
& \le p(\xb(\bar\yb^r, \zb^r)-\xb(\bar \yb^{r+1},
\zb^{r+1}))^\top(\zb^r-\zb^{r+1})\notag\\
& \leq p \|\xb(\bar\yb^r, \zb^r)-\xb(\bar \yb^{r+1},
\zb^{r+1})\|\|\zb^r-\zb^{r+1}\|.
\end{align}
Dividing $(p+\gamma^-)\|\xb(\bar\yb^r, \zb^r)-\xb(\bar \yb^{r+1}, \zb^{r+1})\|$ from
both sides of the above inequality, we have
\begin{align}
  \|\xb(\bar\yb^r, \zb^r)-\xb(\bar \yb^{r+1}, \zb^{r+1})\|
  \le \frac{p}{p+\gamma^-}\|\zb^r-\zb^{r+1}\|
\end{align} which is \eqref{sigma1}.
		\hfill $\blacksquare$

\section{Proof of Lemma \ref{fact: perturb of y 2}} \label{appendix: perturbation of
y 2}

To prove the dual error bounds, we recap the Hoffman error bound \cite[Lemma 3.2.3]{facchinei2007finite} given below.

\begin{Lemma}\label{Hoffman}
Let $\Ab \in \mathbb{R}^{m \times n}$, $\Cb \in \mathbb{R}^{k \times n}$, $\bb \in
\mathbb{R}^{m}$, $\db \in \mathbb{R}^{k }$, and the linear set
\begin{align*}
\Sc = \{\xb | \Ab \xb \leq \bb, \Cb \xb = \db\}.
\end{align*}
The distance from a point $\bar{\xb}
\in \mathbb{R}^{n}$ to $\Sc$ is bounded by:
\begin{align*}
\mathrm{dist}^2(\bar{\xb}, \Sc) \leq \theta^2 (\|(\Ab\bar{\xb} - \bb)_+\|^2+
\|\Cb\bar{\xb} - \db\|^2),
\end{align*}
where $(\cdot)_+=\max\{\cdot,0\}$ denotes the projection to the nonnegative orthant
and $\theta$ is
a positive constant depending on $\Ab$ and $\Cb$ only.
\end{Lemma}

To prove Lemma \ref{fact: perturb of y 2}, let us define $F_i(\xb_i) = f_i (\xb_i) +
\frac{p}{2} \|\xb_i -
    \zb_i^r\|^2$  for all $i\in [N]$,  and thus $\nabla F_i(\xb_i) = \nabla
    f_i(\xb_i) + p (\xb_i -
    \zb^r_i)$. Note that
    function $F_i$ is strongly convex with modulus $(p+\gamma^-)$ and its gradient
    $\nabla F_i$ is
    $(p + \gamma^+)$-Lipschitz continuous. Then by\cite[Eqn.
    (2.1.8)]{nesterov2013introductory}), we have
    \begin{align} \label{gradient dif bound of F}
      & \|\nabla F_i(\xb_i) - \nabla F_i(\xb'_i)\|^2
         \leq (p + \gamma^+)\la\nabla
        F_i(\xb_i) - \nabla F_i(\xb_i'), \xb_i- \xb_i' \ra.
    \end{align}
Besides,  since function $\phi_i$ in \eqref{eq: psi} is concave over $\yb$ and
smooth with Lipschitz constant $\textstyle \frac{1}{p+\gamma^-}$, we have
    \begin{align} \label{gradient dif bound of phi}
      & \|\nabla \phi_i(\xb_i) - \nabla \phi_i(\xb'_i)\|^2
        \leq -\frac{1}{p+\gamma^-}\la\nabla
        \phi_i(\xb_i) - \nabla \phi_i(\xb_i'), \xb_i- \xb_i' \ra.
    \end{align}

    Recall \eqref{eqn: dual function of DP2} and let
	\begin{align}\label{proof of fact: perturb of y2 1}
\mathcal{Y}(\mub^{r+1}, \zb^r)
 &= \arg \max_{\yb}~  \sum_{i=1}^N
\bigg(\phi_i(\yb_i,\zb_i^r) -
                    \frac{1}{N}\yb_i^\top\qb\bigg)  -\la \mub^\rpo, \Ab \yb\ra  -\frac{\rho}{2}\|\Ab \yb
                    \|^2.
	\end{align}
For $\bar \yb ^r \in \mathcal{Y}(\mub^{r+1}, \zb^r)$, we have the following KKT
conditions
	\begin{align}
	 & \begin{bmatrix}
	       \nabla \phi_1(\bar \yb_1^r , \zb_1^r) - \frac{\qb}{N} \\
	       \vdots \\
	        \nabla \phi_N(\bar \yb_N^r , \zb_N^r)  - \frac{\qb}{N}
	  \end{bmatrix}
	  -\Ab^\top \mub^{r+1} - \rho \Lb^- \bar\yb^r =\zerob,
        \label{KKT y2 C1} \\
	  &  \nabla  F_i ( \xb_i(\bar \yb_i^r ,\zb_i^{r}) )
	   + \Bb_i^\top \bar\yb_i^r =\zerob, \forall i\in[N],\label{KKT y2 C2}
	\end{align}
where  $\nabla \phi_i(\bar \yb_i^r, \zb_i^r) = \Bb_i \xb_i(\bar \yb_i^r, \zb_i^r)$
for all $i$ obtained by the Danskin's Theorem \cite[Proposition
B.22]{BK:Bertsekas2003_NP}.
Analogously, from \eqref{eqn: alg original form y}
	the KKT conditions for $\yb^{r+1}$ are
		\begin{align}
	& \begin{bmatrix}
	\nabla \phi_1(\yb_1^{r+1} , \zb_1^r) -  \frac{\qb}{N}  \\
	\vdots \\
	\nabla \phi_N(\yb_N^{r+1} , \zb_N^r) -  \frac{\qb}{N}
	\end{bmatrix}
	-\Ab^\top \mub^{r+1} - \rho \Lb^- \yb^{r+1} - \rho \Lb^+ (\yb^\rpo -\yb^r) =\zerob,
\label{KKT y2 r+1 C1} \\
	& \nabla  F_i ( \xb_i(\yb_i^{r+1},\zb_i^r) ) + \Bb_i^\top \yb_i^{r+1}
=\zerob, \forall i\in [N],
	 \label{KKT y2 r+1 C2}
	\end{align}
where {$\nabla \phi_i(\yb_i^{r+1} , \zb_i^r) = \Bb_i \xb_i(\yb_i^{r+1}, \zb_i^r)$
for all $i$.}
	
Let $\bar \yb^r$ be the projection point of  $\yb^{r+1}$ onto
$\Yc(\mub^{r+1},\zb^{r})$, i.e., ${\rm
dist}(\yb^{r+1},\Yc(\mub^{r+1},\zb^{r})) = \|\yb^{r+1}-\bar \yb^r\|$.
	With fixed $\{\xb_i\big(\bar \yb_i^r,\zb_i^{r}\big)\}_{i=1}^N$,
\eqref{KKT y2 C1}-\eqref{KKT y2 C2} forms a linear system of
	variable $\rho \bar \yb^r$, and any $\rho \bar \yb^r$ satisfying \eqref{KKT y2
C1}-\eqref{KKT y2 C2} is a point in $ \mathcal{Y}(\mub^{r+1}, \zb^r)$. By \eqref{KKT
y2 r+1 C1}-\eqref{KKT y2
r+1 C2} and the Hoffman bound in Lemma \ref{Hoffman}, we can bound ${\rm
dist}(\yb^{r+1},\Yc(\mub^{r+1},\zb^{r}))$ for some constant $\theta_1 >0$ as follows
	\begin{align}\label{bound1}
	    &\rho^2 {\rm dist}^2(\yb^{r+1},\Yc(\mub^{r+1},\zb^{r})) \notag \\
	    &\leq
            {\theta_1^2} \sum_{i=1}^N  \| \Bb_i ( \xb_i(\yb_i^{r+1},\zb_i^r) ) -
            \xb_i(\bar \yb_i^r,\zb_i^{r}) ) \|^2
            +
             {\theta_1^2} \rho^2 \sum_{i=1}^N  \|
            \nabla F_i (\xb_i(\yb_i^{r+1},\zb_i^r))- \nabla F_i ( \xb_i(\bar
            \yb_i^r,\zb_i^{r}) )
            \|^2
           \notag \\
            &~~~ + {\theta_1^2} \rho^2 \|\Lb^+(\yb^{r+1} - \yb^r)\|^2
	    \notag
	    \\
        & \overset{\mathrm{(i)}}{\leq}
           {-\frac{\theta_1^2 }{p + \gamma^-} } \sum_{i=1}^N  \la \Bb_i
            \xb_i(\yb_i^{r+1},\zb_i^r)- \Bb_i
            \xb_i(\bar \yb_i^r,\zb_i^{r})  , \yb_i^{r+1} -\bar \yb_i^r\ra
            \notag \\
             &
             ~~~
            +
                \theta_1^2 \rho^2 (p + \gamma^+) \sum_{i=1}^N
                \left\langle \nabla F_i (\xb_i(\yb_i^{r+1},\zb_i^r))\right. \left. - \nabla F_i ( \xb_i(\bar \yb_i^r,\zb_i^{r}) ),
\xb_i(\yb_i^{r+1},\zb_i^r)-
                    \xb_i(\bar \yb_i^r,\zb_i^{r}) \right\rangle
             \notag \\
             &
             ~~~+\theta_1^2 \rho^2 \|\Lb^+(\yb^{r+1} - \yb^r)\|^2\notag \\
	    &  \overset{\mathrm{(ii)}}{=}
        -\frac{\theta_1^2 }{p + \gamma^-} \sum_{i=1}^N  \la \Bb_i
        \xb_i(\yb_i^{r+1},\zb_i^r)- \Bb_i
        \xb_i(\bar \yb_i^r,\zb_i^{r})  , \yb_i^{r+1} - \bar \yb_i^r\ra
    	       \notag \\
             &
             ~~~
    -
            \theta_1^2 \rho^2 (p + \gamma^+)\sum_{i=1}^N
            \la \Bb_i^T \yb_i^{r+1}-\Bb_i^T \bar \yb_i^r,
             \xb_i(\yb_i^{r+1},\zb_i^r) -
\xb_i(\bar \yb_i^r,\zb_i^{r})
            \ra
            \notag \\
            &~~~+\theta_1^2 \rho^2 \|\Lb^+(\yb^{r+1} - \yb^r)\|^2 ,
	\end{align}
where the first term of $\mathrm{(i)}$ is due to \eqref{gradient dif bound of phi},
the third term of $\mathrm{(i)}$ is directly obtained by
\eqref{gradient dif bound of F}, and $\mathrm{(ii)}$ is owing to \eqref{KKT y2 C2}
and \eqref{KKT y2 r+1 C2}.
	
Rearranging the right hand side (RHS) of \eqref{bound1},
we further have the upper bound as
    \begin{align}
    &\rho^2 {\rm dist}^2(\yb^{r+1},\Yc(\mub^{r+1},\zb^{r})) \notag \\
    & \leq -\bigg( \frac{\theta_1^2 }{p + \gamma^-} +\theta_1^2 \rho^2 (p +
    \gamma^+) \bigg)
         \sum_{i=1}^N  \la \Bb_i
        \xb_i(\yb_i^{r+1},\zb_i^r)
         - \Bb_i
        \xb_i(\bar \yb_i^r,\zb_i^{r})  ,
                \yb_i^{r+1} - \bar \yb_i^r\ra
    	       \notag \\
            &~~~ +\theta_1^2 \rho^2 \|\Lb^+(\yb^{r+1} - \yb^r)\|^2
             \notag \\
	& \overset{\mathrm{(i)}} {=}  -\bigg( \frac{\theta_1^2 }{p + \gamma^-}
+\theta_1^2 \rho^2 (p + \gamma^+) \bigg)
    \sum_{i=1}^N
    \la \nabla \phi_i (\yb^\rpo, \zb^r)
     - \nabla \phi_i
    (\bar\yb_i^r,\zb_i^{r}),
    \yb_i^{r+1} -\bar \yb_i^r \ra  +\theta_1^2 \rho^2 \|\Lb^+(\yb^{r+1} - \yb^r)\|^2
    \notag\\
    & \overset{\mathrm{(ii)}} {=} - \bigg( \frac{\theta_1^2 }{p + \gamma^-}
    +\theta_1^2 \rho^2 (p + \gamma^+) \bigg)
    \la \nabla \Lc_{\rho} (\yb^\rpo, \mu^\rpo; \zb^r)
- \nabla \Lc_{\rho}
    (\bar \yb^r,\mu^\rpo;\zb^{r})+ \rho \Lb^- (\yb^{r+1} -
    \yb(\mub^{r+1},\zb^{r})),\notag \\
            &~~~~~~
    \yb^{r+1} - \bar \yb^r \ra+\theta_1^2 \rho^2 \|\Lb^+(\yb^{r+1} - \yb^r)\|^2
    \notag\\
    & = - \bigg( \frac{\theta_1^2 }{p + \gamma^-} +\theta_1^2 \rho^2 (p + \gamma^+)
    \bigg)
    \la \nabla \Lc_{\rho} (\yb^\rpo, \mu^\rpo; \zb^r)
     - \nabla \Lc_{\rho}
    (\bar \yb^r,\mu^\rpo;\zb^{r}),
    \yb^{r+1} -  \bar \yb^r \ra
    \notag \\
            &~~~
        - \bigg( \frac{\theta_1^2 }{p + \gamma^-} +\theta_1 \rho^2 (p + \gamma^+)
        \bigg)\rho
        \| \yb^{r+1} -  \bar \yb_i^r\|_{\Lb^-}^2
             +\theta_1^2 \rho^2 \|\Lb^+(\yb^{r+1} - \yb^r)\|^2,
    \label{due error bound y2 C3}
    \end{align}
	where $\mathrm{(i)}$ is due to the definition of $\nabla\phi_i$ , and
$\mathrm{(ii)}$ is obtained by substituting the gradient of $\Lc_\rho$ in
\eqref{eqn: AL of DP2}.

On the other hand, by the optimality conditions of \eqref{eqn: dual function of DP2}
for $\bar \yb^r$ and $\bar \yb^{r+1}$ respectively, we have
	\begin{align}
		&\nabla_\yb \Lc_\rho(\bar \yb^r,\mub^{r+1};\zb^{r})
=\zerob,\label{eqn: alg original form y2 1 opt} \\
		&\nabla_\yb \Lc_\rho(\yb^{r+1} ,\mub^{r+1};\zb^r) - \rho \Lb^+ (\yb^{r+1} -
\yb^r) =\zerob.\label{eqn: alg original form y2 2 opt}
	\end{align}
	Thus,
	\begin{align}
	  &\la \yb^{r+1}-\bar \yb^r,  \nabla_\yb
\Lc_\rho(\bar \yb^r,\mub^{r+1};\zb^{r})
    - \nabla_\yb \Lc_\rho(\bar \yb^r;\zb^r)\ra
    = \la \yb^{r+1}-\bar \yb^r,  - \rho \Lb^+ (\yb^{r+1} -
\yb^r)\ra.
	   \label{due error bound y2 C6}
	\end{align}
Substituting \eqref{due error bound y2 C6} into \eqref{due error bound y2 C3}, we
have
    \begin{align} \label{y^r+1 - y(mu^r+1, z^r) 1}
        &\rho^2 \|\yb^{r+1} -  \bar \yb^r\|^2
        \notag \\
        &  \leq - \bigg( \frac{\theta_1^2 }{p + \gamma^-} +\theta_1^2 \rho^2 (p +
        \gamma^+) \bigg) \la \yb^{r+1}-\bar \yb^r,
        - \rho \Lb^+ (\yb^{r+1} -\yb^r)\ra
            \notag \\
            &~~~
            - \bigg( \frac{\theta_1^2 }{p + \gamma^-} +\theta_1^2 \rho^2 (p +
            \gamma^+) \bigg)\rho
            \| \yb^{r+1} - \bar \yb^r\|_{\Lb^-}^2
              +\theta_1^2 \rho^2 \|\Lb^+(\yb^{r+1} - \yb^r)\|^2 \notag\\
        & \leq  \frac{\theta_1^2}{2\eta_1} \bigg( \frac{1 }{p + \gamma^-} + \rho^2
        (p + \gamma^+) \bigg)
        \| \yb^{r+1} -  \bar \yb^r\|^2
        \notag \\
            &~~~ \!+\! \bigg[\frac{\eta_1}{2}  \bigg( \frac{1 }{p + \gamma^-} \!+\!
            \rho^2 (p + \gamma^+) \bigg)\!+\!1\bigg]\theta_1^2\rho^2
            \|\Lb^+(\yb^{r+1} - \yb^r)\|^2 \notag \\
            &~~~- \bigg( \frac{\theta_1^2 }{p + \gamma^-} +\theta_1^2 \rho^2 (p +
            \gamma^+) \bigg)\rho
            \| \yb^{r+1} -  \bar \yb^r\|_{\Lb^-}^2,
    \end{align}
    where the last inequality is by Young's inequality and $\eta_1 >0$.

    Choose {$\eta_1 = \textstyle \frac{\theta_1^2 }{\rho^2}\big( \frac{1 }{p +
    \gamma^-} + \rho^2 (p + \gamma^+) \big) $}, then
    $\frac{\theta_1^2 }{2\eta_1}\big( \frac{1 }{p + \gamma^-} + \rho^2 (p +
    \gamma^+) \big) =
\frac{\rho^2}{2}
        < \rho^2$. Then, we can obtain \eqref{eq: y r+1 - y mu r+1 z r L-} from
        \eqref{y^r+1 - y(mu^r+1, z^r) 1}.

In addition, by substituting the chosen $\eta_1$ into \eqref{y^r+1 - y(mu^r+1, z^r)
1}, we further obtain the desired result in \eqref{y^r+1 - y(mu^r+1, z^r) 2}.
	\hfill $\blacksquare$
\section{Proof of Lemma \ref{fact: perturb of y}} \label{appendix: perturbation of y}
For any $\bar \yb ^r \in \mathcal{Y}(\mub^{r+1}, \zb^r)$, we have the KKT condition in \eqref{KKT y2 C1} and \eqref{KKT y2 C2}.

 Analogously, recall from \eqref{eqn: dual function of DP2} that
	\begin{align}\label{proof of fact: perturb of y2 1 2}
	    \mathcal{Y}(\mub^\rpo, \zb^{r+1})  & = \arg \max_{\yb}~ \sum_{i=1}^N
\bigg(\phi_i(\yb_i,\zb_i^{r+1}) -
\yb_i^\top\frac{\qb}{N}\bigg)
-\la \mub^\rpo, \Ab \yb\ra  -\frac{\rho}{2}\|\Ab \yb \|^2 .
	\end{align}
For $\bar \yb^{r+1} \in \mathcal{Y}(\mub^{r+1}, \zb^{r+1})$, we have the following KKT condition
	\begin{align}
	 & \begin{bmatrix}
	       \Bb_1 \xb_1(\bar \yb_1^{r+1},\zb_1^{r+1}) - \frac{\qb}{N}\\
	       \vdots \\
	        \Bb_N \xb_N(\bar \yb_N^{r+1},\zb_N^{r+1}) -\frac{\qb}{N}
	  \end{bmatrix}
	  -\Ab^\top \mub^{r+1}
    - \rho \Lb^- \bar \yb^{r+1} =\zerob,
        \label{KKT y C1} \\
	  &  \nabla  F_i ( \xb_i(\bar \yb^{r+1},\zb_i^{r+1}) )
	   + \Bb_i^\top \bar \yb^{r+1} =\zerob,
    \forall i \in [N].
    \label{KKT y C2}
	\end{align}

 With fixed $\{\xb_i\big(\bar \yb_i^r,\zb_i^{r}\big)\}_{i=1}^N$,
\eqref{KKT y2 C1} and \eqref{KKT y2 C2} form a linear system of
	variable $\rho \bar \yb^r$ and $\mub^{r+1}$. Then, by applying the Hoffman bound in Lemma \ref{Hoffman}, the distance from $\bar \yb^{r+1}$  satisfying \eqref{KKT y C1}  and \eqref{KKT y C2} to the set $\Yc(\mub^{r+1},\zb^{r})$ can be bounded as
\begin{align} \label{eq: bound for y(mu r+1, zr) - y(mu r+1 , z r+1 ) square}
&\mathrm{dist}^2(\bar \yb^r, \mathcal{Y}(\mub^{r+1}, \zb^{r+1})) \notag \\
&\le \theta_2^2(\|\nabla_{\xb} F(\bar \yb^r, \zb^r))
-\nabla_{\xb} F(\xb(\bar \yb^{r+1}, \zb^{r+1}))\|^2
+\frac{1}{\rho^2}\|\tilde{\Bb}(\xb(\bar \yb^r-\xb(\bar \yb^{r+1}, \zb^{r+1}))\|^2)
\notag \\
&\le \theta_2^2((p+\gamma^+)^2\|\xb(\bar \yb^r, \zb^r)
-\xb(\bar \yb^{r+1}, \zb^{r+1})\|^2
+\frac{B_{\max}^2}{\rho^2}\|\xb(\bar \yb^r, \zb^r)-\xb(\bar \yb^{r+1}, \zb^{r+1})\|^2)
\notag\\
&\overset{\mathrm{(i)}}{\le} \underbrace{ \theta_2^2\bigg((p+\gamma^+)^2 \sigma_3 ^2+\frac{B_{\max}^2 \sigma_3 ^2}{\rho^2}\bigg)}_{\triangleq a_3}
\|\zb^r-\zb^{r+1}\|^2,
\end{align}
where $\tilde{\Bb} \triangleq \mathrm{diag} \{\Bb_1, \cdots,
\Bb_N\}$, and $\mathrm{(i)}$ is because of \eqref{sigma1}.
	\hfill $\blacksquare$

\section{Proof of lemma \ref{Fact y(z^r) - y(mu^{r+1, z^r})}} \label{appendix:
y(z^r) - y(mu^{r+1, z^r})}

   For $\bar \yb^r \in \mathcal{Y}(\mub^{r+1}, \zb^{r})$,
recall the KKT conditions in \eqref{KKT y2 C1} and \eqref{KKT y2 C2}.

For ${\yb}(\zb^{r}) \in \mathcal{Y}(\zb^r)$, recall from \eqref{eq: y(z)} and we
have
its KKT conditions as
\begin{align}
	 & \begin{bmatrix}
	       \Bb_1 \xb_1(\zb^{r}) - \qb/N \\
	       \vdots \\
	        \Bb_N \xb_N(\zb^{r}) - \qb/N
	  \end{bmatrix}
	  -\Ab^\top \mub(\zb^r) - \rho \Lb^- {\yb}(\zb^{r}) =\zerob,\label{KKT y(z^r) 1}
\\ \label{KKT y(z^r)}
	  & \nabla  f_i ( \xb_i(\zb^{r}) ) + p (\xb_i(\zb^{r}) - \zb_i^{r})
	+ \Bb_i^\top {\yb}_i(\zb^{r}) =\zerob, \forall i \in [N],\\
&      \Ab{\yb}(\zb^{r}) = \mathbf{0}. \label{KKT y(z^r) Ay=0}
	\end{align}
    With fixed $\{\xb_i\big(\bar \yb_i^r,\zb_i^{r}\big)\}_{i=1}^N$,
\eqref{KKT y(z^r) 1}-\eqref{KKT y(z^r) Ay=0} forms a linear system of
	variable $\rho \bar \yb^r$ and $\mub^{r+1}$. Then, by using the Hoffman bound in
Lemma \ref{Hoffman}, the distance from $\yb(\zb^r)$ and $\mub(\zb^r)$ satisfying
\eqref{KKT y(z^r) 1}-\eqref{KKT y(z^r) Ay=0} to the set $\Yc(\mub^{r+1},\zb^{r})$
can be bounded as
    \begin{align}\label{eq: y(z^r) - y(mu^{r+1, z^r} process 1}
    & \rho^2 \| \bar \yb^r - {\yb}(\zb^{r}) \|^2 + \|\mub^{r+1} -
    \mub(\zb^r)\|^2
    \notag \\
    & \leq
        {\theta_3^2 \sum_{i=1}^N \|\Bb_i \xb_i(\bar \yb_i^r, \zb_i^r) -
        \Bb_i\xb_i(\zb^r)\|^2 }
        \notag \\
        &~~~
        + \rho^2 \theta_3^2 \sum_{i=1}^N \| \nabla  f_i ( \xb_i(\bar
        \yb_i^r,\zb_i^{r})) -
        \nabla  f_i ( \xb_i(\zb^{r}))
         +p[\xb_i(\bar \yb_i^r,\zb_i^{r})-\xb_i(\zb^{r}) ]
        \|^2  +\theta_3^2 \|\Ab\bar \yb^r\|^2
    \notag\\
    & \overset{\mathrm{(i)}}{\leq} - \bigg(\frac{\theta_3^2}{p+\gamma^-}\bigg)
    \sum_{i=1}^N \la \Bb_i\xb_i(\bar \yb_i^r, \zb_i^r)- \Bb_i\xb_i(\zb^r)),   \bar
    \yb_i^r- {\yb}_i(\zb^r) \ra
        \notag \\
        &~~~
        + \rho^2 \theta_3^2 (p + \gamma^+ )\sum_{i=1}^N \la
        \nabla  f_i ( \xb_i(\bar \yb_i^r,\zb_i^{r}))-\nabla  f_i
        (\xb_i(\zb^{r}))
      + p[\xb_i(\bar \yb_i^r,\zb_i^{r}) - \xb_i(\zb^{r})],
\xb_i(\bar \yb_i^r,\zb_i^{r}) - \xb_i(\zb^{r})  \ra
       \notag\\
        &~~~+\theta_3^2 \|\Ab\bar \yb^r\|^2
    \notag\\
    & \overset{\mathrm{(ii)}}{\leq} - \bigg(\frac{\theta_3^2}{p+\gamma^-}\bigg)
    \sum_{i=1}^N \la \Bb_i\xb_i(\bar \yb_i^r, \zb_i^r)
- \Bb_i\xb_i(\zb^r)),  \bar \yb_i^r- {\yb}_i(\zb^r) \ra
        \notag \\
        &~~~
        +  \rho^2 \theta_3^2 (p + \gamma^+ )\sum_{i=1}^N \la -\Bb_i^\top\bar \yb_i^r
        + \Bb_i^\top{\yb}_i(\zb^r),
        \xb_i(\bar \yb_i^r,\zb_i^{r}) - \xb_i(\zb^{r})  \ra
        +\theta_3^2 \|\Ab\bar \yb^r\|^2 \notag \\
    & =
        - \bigg(\frac{\theta_3^2}{p+\gamma^-} + \rho^2 \theta_3^2 (p+\gamma^+)
        \bigg) \sum_{i=1}^N
       \la \Bb_i\xb_i(\bar \yb_i^r, \zb_i^r)- \Bb_i\xb_i(\zb^r)),   \bar \yb_i^r- {\yb}_i(\zb^r) \ra
        +\theta_3^2 \|\Ab\bar \yb^r\|^2
    \end{align}
    where $\mathrm{(i)}$ is obtained by \eqref{gradient dif bound of phi},
    $\mathrm{(ii)}$ is due to \eqref{KKT y2 C2} and \eqref{KKT y(z^r)}.
Further, note that we have
    \begin{align}\label{eq: y(z^r) - y(mu^{r+1, z^r} process 2}
    &-\sum_{i=1}^N\la \Bb_i\xb_i(\bar \yb_i^r, \zb_i^r) - \Bb_i\xb_i(\zb^r)),   \bar
    \yb_i^r- {\yb}_i(\zb^r) \ra \notag \\
    & \overset{\mathrm{(i)}}{=} -\la \Ab^\top( \mub^{r+1} - \mub(\zb^r)) + \rho
    \Lb^- (\bar \yb^r- {\yb}(\zb^r) ),
      \bar \yb^r- {\yb}(\zb^r) \ra
    \notag \\
    & \overset{\mathrm{(ii)}}{=} -\la  \mub^{r+1} - \mub(\zb^r) ,\Ab\bar \yb^r\ra
     - \| \bar \yb^r- {\yb}(\zb^r)  \|_{\rho \Lb^-}^2
    \notag \\
    & \leq  -\la  \mub^{r+1} - \mub(\zb^r) ,\Ab\bar \yb^r \ra,
   \end{align}
   where $\mathrm{(i)}$ can be deduced by taking the difference between \eqref{KKT
   y2 C1} and \eqref{KKT y(z^r) 1}, and $\mathrm{(ii)}$ is due to \eqref{KKT
   y(z^r) Ay=0}.
 By inserting \eqref{eq: y(z^r) - y(mu^{r+1, z^r} process 2} into \eqref{eq: y(z^r)
 - y(mu^{r+1, z^r} process 1}, we have
   \begin{align}
   &\rho^2 \|\bar \yb^r -{\yb}(\zb^{r}) \|^2 + \|\mub^{r+1} -
   \mub(\zb^r)\|^2
    \notag \\
   &
   \leq - \bigg(\frac{\theta_3^2}{p+\gamma^-} + \rho^2 \theta_3^2 (p+\gamma^+)
   \bigg)
\la \Ab \bar \yb^r,
     \mub^{r+1} - \mub(\zb^r)  \ra + \theta_3^2 \|\Ab \bar \yb^r\|^2.
   \end{align}
   Rearranging the above inequality, we have
   \begin{align}
   & \|\bar \yb^r -{\yb}(\zb^{r}) \|^2
    \notag \\
   &
   \leq \frac{\theta_3^2}{\rho^2} \|\Ab \bar \yb^r\|^2- \frac{1}{\rho^2}\|\mub^{r+1}
   -
   \mub(\zb^r)\|^2
    - \bigg(\frac{\theta_3^2}{\rho^2(p+\gamma^-)} +  \theta_3^2 (p+\gamma^+)
   \bigg)
\la \Ab \bar \yb^r,
     \mub^{r+1} - \mub(\zb^r)  \ra,
     \notag
     \\
     & \overset{\mathrm{(i)}}{\leq} \frac{\theta_3^2}{\rho^2} \|\Ab \bar \yb^r\|^2-
     \frac{1}{\rho^2}\|\mub^{r+1} - \mub(\zb^r)\|^2
    + \bigg(\frac{\theta_3^2}{\rho^2(p+\gamma^-)} +  \theta_3^2 (p+\gamma^+)
   \bigg)
\| \Ab \bar \yb^r\|  \|\mub^{r+1} - \mub(\zb^r)\|, \notag  \\
&\overset{\mathrm{(ii)}}{\leq} \frac{\theta_3^2}{\rho^2} \|\Ab \bar \yb^r\|^2
        +\frac{\tau}{2} \bigg(\frac{\theta_3^2}{\rho^2(p+\gamma^-)} +  \theta_3^2
        (p+\gamma^+) \bigg) \| \Ab \bar \yb^r\|^2
    \notag \\
   & ~~~
    +\bigg[\frac{\bigg(\frac{\theta_3^2}{\rho^2(p+\gamma^-)} +  \theta_3^2
    (p+\gamma^+) \bigg)}{2\tau }- \frac{1}{\rho^2}\bigg] \|\mub^{r+1} -
    \mub(\zb^r)\|^2,\label{eqn: last}
   \end{align}
   where $\mathrm{(i)}$ is obtained by Cauchy-Schwartz inequality, and
   $\mathrm{(ii)}$ is due to Young's inequality with $\tau>0$.

   If we choose $\textstyle \tau = \frac{\theta_3^2}{p+\gamma^-} + \rho^2 \theta_3^2
   (p+\gamma^+)  $, then $\textstyle \frac{1}{2\tau
   }\bigg(\frac{\theta_3^2}{\rho^2(p+\gamma^-)} +  \theta_3^2 (p+\gamma^+) \bigg)-
   \frac{1}{\rho^2}  <0$. By substituting this choice of $\tau$ into \eqref{eqn:
   last},  we obtain the desired \eqref{y(z^r) - y(mu^{r+1, z^r})}.
        \hfill $\blacksquare$


\section{Proof of Lemma \ref{Fact descent potential}(b)}\label{appendix: lemma 8b}

    According to Lemma \ref{Fact descent potential}(a), the upper bound for $\Phi^{r+1}
    - \Phi^r$ is non-positive. Therefore we can have the conclusion that the
    potential function $\Phi$ is non-increasing with $r$.

Recall from \eqref{eq: potential function G} that:
    \begin{align}
      & G(\xb^{r+1}, \yb^{r+1}, \mub^{r+1}, \zb^{r+1})\notag \\
      &=\sum_{i=1}^N \bigg( f_i(\xb_i^{r+1})  +
      \frac{p}{2}\|\xb_i^{r+1}-\zb_i^{r+1}\|^2 + (\yb_i^{r+1})^\top\Bb_i\xb_i^{r+1}  - \frac{1}{N}(\yb_i^{r+1})^\top\bb\bigg)
-(\mub^{r+1})^\top\Ab \yb^{r+1} -\frac{\rho}{2}\|\Ab \yb^{r+1} \|^2.
    \end{align}
    Note that,
    \begin{align} \label{G part bound}
        &\sum_{i=1}^N \bigg(  (\yb_i^{r+1})^\top\Bb_i\xb_i^{r+1} -
        \frac{1}{N}(\yb_i^{r+1})^\top\bb\bigg) -(\mub^{r+1})^\top\Ab \yb^{r+1}
        \notag \\
        &=\sum_{i=1}^N (\yb_i^{r+1})^\top (  \Bb_i\xb_i^{r+1} - \bb/N - \Ab_i^\top
        \mub^{r+1} ) \notag \\
        &\overset{\mathrm{(i)}}{=} \sum_{i=1}^N (\yb_i^{r+1})^\top ( 2\rho |\Nc_i| \yb_i^{r+1} - \rho  \Lb^+_i
        \yb^r ) \notag \\
        & = \rho (\yb^{r+1})^\top ( 2 \Db \yb^{r+1} -  \Lb^+ \yb^r )
        \notag \\
        & \overset{\mathrm{(ii)}}{=} \rho (\yb^{r+1})^\top (  \Lb^+ (\yb^{r+1} -  \yb^r)  +  \Lb^- \yb^{r+1}).
        \notag\\
        & = \rho \|\yb^{r+1}\|^2_{\Lb^-} + \frac{\rho}{2}
        \|\yb^{r+1}-\yb^r\|^2_{\Lb^+} +  \frac{\rho}{2}\|\yb^{r+1}\|^2_{\Lb^+} -
        \frac{\rho}{2}\|\yb^r\|^2_{\Lb^+},
    \end{align}
     where $\mathrm{(i)}$ is obtained by \eqref{y: update} and $\mathrm{(ii)}$ is due to $\Lb^+ = 2\Db - \Lb^-$.

    Thus, recall $\tilde G$, then combine the results in \eqref{eq: potential function G} and \eqref{G part bound}. We have
    \begin{align}
            \tilde{G}^{r+1}
          &\triangleq G(\xb^{r+1}, \yb^{r+1}, \mub^{r+1}, \zb^{r+1}) +
          \frac{\alpha}{2}\|\yb^{r+1}\|^2_{\Lb^-}+  \frac{\rho}{2}\|\yb^{r+1} - \yb^r\|^2_{\Lb^+ } \notag \\
          & =\sum_{i=1}^N  f_i(\xb_i^{r+1})  +  \frac{p}{2}\|\xb^{r+1}-\zb^{r+1}\|^2
          + \bigg( \frac{\alpha}{2} + \frac{\rho}{2} \bigg) \|\yb^{r+1}\|^2_{\Lb^-}
          \notag \\
          &~~~+ {\rho}\|\yb^{r+1}-\yb^r\|^2_{\Lb^+}
          + \frac{\rho}{2}\|\yb^{r+1}\|^2_{\Lb^+} -
          \frac{\rho}{2}\|\yb^r\|^2_{\Lb^+},
    \end{align}
    which implies that
    \begin{align}
        \sum_{r=0}^T (\tilde{G}^{r+1}-\underline{f}) > -\frac{\rho}{2}\|\yb^0\|_{\Lb^+}^2>-\infty,
    \end{align} for any $T$.
	Due to the weak duality, $d(\mub^r;\zb^r) \geq \Lc_\rho(\yb^r,\mub^r;\zb^r)$.
We further obtain
$\sum_{r=0}^T (\Phi^{r+1}-\underline{f}) > -\frac{\rho}{2}\|\yb^0\|^2_{\Lb^+}> -\infty.$
    By this fact and the conclusion of (a) that $\Phi^{r+1} $ is non-increasing,
    $\Phi^{r+1}$
    must be lower bounded, which is denoted as $\underline{\Phi}$.
\hfill $\blacksquare$


\section{Proof of lemma \ref{lemma9}} \label{appendix: proof of lemma 9}
Let $(\xb(\zb), \yb_0(\zb))$ be the primal-dual solution of problem  \eqref{consensus problem 2 PX}.  They satisfy the following KKT conditions
        \begin{subequations}
        \begin{align}
            &   \nabla_{\xb} f_i(\xb_i(\zb_i)) +p( \xb_i(\zb_i)-  \zb_i)
             +\Bb_i^\top  \yb_0(  \zb ) = 0, \forall i \in [N],\label{eq: constrained strongly convex
         KKT 1}\\
            &\sum_{i=1}^N \Bb_i  \xb_i(\zb_i)=    \qb.\label{eq: constrained strongly convex KKT 2}
        \end{align}
        \end{subequations}
 Suppose that $\|\zb-  \xb(  \zb)\|^2\leq \epsilon$ for some $\zb$.  We analyze the gap between $\zb$ and a KKT solution of problem {\sf (P)} (see \eqref{eqn: KKT of P}).  Firstly,  consider the following chain of bounds

   \vspace{-0.3cm}
        {\small \begin{align}\label{eq: gradient res upper bound}
        & \sum_{i=1}^N \left\| \nabla_{\xb}f_i(\zb_i)+ \Bb_i^\top  \yb_0(\zb) \right\|^2 \notag \\
        & \leq \sum_{i=1}^N \left\|
                \nabla_{\xb}f_i(\zb_i)
                \! -\!\nabla_{\xb}f_i(\xb_i(\zb_i))\!+\!\nabla_{\xb}f_i(\xb_i(\zb_i))\!+\!  \Bb_i^\top  \yb_0(\zb ) \right\|^2 \notag \\
        & \overset{{\mathrm{(i)}}}{=} \sum_{i=1}^N \left\| \nabla_{\xb}f_i(\zb_i)-\nabla_{\xb}f_i(\xb_i(\zb_i)) + p(\zb_i - \xb_i(\zb_i))\right\|^2\notag \\
        & {\leq} \sum_{i=1}^N 2\left\| \nabla_{\xb}f_i(\zb_i)-\nabla_{\xb}f_i(\xb_i(\zb_i))\right\|^2
               \notag\\
               &~~~~ +\sum_{i=1}^N 2\left\| p(\zb_i - \xb_i(\zb_i))\right\|^2\notag \\
        &\overset{{\mathrm{(ii)}}}{\leq}  2(p^2+(\gamma^+)^2)  \sum_{i=1}^N \left\|  \zb_i - \xb_i(\zb_i) \right\|^2 \notag \\
        & \leq  2(p^2+(\gamma^+)^2) \epsilon,
        \end{align}}

        \vspace{-0.3cm}
\noindent        where $\mathrm{(i)}$ is by \eqref{eq: constrained strongly convex
         KKT 1},
         and $\mathrm{(ii)}$ is by Assumption \ref{assumption smooth}.

        Secondly, we have

   \vspace{-0.3cm}
        {\small         \begin{align}\label{eq: infeasibility upper bound}
         \left\|\sum_{i=1}^N \Bb_i  \zb_i -  \qb\right\|^2
         & =  \left \|\sum_{i=1}^N \Bb_i  \zb_i - \sum_{i=1}^N \Bb_i \xb_i(\zb_i)
        +\sum_{i=1}^N \Bb_i \xb_i(\zb_i)-  \qb\right\|^2 \notag\\
         & \overset{{\mathrm{(i)}}}{=}   \left\|\sum_{i=1}^N \Bb_i  (\zb_i-  \xb_i(\zb_i))\right\|^2 \notag\\
          & \overset{{\mathrm{(ii)}}}{\leq}   NB_{\max}^2  \sum_{i=1}^N \|\zb_i-  \xb_i(\zb_i)\|^2 \notag\\
        & \leq   N B_{\max}^2 \epsilon,
        \end{align}}

           \vspace{-0.3cm}
 \noindent       where ${\mathrm{(i)}}$ is by\eqref{eq: constrained strongly convex KKT 2} and ${\mathrm{(ii)}}$ is by $B_{\max} \triangleq \max_{i\in[N]}\|\Bb_i\|$.
%
        Thus, we can conclude that $\zb$ is a $\kappa \epsilon$-KKT solution of problem {\sf (P)} with
        \begin{align}
            \kappa \triangleq  \max\{ 2(p^2+(\gamma^+)^2, NB_{\max}^2\}.\notag
        \end{align}
    \hfill $\blacksquare$

\section{Proof of Corollary 1}\label{app: corollary 1}
   To determine the dependence of the graph topology in the convergence rate, we make the analysis by 3 steps, including  (1) simplifying the upper bound of the term $\|\zb^t - \xb(\zb^t)\|$, (2) deriving the bound for $\theta_1,\theta_2, \theta_3$ which appear in Lemmas \ref{fact: perturb of y 2}, \ref{fact: perturb of y} and \ref{Fact y(z^r) - y(mu^{r+1, z^r})}, and  (3) determining the dependence of the bound on the topology and data matrix.

    \begin{itemize}
        \item [{(1)}] {\bf{Simplify the upper bound of $\|\zb^t - \xb(\zb^t)\|$:}}
        Recall in \eqref{eq: bound for z x(z) 1}  that
        \begin{align}
    	    \|\zb^t-\xb(\zb^t)\|& \leq \|\zb^t-\xb^{t+1}\| +\| \xb^{t+1} - \xb(\zb^t)\|
            \notag  \\
            	&\leq {\textstyle \sqrt{\frac{\Phi^0-\underline{\Phi}}{r}}} \bigg(\frac{1}{\beta
            \sqrt{C_3}} + \frac{\sigma_2\sqrt{a_2}\lambda_{\max}}{\sqrt{C_2}}  +
            \frac{\sigma_2
            \sqrt{a_4}}{\sqrt{\alpha}}\bigg),
        \end{align}
        where $\underline{\Phi}$ is the lower bound of $\sum_{r=0}^T \Phi^{r+1}$ for any iteration $T$, $\lambda_{\max}$ is the largest eigenvalue of $\Lb^+$, $\sigma_2$ and $a_2$ are defined in Lemma \ref{fact: pertubration of x_i} and Lemma \ref{fact: perturb of y 2} respectively,  and $C_2$ and $C_3$ are in Lemma \ref{Fact descent potential} of the revised manuscript.

        Since $C_2 \geq \frac{\rho}{2}$ shown in the proof of Lemma \ref{Fact descent potential}(a),  the second term of the above upper bound can be simplified as $ \frac{\sigma_2\sqrt{a_2}\lambda_{\max}}{\sqrt{C_2}} \leq \frac{\sqrt{2\rho a_2} \sigma_2\lambda_{\max} }{\rho}$.
        Further inserting values of $C_3$, $\sigma_2$, $a_2,a_4$ defined in Lemmas \ref{Fact descent potential}, \ref{fact: pertubration of x_i},  \ref{fact: perturb of y 2},  and \ref{Fact y(z^r) - y(mu^{r+1, z^r})},  respectively,  we obtain
        \begin{align}\label{eq: z-x(z) bound 1}
            &\|\zb^t-\xb(\zb^t)\|\notag \\
            & \leq
             {\textstyle \sqrt{\frac{\Phi^0-\underline{\Phi}}{r}}} \bigg\{
             \frac{1}{\beta \sqrt{p \left(-\frac{1}{2} + \frac{1}{\beta} - \frac{20p\sigma_2^2 a_2 \lambda_{\max}^2 }{\rho}
              - \frac{\rho (\sigma_1^2 + 2\sigma_2^2 a_3)}{10p \sigma_2^2 a_2 \lambda_{\max}^2} \right)}} \bigg.\notag\\
              &~~~ ~~~ ~~~ \bigg.
             + \frac{\sqrt{2\rho} }{\rho} \cdot \frac{B_{\max}(p+\gamma^+ + 1)}{p+\gamma^-} \cdot \theta_1 \sqrt {\frac{\theta_1^2}{\rho^2} \left(\frac{1}{p+\gamma^-} + \rho^2 (p+\gamma^+) \right)^2 +2}\cdot \lambda_{\max}\bigg.\notag\\
              &~~~ ~~~ ~~~ \bigg.
             + \frac{B_{\max}(p+\gamma^++1)}{p+\gamma^-} \cdot \theta_3 \sqrt{\frac{1}{2} \left[\frac{1}{p+\gamma^-} + \rho^2 (p+\gamma^+)\right]^2 \theta_3^2 + \frac{1}{\rho^2}}  \cdot
             \frac{1}{\sqrt{\alpha}} \bigg\} .
             \end{align}
             \begin{itemize}
                 \item [a)] \underline{\emph{{\textbf{Choose parameter $\beta, p, \rho, \alpha$.}}}}

                Let
                \begin{align}
                 & \beta = \frac{1}{ 2\left(\frac{1}{2}   + \frac{20p\sigma_2^2 a_2 \lambda_{\max}^2 }{\rho}
               +\frac{\rho (\sigma_1^2 + 2\sigma_2^2 a_3)}{10p \sigma_2^2 a_2 \lambda_{\max}^2} \right)},\\
               & p = \rho = -2\gamma^-,\\
               & \alpha = \min\left\{\frac{\rho}{5}, \frac{\rho}{8a_1\lambda_{\max}^2}\right\}.
               \end{align}

             \item [b)]  \underline{\emph{{\textbf{Inserting $p, \rho, \alpha$ into \eqref{eq: z-x(z) bound 1}}}}} followed by making square,
                    we have
                  \small{
                \begin{align}\label{eq: z-x(z) bound 2}
                    &\|\zb^t-\xb(\zb^t)\|^2 \notag \\
                    &  \leq   \frac{\Phi^0-\underline{\Phi}}{r}  \left\{
                      D_1 B_{\max}  \lambda_{\max} \sqrt{D_2  \theta_1^4 +D_{3} \theta_1^2  }\right.\notag \\
                    &~~~\left.
                    +\sqrt{D_4 +  B_{\max}^2 \lambda_{\max}^2 (D_5 \theta_1^4 + D_6\theta_1^2)
                    + \frac{D_7 + D_8 B_{\max}^2\theta_2^2}{B_{\max}^2\lambda_{\max}^2(D_9 \theta_1^4 + D_{10}\theta_1^2) }} \right.\notag \\
                    &~~~\left.
                        +D_{11} B_{\max}\theta_3 \sqrt{ D_{12} \theta_3^2 + D_{13}  } \cdot \max\{D_{14}, \lambda_{\max} \sqrt{D_{15}   \theta_1^2 + D_{16}}\}\right\}^2,
                \end{align}}
                where $D_1$-$D_{16}$ are positive constants related to $\gamma^+, \gamma^-$.

        \end{itemize}
    \item      [{(2)}] {\bf{Determine $\theta_1, \theta_2, \theta_3$ by Hoffman bound:}}

\vspace{0.3cm}
    {\bf{Notation:}}
    For a matrix $\bar{\Mb}$, define $\sigma_{\max}(\bar{\Mb})$ to be the largest singular value of $\bar{\Mb}$ and define $\sigma_{\min}(\bar{\Mb})$ to be the smallest nonzero singular value of $\bar{\Mb}$.

\vspace{0.3cm}
    We can give an upper bound of the Hoffman constant $\theta_1, \theta_2, \theta_3$. Specializing Lemma 3.2.3 in \cite{facchinei2003finite} to the case where the linear system is defined only by linear equality, we have the following Lemma:

      \begin{Lemma} \label{lem: hoffman bound constant def}
        Let $S=\{\xb\in \mathbb{R}^n\mid \Mb \xb=\vb\}$ be a linear system. Then for any $\bar{\xb}\in \mathbb{R}^n$, the distance from $\bar{\xb}$ to the set $S$ is bounded as:
        $$\mathrm{dist}(\bar{\xb}, S)\le \frac{\sigma_{\max}(\bar{\Mb})}{\sigma^2_{\min}(\bar{\Mb})}\|\Mb\bar{\xb}-\vb\|,$$
        where $\bar{\Mb}$ is any submatrix whose rows span the row space of $\Mb$.\hfill$ \blacksquare$
        \end{Lemma}

        Recall in the proof of Lemma \ref{fact: perturb of y 2} that KKT conditions \eqref{KKT y2 C1}-\eqref{KKT y2 C2} form a linear system of $\rho \bar{\yb}^r$, with fixed $\{\xb_i(\bar{\yb}_i^r, \zb_i^r)\}_{i=1}^N$ as follows
        \begin{align}
            \underbrace{\begin{bmatrix}\Lb^-\\ \Bb^T_{\mathrm{diag}} \end{bmatrix}}_{\triangleq \Mb_1} \rho \bar{\yb}^r = \ub_1 ,
        \end{align}
        where $\ub_1$ is the corresponding vector defined by \eqref{KKT y2 C1}-\eqref{KKT y2 C2} and $\Bb_{\mathrm{diag}}\triangleq \mathrm{diag}(\Bb_1,\cdots, \Bb_N)$. Then, by Lemma \ref{lem: hoffman bound constant def}, we have $\theta_1$ in Lemma \ref{fact: perturb of y 2} is $\theta_1 = \theta(\Mb_1) = \sigma_{\max}(\Mb_1)/\sigma_{\min}^2(\Mb_1)$. Moreover, since the proof of Lemma \ref{fact: perturb of y} is based on the same linear system, we obtain $\theta_1 = \theta_2$.

        Similarly, to prove Lemma \ref{Fact y(z^r) - y(mu^{r+1, z^r})}, KKT conditions \eqref{KKT y(z^r) 1}-\eqref{KKT y(z^r) Ay=0} can be viewed as a linear system of variables $\rho\bar{\yb}^r$ and $\mub^{r+1}$
        \begin{align}
            \underbrace{\begin{bmatrix}
            \Ab^\top & \Lb^-\\
            \mathbf{0} & \Ab\\
            \mathbf{0} &\Bb^T_{\mathrm{diag}}
            \end{bmatrix}}_{\triangleq  \Mb_2}
            \begin{bmatrix}
            \mub^{r+1}\\
            \rho \bar{\yb}^r
            \end{bmatrix}  = \ub_2 ,
        \end{align}
        where $\ub_2$ is the corresponding vector defined by \eqref{KKT y(z^r) 1}-\eqref{KKT y(z^r) Ay=0}. Therefore, we obtain $\theta_3 = \theta(\Mb_2) = \sigma_{\max}(\Mb_2)/\sigma_{\min}^2(\Mb_2) $.

        \vspace{0.2cm}
        Notice that $\Mb_1$ and $\Mb_2$ are related to the network topology and the data matrix $\Bb$. In general, we cannot write $\theta_1, \theta_2, \theta_3$ in terms of the eigenvalues of $\Lb^-$. However, if $\Bb=[\Bb_1,\cdots,\Bb_N]$ is of full row rank, we can write $\theta_1, \theta_2,\theta_3$ concerning the singular values of $\Lb^-$.

        \vspace{0.5cm}
        Next, for $\Mb_1$ and $\Mb_2$, we obtain their lower bound of $\sigma_{\min}$ and the upper bound of $\sigma_{\max}$,  respectively by the following two Lemmas.

       \begin{Lemma}\label{lem: sigma min M}
            If $\Bb$ is of full row rank, we have
            \begin{align}
            & \sigma_{\min}(\Mb_1)\ge \max\left\{
            \frac{\sigma_{\min}(\Lb^-) }{ 1+ \zeta_B },
            \frac{ \sigma_{\max}(\Bb_{\mathrm{diag}})}{2(1+\zeta_B)}
        \right\}.\\
            & \sigma_{\min}(\Mb_2)
                  \geq \frac{\max\left\{
                    \frac{\sqrt{\sigma_{\min}(\Lb^-)} }{ 1+ \zeta_B },
            \frac{ \sigma_{\max}(\Bb_{\mathrm{diag}})}{2(1+\zeta_B)}\right\}}
                    {3+\sqrt{\sigma_{\max}(\Lb^-)} },
            \end{align}
            where $\zeta_B \triangleq \frac{ 2\sqrt{N}\cdot \sigma_{\max}(\Bb_{\mathrm{diag}})}{\sigma_{\min}(\Bb)} $.
        \end{Lemma}
\vspace{0.2cm}
 {\bf{Proof:}} See Appendix \ref{app: proof of sigma min M}.\hfill$ \blacksquare$
\vspace{0.3cm}

        \begin{Lemma}\label{lem: sigma max M}
            We have
            \begin{align}
                & \sigma_{\max}(\Mb_1)\leq  \sigma_{\max}(\Lb^-)+\sigma_{\max}(\Bb_{\mathrm{diag}}) ,\\
            &\sigma_{\max}(\Mb_2)\leq  2\sqrt{\sigma_{\max}(\Lb^-)}+\sigma_{\max}(\Lb^-)+\sigma_{\max}(\Bb_{\mathrm{diag}}).
            \end{align}
        \end{Lemma}
{\bf{Proof:}} See Appendix \ref{app: proof of sigma max M}.\hfill$ \blacksquare$

\vspace{0.3cm}
        Thus, combining these bounds and Lemma \ref{lem: hoffman bound constant def} yields the desired bound for $\theta_1, \theta_2$ and $\theta_3$ as shown below
        \begin{align}
         \label{eq: upper bound of theta 12}
            & \theta_1 = \theta_2  \leq \frac{ \sigma_{\max}(\Lb^-)+\sigma_{\max}(\Bb_{\mathrm{diag}}) }
            {\max\left\{
                    \left(\frac{\sigma_{\min}(\Lb^-) }{ 1+ \zeta_B }\right)^2,
                    \left(\frac{ \sigma_{\max}(\Bb_{\mathrm{diag}})}{2(1+\zeta_B)}\right)^2
                \right\}}\notag\\
                &~~~ ~~~ ~~~
                 \overset{\mathrm{(i)}}{\leq}
                     (1+ \zeta_B)^2 \frac{\sigma_{\max}(\Lb^-)}{\sigma_{\min}^2(\Lb^-) }
                     + \frac{4(1+ \zeta_B)^2}{ \sigma_{\max}(\Bb_{\mathrm{diag}})}  ,\\
             \label{eq: upper bound of theta 3}
            &\theta _3 \leq \frac{ 2\sqrt{\sigma_{\max}(\Lb^-)}+\sigma_{\max}(\Lb^-)+\sigma_{\max}(\Bb_{\mathrm{diag}})}
                {\frac{\max\left\{
                            \frac{ \sigma_{\min}(\Lb^-)  }{ (1+ \zeta_B)^2 },
                    \frac{ \sigma_{\max}^2(\Bb_{\mathrm{diag}})}{4(1+\zeta_B)^2}\right\}}
                            {(3+\sqrt{\sigma_{\max}(\Lb^-)} )^2}} \notag\\
                &~~~ \overset{\mathrm{(ii)}}{\leq}
                    \frac{\left(2\sqrt{\sigma_{\max}(\Lb^-)}+\sigma_{\max}(\Lb^-)\right)(3+\sqrt{\sigma_{\max}(\Lb^-)} )^2 (1+ \zeta_B)^2 }
                    {\sigma_{\min}(\Lb^-) }\notag\\
                    &~~~ ~~~ ~~~
                    +
                        \frac{4(1+\zeta_B)^2 }{\sigma_{\max} (\Bb_{\diag})}\left(3+\sqrt{\sigma_{\max}(\Lb^-)} \right)^2
        \end{align}
        where $\zeta_B \triangleq \frac{ 2\sqrt{N}\cdot \sigma_{\max}(\Bb_{\diag})}{\sigma_{\min}(\Bb)} $, $\mathrm{(i)}$ and $\mathrm{(ii)}$ are due to the fact that $\frac{ a+b }{\max\{c, d\}} \leq \frac{a}{c} + \frac{b}{d}$.

            \item [{(3)}] {\bf{Determine the dependence on topology of convergence rate:}}

                Combining \eqref{eq: z-x(z) bound 2} , \eqref{eq: upper bound of theta 12}, \eqref{eq: upper bound of theta 3}, we obtain the result in Corollary 1.
                \hfill$ \blacksquare$
            \end{itemize}

\section{Proof of Lemma \ref{lem: sigma min M}}\label{app: proof of sigma min M}

    \begin{itemize}
        \item [a)] \underline{\emph{{\textbf{Solve the lower bound for $\sigma_{\min}(\Mb_1)$}}}}

    Let $\vb \in \mathbb{R}^{NM}$ be a vector. We write $\vb =\vb_1+\vb_2$, where $\vb_1$ is in the null space of $\Lb^-$ and $\vb_2$ is orthogonal to $\vb_1$. Then, $\Lb^-\vb_1 = \mathbf{0}$.

    Recall $\zeta_B \triangleq \frac{ 2\sqrt{N}\cdot \sigma_{\max}(\Bb_{\mathrm{diag}})}{\sigma_{\min}(\Bb)} $. Then, we can divide into it two cases to analyze the lower bound of $\sigma_{\min}(\Mb_1)$.

\vspace{0.2cm}
    \textbf{Case 1.} If $\|\vb_1\|\le  \zeta_B\|\vb_2\|$, we have
    \begin{align}
        \| \Mb_1\vb\|
        &\ge \|\Lb^-\vb\| \overset{\mathrm{(i)}}{=}  \|\Lb^-\vb_2\| \overset{\mathrm{(ii)}}{\ge}  \sigma_{\min}(\Lb^-)\|\vb_2\| \overset{\mathrm{(iii)}}{\ge} \frac{\sigma_{\min}(\Lb^-)\|\vb\|}{ 1+ \zeta_B },
    \end{align}
    where $\mathrm{(i)}$ is because of $\Lb^-\vb_1 = \mathbf{0}$,  $\mathrm{(ii)}$  is due to the definition of  the smallest singular value $\sigma_{\min}$, and  $\mathrm{(iii)}$  is obtained by triangle inequality and the assumption $\|\vb_1\|\le \zeta_B \|\vb_2\|$, i.e., $\|\vb\| = \|\vb_1 + \vb_2\| \leq \|\vb_1\| + \|\vb_2\| \leq (1+ \zeta_B)\|\vb_2\|$.

\vspace{0.2cm}
    \textbf{Case 2.} Suppose $\|\vb_1\|\ge \zeta_B\|\vb_2\|$.
    Notice that $\vb_1=[\yb^\top,\yb^\top,\cdots,\yb^\top]^\top\in \mathbb{R}^{NM}$ for some $\yb\in \mathbb{R}^M$.
    We have
    \begin{align}
    \|\Mb_1 \vb\|
    &\ge \|\Bb^\top_{\diag}\vb\|\notag \\
    & = \|\Bb^\top_{\diag}\vb_1 + \Bb^\top_{\diag}\vb_2\|\notag \\
    & \overset{\mathrm{(i)}}{\ge} \|\Bb^\top_{\diag}\vb_1 \|- \| \Bb^\top_{\diag}\vb_2\|\notag \\
    & \geq  \|\Bb^\top \yb \|- \sigma_{\max}(\Bb_{\diag})\|  \vb_2\|\notag \\
    & \overset{\mathrm{(ii)}}{\ge} \sigma_{\min} (\Bb) \|\yb\| -  \sigma_{\max}(\Bb_{\diag}) \cdot \frac{1}{\zeta_B} \|\vb_1\| \notag\\
    & \overset{\mathrm{(iii)}}{=} \sigma_{\min} (\Bb) \frac{1}{\sqrt{N}}\|\vb_1\| - \sigma_{\min} (\Bb) \frac{1}{2\sqrt{N}}\|\vb_1\|\notag\\
    &= \sigma_{\min} (\Bb) \frac{1}{\sqrt{N}}\|\vb_1\| \notag\\
    &\overset{\mathrm{(iv)}}{\ge}  \frac{ \sigma_{\max}(\Bb_{\diag})}{2\zeta_B} \|\vb_1\|,\notag\\
    &\overset{\mathrm{(v)}}{\ge}  \frac{ \sigma_{\max}(\Bb_{\diag})}{2(1+\zeta_B)} \|\vb\|,
    \end{align}
where $\mathrm{(i)}$ is by triangular inequality, $\mathrm{(ii)}$ is because of assumption $\|\vb_1\|\ge \zeta_B\|\vb_2\|$, $\mathrm{(iii)}$ and $\mathrm{(iv)}$ are both due to the definition of $\zeta_B$, and $\mathrm{(v)}$ is obtained by triangle inequality and assumption $\|\vb_1\|\ge \zeta_B\|\vb_2\|$, i.e., $\|\vb\| = \|\vb_1 + \vb_2\| \leq \|\vb_1\| + \|\vb_2\| \leq (1+ 1/\zeta_B)\|\vb_1\|$.

    Thus, we have
    \begin{align}
        \sigma_{\min}(\Mb_1)\ge \max\left\{
            \frac{\sigma_{\min}(\Lb^-) }{ 1+ \zeta_B },
            \frac{ \sigma_{\max}(\Bb_{\diag})}{2(1+\zeta_B)}
        \right\}.
    \end{align}
        \item [b)] \underline{\emph{{\textbf{Solve the lower bound for $\sigma_{\min}(\Mb_2)$}}}}

            Before solving $\sigma_{\min}(\Mb_2)$, we give the following lemmas.

            \begin{Lemma}\label{lem: block matrix singular value}
            For block matrix $\hat{\Cb} = [\Cb_1^\top, \Cb_2^\top]^\top$ and $\tilde{\Cb} = [\Cb_3, \Cb_4]$, we have the following conclusions.
                        \begin{align}
                            & \sigma_{\max} (\hat{\Cb})
                            \leq \sigma_{\max}(\Cb_1) + \sigma_{\max}(\Cb_2) . \label{eq: block col singular val}  \\
                            & \sigma_{\max} (\tilde{\Cb})
                            \leq  \sigma_{\max}(\Cb_3) + \sigma_{\max}(\Cb_4) .\label{eq: block row singular val}
                        \end{align}
            \end{Lemma}
            {\bf{Proof:}} See Appendix \ref{app: proof of block matrix singular value}.\hfill$ \blacksquare$

            \vspace{0.5cm}
            Let $\Pb=\begin{bmatrix}
            \mathbf{I}&-\Ab^\top &\mathbf{0}\\
            \mathbf{0}&\mathbf{I}&\mathbf{0}\\
            \mathbf{0}&\mathbf{0}&\mathbf{I}
            \end{bmatrix}$.
            Then $\Pb\Mb_2=\begin{bmatrix}
            \Ab^\top& \mathbf{0} \\
            \mathbf{0}& \Ab\\
              \mathbf{0}&\Bb_{\diag}^\top
            \end{bmatrix}$.
            Since $\Pb$ is nonsingular, the null space of $\Mb_2$ is equal to that of $\Pb\Mb_2$.
            Therefore, any $\vb$ orthogonal to the null space of $\Mb_2$ is also orthogonal to the null space of $\Pb\Mb_2$.

            Pick $\vb$ satisfying $\|\Mb_2\vb\| = \sigma_{\min}(\Mb_2) \|\vb\|$. We then have
            \begin{align}\label{eq: PM2}
            \sigma_{\min}(\Pb\Mb_2)\|\vb\|
             \le  \|\Pb\Mb_2\vb\| \le \sigma_{\max}(\Pb)\|\Mb_2\vb\|   = \sigma_{\max}(\Pb) \sigma_{\min}(\Mb_2) \|\vb\| .
            \end{align}
            Therefore,
            \begin{align}\label{eq: M_2 min singular val 1}
                \sigma_{\min}(\Mb_2) \geq \frac{\sigma_{\min}(\Pb\Mb_2) }{\sigma_{\max}(\Pb)}.
            \end{align}

            By Lemma \ref{lem: block matrix singular value}, we have $\sigma_{\max}(\Pb)\leq 3+\sqrt{\sigma_{\max}(\Lb^-)} $.

            Moreover, by viewing $\Pb\Mb_2$ as a block diagonal matrix, it is obvious that $\sigma_{\min}(\Pb\Mb_2)$ $= \min\{\sigma_{\min}(\Ab),$ $\sigma_{\min}(\Mb_3)\}$, where $\Mb_3\triangleq [\Ab^\top, \Bb_{\diag}]^\top$.
            Since $\Lb^- = \Ab^\top \Ab$, then we can easily prove that $\text{Null}(\Lb^-) = \text{Null}(\Ab ) $. Then, the lower bound of $\sigma_{\min}(\Pb\Mb_2)$ can be solved by following the similar method to solve that of $\sigma_{\min}( \Mb_1)$. Thus, we have
            \begin{align}\label{eq: M3}
                \sigma_{\min}(\Mb_3)\geq \max\left\{
                    \frac{\sqrt{\sigma_{\min}(\Lb^-)} }{ 1+ \zeta_B },
            \frac{ \sigma_{\max}(\Bb_{\diag})}{2(1+\zeta_B)}\right\}.
            \end{align}

            Then, by inserting \eqref{eq: PM2} and \eqref{eq: M3} into \eqref{eq: M_2 min singular val 1}, we have
            \begin{align}
                \sigma_{\min}(\Mb_2)
                & \geq \frac{\max\left\{
                    \frac{\sqrt{\sigma_{\min}(\Lb^-)} }{ 1+ \zeta_B },
            \frac{ \sigma_{\max}(\Bb_{\diag})}{2(1+\zeta_B)}\right\}}
                    {3+\sqrt{\sigma_{\max}(\Lb^-)} }.
            \end{align}

\hfill$ \blacksquare$
\end{itemize}

\section{Proof of Lemma \ref{lem: block matrix singular value}}\label{app: proof of block matrix singular value}

            \begin{itemize}
                \item [(a)] For any vector $\vb$,
                    \begin{align}
                        \left\| \hat{\Cb} \vb\right\|
                        & = \left\| \begin{bmatrix}
                                    \Cb_1\vb\\
                                    \Cb_2\vb
                                \end{bmatrix}\right\|   \leq \|\Cb_1\vb\| + \|\Cb_1\vb\|   \leq \left(\sigma_{\max}(\Cb_1) + \sigma_{\max}(\Cb_2)\right) \|\vb\|.
                    \end{align}
                    Therefore, we have $\sigma_{\max}(\hat{\Cb}) \leq \sigma_{\max}(\Cb_1) + \sigma_{\max}(\Cb_2) $.
                \item [(b)] Notice that $\tilde{\Cb} =  [\Cb_3^\top, \Cb_4^\top]^\top $. Then $\sigma_{\max}(\tilde{\Cb}) = \sigma_{\max}([\Cb_3^\top, \Cb_4^\top]^\top)$ by the fact that $\sigma_{\max}(\Eb) = \sigma_{\max}(\Eb^\top) $ for any $\Eb$.
                    From \eqref{eq: block col singular val}, we have $  \sigma_{\max}([\Cb_3^\top, \Cb_4^\top]^\top) \leq \sigma_{\max}(\Cb_3)+\sigma_{\max}(\Cb_4)$.

                    Thus, we have $\sigma_{\max}(\tilde{\Cb}) \leq \sigma_{\max}(\Cb_3)+\sigma_{\max}(\Cb_4)$.  \hfill$ \blacksquare$
            \end{itemize}

\section{Proof of Lemma \ref{lem: sigma max M}}\label{app: proof of sigma max M}
The bounds for $\sigma_{\max}(\Mb_1)$ and $\sigma_{\max}(\Mb_2)$ can be easily obtained by applying Lemma 14.\hfill$ \blacksquare$

\section{Proof of Corollary 2}\label{app: proof of corollary 2}

Recalling the upper bound in \eqref{eqn: further bound},
 {\small
        \begin{align}
           \textstyle & \|\zb^t-\xb(\zb^t)\|^2 \notag\\& \leq
            \textstyle \frac{\Phi^0-\underline{\Phi}}{r}\textstyle \left(
              D_1 B_{\max}^2 \lambda_{\max} \sqrt{D_2  \theta_1^4 +D_{3} \theta_1^2  } \right.\notag \\
            &~~~\left.
             +
             \sqrt{D_4 + \textstyle  B_{\max}^2 \lambda_{\max}^2 (D_5 \theta_1^4 + D_6\theta_1^2)
            \textstyle+ \textstyle \frac{D_7 + D_8 B_{\max}^2\theta_2^2}{B_{\max}^2\lambda_{\max}^2(D_9 \theta_1^4 + D_{10}\theta_1^2) }} \right.\notag \\
            &~~~\left.\textstyle
                +\textstyle D_{11} B_{\max}\theta_1 \sqrt{ D_{12} \theta_3^2 + D_{13}  } \right.\notag \\
            &~~~ ~~~ ~~~\left.\times \max\bigg\{D_{14}, \sqrt{D_{15} \lambda_{\max}^2 \theta_1^2 + D_{16}}\bigg\}\right)^2,
        \end{align} }
 there are four terms to be analyzed, including $\lambda_{\max}$, $\theta_1$, $\theta_2$ and $\theta_3$.
 \begin{itemize}
    \item [a)] { \underline{\bf{$\lambda_{\max}$: }} }
Notice that the degree of each node in the cycle graph is always 2. Then, by \cite[Theorem 4.6]{cvetkovic2007signless}, we have $\lambda_{\max}\leq 2\max_{i\in [N]}d_i = 4$, where $d_i$ denotes the degree of node $i$. It means that the upper bound $\lambda_{\max} \leq 4$ is independent of network size $N$ and connectivity.

    \item [b)]  { \underline{\bf{$\theta_1$, $\theta_2$, $\theta_3$: }} }
By \cite[Lemma 3.7.1]{ringgraph}, the eigenvalues of $\Lb^-$ for the cycle graph are $2- 2\cos(2\pi k/N)$ with $k = 0, 1, \ldots, N/2$. Then, the largest eigenvalue is $\sigma_{\max}(\Lb^-) = 2$.

Since $2- 2\cos(2\pi k/N)$ is an monotonically increasing function with respect to $k$, then the smallest nonzero eigenvalue is attained at $k=1$. Therefore, we have
\begin{align}
    \sigma_{\min}(\Lb^-) =2- 2\cos\left(\frac{2\pi }{N}\right)
        = 2\left(\frac{\left(\frac{2\pi }{N}\right)^2}{2} + o\left(\left(\frac{2\pi }{N}\right)^2\right)\right)
= \frac{4\pi^2}{N^2} +  2\times o\left(\frac{4\pi^2 }{N^2}\right),
\end{align} where $o\left(\frac{4\pi^2 }{N^2}\right)$ is the high-order infinitely small quantities of $\frac{4\pi^2 }{N^2}$. Then, when $N$ increases, $\sigma_{\min}(\Lb^-)$ decreases.

Recall $\theta_1$, $\theta_2$, $\theta_3$ in \eqref{eq: theta_1} and \eqref{eq: theta_3},
\begin{align*}
    & \textstyle \theta_1=  \theta_2  \leq\textstyle
             (1+ \zeta_B)^2 \frac{\sigma_{\max}(\Lb^-)}{\sigma_{\min}^2(\Lb^-) }
             + \frac{4(1+ \zeta_B)^2}{ \sigma_{\max}(\Bb_{\mathrm{diag}})}  , \\
    &\textstyle \theta _3 \leq
            \textstyle \frac{\left(2\sqrt{\sigma_{\max}(\Lb^-)}+\sigma_{\max}(\Lb^-)\right)(3+\sqrt{\sigma_{\max}(\Lb^-)} )^2 (1+ \zeta_B)^2 }
            {\sigma_{\min}(\Lb^-) }
            \notag\\
            &~~~ ~~~ ~~~
            + \textstyle
                \frac{ 4(1+\zeta_B)^2}{ \sigma_{\max}^2(\Bb_{\mathrm{diag}})}\left(3+\sqrt{\sigma_{\max}(\Lb^-)} \right)^2 .
\end{align*}
        Obviously, after inserting $\sigma_{\max}(\Lb^-) = 2$, the decreases of $\sigma_{\min}(\Lb^-)$ leads to the rise of $\theta_1$, $\theta_2$, and $\theta_3$. Therefore, $\theta_1$, $\theta_2$, and $\theta_3$ increases with $N$.
 \end{itemize}
    By further inserting $\lambda_{\max}\leq 4$, one can conclude that the right hand side of  \eqref{eqn: further bound} increases with $N$.
   Since for the cycle graph,  the increase of $N$ implies poorer connectivity,  it implies that the convergence speed may be degraded when the network connectivity decreases.

\section{Proof of Theorem 2} \label{app: inexact}

For the inexact case, since for agent $i$ the inexact solution $\xb_i^\rpo$ is not equal to $\xb_i(\yb_i^\rpo, \zb_i^r)$ anymore, then we have different upper bounds for $\Lc_\rho$, $\tilde{G}$, and perturbation of $\yb$ from those in the exact PDC case.

\subsection{ Upper bound for $\Lc_\rho $ }
\vspace{-0.2cm}
\begin{Lemma}\label{fact: potential fun L inexact}
    For the AL function in \eqref{eqn: AL of DP2},we have	
    		\begin{align}\label{eq: potential function L inexact}
    		&-\Lc_\rho(\yb^{r+1},\mub^{r+1};\zb^{r+1}) + \Lc_\rho(\yb^r,\mub^r;\zb^r)
    \notag \\
    		& \leq   {   \sf (A1)}
    -\frac{1}{2} \|\yb^\rpo - \yb^r\|^2_{\rho\Lb^+ } -\frac{\rho}{2} \|\yb^\rpo -
    \yb^r\|^2
    		\notag \\
    		&~~~+ { \sf (A2)} + \frac{  B_{\max}^2}{2} \sum_{i=1}^N\|\xb_i^{r+1} -  \xb_i(\yb_i^{r+1} , \zb_i^r)\|^2 + \frac{1}{2}  \|\yb^\rpo - \yb^r\|^2.
    		\end{align}	
\end{Lemma}
{\bf Proof of Lemma \ref{fact: potential fun L inexact}:}
For proving the descent of $-\Lc_{\rho}$, the IPDC algorithm differs from its exact counterpart (Lemma \ref{fact: potential fun L}) only in the update of $\yb$, as shown below.

Bound of $-\Lc_\rho(\yb^{r+1},\mub^{r+1};\zb^{r}) +
\Lc_\rho(\yb^r,\mub^{r+1};\zb^r)$:
	Rewrite \eqref{eq: nabla y exact} explicitly, we have
	\begin{align}\label{eq: gradient of Lrho}
	\begin{bmatrix}
        \Bb_1 \xb_1(\yb_1^{r+1} , \zb_i^r)\\
        \vdots,\\
         \Bb_N \xb_N(\yb_N^{r+1} , \zb_i^r)
         \end{bmatrix}- \qb/N -\Ab^\top \mub - \rho\Ab^\top\Ab\yb -  \rho \Lb^+  (\yb^\rpo - \yb^r) =\zerob.
	\end{align}
    Then, by the optimality condition of \eqref{eqn: alg original form y} and applying \eqref{eq: gradient of Lrho}, we have
	\begin{align}
	  & \Lc_\rho(\yb^r,\mub^{r+1};\zb^r) - \frac{1}{2}\|\yb^r - \yb^r\|_{\rho \Lb^+
}^2  \notag \\
    &  \leq\sum_{i=1}^N \la \Bb_i (\xb_i^{r+1} -  \xb_i(\yb_i^{r+1} , \zb_i^r)), \yb_i^\rpo - \yb_i^r \ra  +  \Lc_\rho(\yb^\rpo,\mub^{r+1};\zb^r) - \frac{1}{2}\|\yb^\rpo - \yb^r\|_{\rho \Lb^+}^2 \notag\\
        &~~~ -\frac{\rho}{2}\| \yb^\rpo - \yb^r \|^2 \notag\\
    & \leq \frac{  B_{\max}^2}{2} \sum_{i=1}^N\|\xb_i^{r+1} -  \xb_i(\yb_i^{r+1} , \zb_i^r)\|^2 + \frac{1}{2}  \|\yb^\rpo - \yb^r\|^2
        +  \Lc_\rho(\yb^\rpo,\mub^{r+1};\zb^r)- \frac{1}{2}\|\yb^\rpo - \yb^r\|_{\rho \Lb^+}^2 \notag\\
        &~~~ -\frac{\rho}{2}\| \yb^\rpo - \yb^r \|^2
	\label{proof of fact: potential fun L 2 inexact}
	\end{align}

	Summing \eqref{proof of fact: potential fun L 1}, \eqref{proof of fact: potential fun L 2 inexact} and \eqref{proof of fact: potential fun L 3} leads to the desired result.
	\hfill $\blacksquare$

\subsection{ Upper bound for $\tilde{G} $}

\vspace{0.2cm}
\begin{Lemma}\label{fact: potential func G inexact}
We have
\begin{align}\label{eq: potential function G 2 inexact}
	& \tilde{G}_{r+1}(\xb^{r+1}, \yb^{r+1}, \mub^{r+1}; \zb^{r+1})
-\tilde{G}_r(\xb^r, \yb^r, \mub^r; \zb^r)
	\notag \\
	&\leq {-\frac{1}{\alpha}\|\mub^\rpo - \mub^r\|^2} - \sum_{i=1}^N
\frac{p+\gamma^-}{2} \|\xb_i^\rpo-\xb_i^r\|^2
-\frac{1}{2}\| (\yb^r-\yb^\rmo) - (\yb^\rpo-\yb^r) \|^2_{\rho \Lb^+ }
	\notag \\
	&~~~
-\bigg( \frac{ \rho}{2}- \frac{\alpha}{2}  \bigg) \|\yb^\rpo -
\yb^r\|^2_{\Lb^-}
+  \| \yb^\rpo - \yb^r\|^2_{\rho  \Lb^+ }
+
\frac{p}{2}(1-\frac{2}{\beta})\|\zb^\rpo
- \zb^r\|^2
\notag\\
&~~~ + \sum_{i=1}^N( -\frac{1}{\zeta} + \frac{p + \gamma^+}{2})\|\xb_i^{r+1}
-\xb_i^r\|^2.
	\end{align}
\end{Lemma}

{\bf Proof of Lemma \ref{fact: potential func G inexact}:}
The only difference
for proving the descent of $\tilde  G$ between the PDC algorithm (Lemma \ref{fact: potential func G}) and the IPDC algorithm is the bound of $G(\xb^{r+1}, \yb^{r}, \mub^{r+1}, \zb^{r}) - G(\xb^{r}, \yb^{r},
\mub^{r+1}, \zb^{r}) $:
    Define
$$H(\xb, \yb^r, \mub; \zb)=\max_{\yb}\bigg(G(\xb, \yb, \mub;
\zb)-\frac{\rho}{2}\|\yb-\yb^r\|^2_{\Lb^+}\bigg).$$

	Thus, we have $G(\xb^{r+1}, \yb^{r}, \mub^{r+1}, \zb^{r}) - G(\xb^{r}, \yb^{r},
\mub^{r+1}, \zb^{r}) $ to be bounded as
	\begin{align}
	&G(\xb^{r+1}, \yb^{r}, \mub^{r+1}, \zb^{r}) - G(\xb^{r}, \yb^{r}, \mub^{r+1},
\zb^{r}) \notag  \\
	&= \sum_{i=1}^N \bigg( f_i(\xb_i^\rpo)  +  \frac{p}{2}\|\xb_i^\rpo-\zb_i^r\|^2
-   f_i(\xb_i^r)  -  \frac{p}{2}\|\xb_i^r-\zb_i^r\|^2 \bigg)
	+  \sum_{i=1}^N  \la  \yb_i^r, \Bb_i ( \xb_i^\rpo -\xb_i^r)  \ra   \notag \\
    & \overset{\mathrm{(i)}}{\leq}  \sum_{i=1}^N \bigg(\big(\nabla f_i(\xb_i^r) + p(\xb_i^r - \zb_i^r)\big)^\top (\xb_i^{r+1} -\xb_i^r) + \frac{p +
    \gamma^+}{2} \|\xb_i^{r+1} -\xb_i^r\|^2 \bigg)
    +  \sum_{i=1}^N  \la   \yb_i^r, \Bb_i ( \xb_i^\rpo -\xb_i^r )  \ra
\notag \\
    & =  \sum_{i=1}^N \bigg(\nabla H(\xb_i^r)^\top (\xb_i^{r+1} -\xb_i^r)
    -\la  \Bb_i^\top \yb_i^{r+1}  ,   \xb_i^\rpo -\xb_i^r   \ra
    + \frac{p +
    \gamma^+}{2} \|\xb_i^{r+1} -\xb_i^r\|^2 \bigg)
    \notag \\
&~~~
     + \sum_{i=1}^N  \la   \yb_i^r, \Bb_i ( \xb_i^\rpo -\xb_i^r )  \ra
\notag \\
	&\overset{\mathrm{(ii)}}{\leq} \sum_{i=1}^N(-\frac{1}{\zeta}+ \frac{p + \gamma^+}{2}) \|\xb_i^{r+1} -\xb_i^r\|^2
 +  \sum_{i=1}^N  \la  \yb_i^{r+1} - \yb_i^r, \Bb_i ( \xb_i^r -\xb_i^\rpo)  \ra,
\label{proof of fact: potential func G 3 inexact}
	\end{align}
where $\mathrm{(i)}$ is due to the Lipschitz continuity of $\nabla F_i$, and $\mathrm{(ii)}$ is obtained by the updating step for $\xb$ in \eqref{eq: dual ADMM x2} of the IPDC algorithm.
\hfill $\blacksquare$
\subsection{ Error bounds and perturbation bounds}
The error bounds and perturbation bound in inexact case are analyzed in this
subsection.
    \begin{Lemma}\label{fact: inexact gap for x}
          Given $\zeta < \frac{2}{p+\gamma^+}$, we consider the upper bound of the solution difference between \eqref{eq: dual ADMM x1} and \eqref{eq: dual ADMM x2}.
        \begin{align}\label{x^r+1 - x(y^r+1, z^r)}
            \|\xb_i^\rpo - \xb_i(\yb_i^\rpo, \zb_i^r)\| \leq  \sigma_4   \|\xb_i^\rpo -
            \xb_i^r\|,
        \end{align}
        where $ \sigma_4   =  \bigg(1 + \frac{3}{\zeta(p + \gamma^-)}\bigg) $.
    \end{Lemma}
    {\bf Proof of Lemma \ref{fact: inexact gap for x}:}
The proof is similar to inequality (3.4) of Lemma 3.6 in \cite{Zhang18} hence is omitted.
		\hfill $\blacksquare$

	\begin{Lemma}\label{fact: perturb of y 2 inexact}
		Consider \eqref{eqn: dual function of DP2} with $p>-\gamma^-$. Then for any $  \bar \yb^r \in \mathcal{Y}(\mub^{r+1},\zb^{r})$,  we have
    \begin{enumerate}[(a)]
    \item
        \begin{align}
            \| \yb^{r+1} -  \bar \yb^r\|_{\Lb^-}^2
             \leq {a_1}
            \|\Lb^+(\yb^{r+1} - \yb^r)\|^2
        + a_5
         \|\xb^\rpo - \xb^r\|^2,
        \end{align}
    \item
        \begin{align} \label{perturb y: y^r+1 - y(mu^r+1, z^r)}
	     \|\yb^{r+1}-\bar \yb^r\|^2
     \leq {a_2}
        \|\Lb^+ (\yb^{r+1} - \yb^r)\|^2
     +  a_6
    \|\xb^\rpo - \xb^r\|^2.
    \end{align}
    \end{enumerate}
	\end{Lemma}

{\bf Proof of Lemma \ref{fact: perturb of y 2 inexact}(a):}
    For any $\bar \yb^r \in \mathcal{Y}(\mub^{r+1}, \zb^{r})$, recall its KKT conditions in \eqref{KKT y2 C1} and \eqref{KKT y2 C2}.

	Analogously, from \eqref{eqn: alg original form y}%
    , we have KKT conditions for $\yb^{r+1}$ as
		\begin{align}
	& \begin{bmatrix}
	\Bb_1 \xb_1(\yb_1^{r+1},\zb_1^{r}) -\frac{\qb}{N}+ \Bb_i(\xb_1^\rpo -
\xb_1(\yb_1^{r+1},\zb_1^{r}))   \\
	\vdots \\
	\Bb_N \xb_N(\yb_N^{r+1},\zb_N^{r}) - \frac{\qb}{N} + \Bb_N(\xb_N^\rpo -
\xb_N(\yb_N^{r+1},\zb_N^{r}))
	\end{bmatrix}-\Ab^\top \mub^{r+1} - \rho \Lb^- \yb^{r+1}
- \rho \Lb^+ (\yb^\rpo -\yb^r) =\zerob,
\label{KKT y2 r+1 C1 inexact} \\
	& \nabla  F_i ( \xb_i(\yb_i^{r+1},\zb_i^r) ) + \Bb_i^\top \yb_i^{r+1}
=\zerob, \forall i\in [N].
	 \label{KKT y2 r+1 C2 inexact}
	\end{align}

With fixed $\{\xb_i\big(\bar \yb_i^r,\zb_i^{r}\big)\}_{i=1}^N$,
\eqref{KKT y2 C1} and \eqref{KKT y2 C2} forms a linear system of
	variable $\rho \bar \yb^r$. Then, by using the Hoffman bound in Lemma \ref{Hoffman}, the distance from $\yb^{r+1}$ satisfying \eqref{KKT y2 r+1 C1 inexact} and \eqref{KKT y2 r+1 C2 inexact} to the set $\Yc(\mub^{r+1},\zb^{r})$ can be bounded as
	\begin{align}
	    &\rho^2 {\rm dist}^2(\yb^{r+1},\Yc(\mub^{r+1},\zb^{r})) \notag \\
	    &\leq
            \theta_1^2 \sum_{i=1}^N  \| \Bb_i ( \xb_i(\yb_i^{r+1},\zb_i^r) ) -
            \xb_i(\bar \yb_i^{r},\zb_i^{r}) ) \|^2
            +\theta_1^2 \rho^2 \|\Lb^+(\yb^{r+1} - \yb^r)\|^2
        \notag \\
	    &~~~+
            \theta_1^2 \rho^2 \sum_{i=1}^N  \|
            \nabla F_i (\xb_i(\yb_i^{r+1},\zb_i^r))
                - \nabla F_i ( \xb_i(\bar \yb_i^{r},\zb_i^{r}) )
            \|^2
          + \theta_1^2 \sum_{i=1}^N  \| \Bb_i ( \xb_i^\rpo -
\xb_i(\yb_i^{r+1},\zb_i^r) ) \|^2
    \notag\\
     & \leq
            \frac{\theta_1^2}{2\eta_1} \bigg( \frac{1 }{p + \gamma^-} + \rho^2 (p + \gamma^+) \bigg)
            \| \yb^{r+1} -  \bar \yb^{r}\|^2
            + \bigg[\frac{\eta_1}{2}  \bigg( \frac{1 }{p + \gamma^-} + \rho^2 (p + \gamma^+) \bigg)+1\bigg]\theta_1^2\rho^2 \|\Lb^+(\yb^{r+1} - \yb^r)\|^2
        \notag \\
            &~~~
                - \bigg( \frac{\theta_1^2 }{p + \gamma^-} +\theta_1^2 \rho^2 (p + \gamma^+) \bigg)\rho
                \| \yb^{r+1} - \bar \yb^{r}\|_{\Lb^-}^2
    +  \theta_1^2 B_{\max}^2  \sigma_4 ^2
    \|\xb^\rpo - \xb^r\|^2,
    \label{due error bound y2 C3 inexact}
    \end{align}
	where the last inequality is due to \eqref{y^r+1 - y(mu^r+1, z^r) 1} and
\eqref{x^r+1 - x(y^r+1, z^r)}.

     Analogous to Lemma \ref{fact: perturb of y 2}, choose $\eta_1 = \textstyle \frac{\theta_1^2 }{\rho^2}\big( \frac{1 }{p + \gamma^-} + \rho^2 (p + \gamma^+) \big) $, we can obtain

    \begin{align}
         \| \yb^{r+1} -  \bar \yb^{r}\|_{\Lb^-}^2 \leq  a_1  \|\Lb^+(\yb^{r+1} - \yb^r)\|^2 +
        \underbrace{
        \frac{ B_{\max}^2 \sigma_4 ^2 }{ \left(\frac{1 }{p + \gamma^-} + \rho^2 (p + \gamma^+)\right)\rho }
        }_{\triangleq a_5}
        \|\xb^\rpo - \xb^r\|^2.
    \end{align}
 {\bf Proof of lemma \ref{fact: perturb of y 2 inexact}(b):}
 By  lemma \ref{fact: perturb of y 2} and \eqref{due error bound y2 C3 inexact},
 we have
 \begin{align}\label{eb1 inexact}
	    \|\yb^{r+1}-\bar \yb^{r}\|^2
     & \leq a_2
        \|\Lb^+ (\yb^{r+1} - \yb^r)\|^2
     + \underbrace{\frac{  2 \theta_1^2 B_{\max}^2 \sigma_4 ^2 }{\rho^2}}_{\triangleq a_6}\|\xb^\rpo - \xb^r\|^2 .
    \end{align}
    By taking the square root on both sides, we have
    \begin{align}
	    \|\yb^{r+1}-\bar \yb^{r}\|
     & \leq \sqrt{  a_2
        \|\Lb^+ (\yb^{r+1} - \yb^r)\|^2 + a_6 \|\xb^\rpo - \xb^r\|^2 } \notag \\
     & \leq  \sqrt{  a_2} \|\Lb^+ (\yb^{r+1} - \yb^r)\| + \sqrt{a_6}\|\xb^\rpo - \xb^r\|.
    \end{align}
	\hfill $\blacksquare$

	\begin{Lemma}\label{fact: perturb of y inexact}
		Consider $p>\gamma^-$, based on Lemma \ref{fact: perturb of y} and Lemma \ref{fact: perturb of y 2 inexact}, we have
        \begin{align}
             \| \yb^{r+1} -  \bar \yb^{r+1}\|^2 \leq 2 a_2 \|\Lb^+(\yb^{r+1} - \yb^r)\|^2 + 2a_6 \|\xb^\rpo -
            \xb^r\|^2
            +2 a_3\|\zb^{r+1} - \zb^r\|^2.
        \end{align}
	\end{Lemma}

{\bf Proof of Lemma \ref{fact: perturb of y inexact}:} We have
\begin{align}
             \| \yb^{r+1} -  \bar \yb^{r+1}\|^2 &= \| \yb^{r+1} - \yb(\mub^{r+1},\zb^{r}) + \bar \yb^{r}-  \bar \yb^{r+1}\|^2 \notag \\
             &~~~\leq 2 \| \yb^{r+1} - \bar \yb^{r}\|^2 + 2\|\bar \yb^{r}-  \bar \yb^{r+1}\|^2 \notag \\
            &~~~\leq 2 a_2 \|\Lb^+(\yb^{r+1} - \yb^r)\|^2 + 2a_6 \|\xb^\rpo -
            \xb^r\|^2
             +2 a_3\|\zb^{r+1} - \zb^r\|^2 ,
        \end{align}
where the last inequality is due to \eqref{perturb y: y^r+1 - y(mu^r+1, z^r)} and Lemma \ref{fact: perturb of y}.

\hfill $\blacksquare$

%
\vspace{0.15cm}
Then, we can bound the following differences.
\vspace{-0.15cm}
\begin{Lemma}\label{eb inexact}
    	For any $\bar \yb^{r}\in \mathcal{Y}(\mub^{r+1}, \zb^r)$, there
    exist $\yb(\zb^r) \in \mathcal{Y}(\zb^r)$  and $\xb(\zb^r)$ such that
    	\begin{align}
        	&\|\yb^{r+1}-\yb(\zb^r)\| \leq \sqrt{  a_2} \lambda_{\max}\|   \yb^{r+1} - \yb^r \| + \sqrt{a_6}\|\xb^\rpo - \xb^r\| +  \sqrt{a_4} \|\Ab \bar \yb^{r}\| , \label{eb2
        inexact}\\
        	&\|\xb^{r+1}-\xb(\zb^r)\| \leq \sigma_2\|\yb^{r+1}-\yb(\zb^r)\| +
          \sigma_4   \|\xb^\rpo - \xb^r\| .\label{eb3 inexact}
    	\end{align}
    \end{Lemma}

{\bf Proof of Lemma \ref{eb inexact}:}
    Firstly,
    \begin{align} \label{eb2-1}
    \|\yb^{r+1}-\yb(\zb^r)\|
    & = \|\yb^{r+1}- \bar \yb^{r} + \bar \yb^{r} - \yb(\zb^r)\|
    \notag \\
          & \leq \|\yb^{r+1}-\bar \yb^{r}\| + \| \bar \yb^{r} -
          \yb(\zb^r)\|
    \notag \\
          & {\leq}   \sqrt{  a_2} \lambda_{\max}\|   \yb^{r+1} - \yb^r \| + \sqrt{a_6}\|\xb^\rpo - \xb^r\| +  \sqrt{a_4} \|\Ab \bar \yb^{r}\|,
    \end{align}
    where the last inequality is due to \eqref{eb1 inexact} and \eqref{y(z^r) - y(mu^{r+1, z^r})}.

    Secondly, we can have
   \begin{align}
    \|\xb^{r+1}-\xb(\zb^r)\|
        & = \|\xb^{r+1} - \xb(\yb^{r+1},\zb^r)+\xb(\yb^{r+1},\zb^r) - \xb(\zb^r)\| \notag \\
        & \leq  \|\xb^{r+1} - \xb(\yb^{r+1},\zb^r)\|+ \|\xb(\yb^{r+1},\zb^r) - \xb(\yb(\zb^r), \zb^r)\|\notag \\
        &
        \leq  \sigma_4  \| \xb^{r+1} - \xb^r \| + \sigma_2 \| \yb^{r+1} - \yb(\zb^r) \|,
   \end{align}
   where the last inequality is due to \eqref{eq: perturb of xi}.
\hfill $\blacksquare$

\vspace{0.15cm}
Similar to the PDC algorithm (Lemma \ref{Fact descent potential}), we also have the non-increasing property of the potential
function for the IPDC algorithm, which is given in the following lemma.

	\begin{Lemma}\label{Fact descent potential inexact}
	Let $ \beta $ satisfies \eqref{eqn: beta condition} in Lemma \ref{Fact descent potential} and $\rho > 1/2$.
Then given $p$ meeting
\begin{align}
 p > \max\bigg\{0, \gamma^+ - 2 \gamma^-, \frac{32 ( 2\rho-1)B_{\max}^2}{5\lambda_{\max}^2 \rho^2} - \gamma^-, 32B_{\max}^2 - \gamma^-, \frac{2\gamma^+ - 7\gamma^-}{5}\bigg\},
\label{eq: IPDC p}
\end{align}
for sufficiently small $\zeta$ and $\alpha$ satisfying
\begin{align}
& \frac{1}{p+\gamma^-} < \zeta < \min \bigg\{ \frac{5}{2(\gamma^+- \gamma^-)},~\frac{2}{p+\gamma^+},~
                    \frac{5a_2 \lambda_{\max}^2 \rho^2}{64(2\rho - 1)\theta_1^2 B_{\max}^2 },~
                    \frac{1}{32B_{\max}^2}\bigg\} ,\label{eq: IPDC zeta}\\
&  \alpha <\min \bigg\{ \frac{\rho}{5},~ \frac{2\rho-1}{16a_1
            \lambda_{\max}^2}, ~ \frac{1}{8a_5\zeta}\bigg\},\label{eq: IPDC alpha}
\end{align}
 and for any $r$,
\begin{enumerate}[(a)]
    \item we have
	\begin{align}
	 \Phi^\rpo -\Phi^r
	& \leq  -C_1
\|\yb^\rpo - \yb^r\|^2_{\Lb^-}
-\frac{1}{2}\| (\yb^r-\yb^\rmo) - (\yb^\rpo-\yb^r) \|^2_{\Lb^+ }
 -\alpha \| \Ab \yb(\mub^\rpo,\zb^r)\|^2
\notag \\
	&~~~
    -\bigg(C_2 -\frac{1}{2} \bigg)  \|\yb^\rpo - \yb^r\|^2 - C_3  \|\zb^\rpo - \zb^r\|^2
- C_4  \|\xb^\rpo - \xb^r\|^2 .\label{eq: Phi descent 2 inexact}
	\end{align}
    where $C_1, C_2,C_3  $ is defined in Lemma \ref{Fact descent potential},
    $C_4$ is
    \begin{align}\label{eq: C4}
     C_4 \triangleq    \frac{1}{\zeta} - \frac{ \gamma^+ - \gamma^-}{2}  -2\alpha a_5 - \frac{4p\sigma_2^2 a_6}{\delta}  -\frac{B_{\max}^2 \sigma_4^2}{2}  ,
    \end{align}
    and $C_1, (C_2-1/2),C_3, C_4 >0  $.
    \item
		The function $\Phi^r$ in \eqref{eq: potential fun} is lower bounded, i.e.,
$\sum_{r=0}^T \Phi^{r+1} > \underline{\Phi}>-\infty$ for some constant $\underline{\Phi}$.
	\end{enumerate}
    It implies that $\Phi^r$ is non-increasing with $r$.
\end{Lemma}

{\bf{ Proof of Lemma \ref{Fact descent potential inexact}}(a):}
This proof is divided into the following three steps:
 \begin{itemize}
 \item {\bf{Step 1:}} Obtain the bound in \eqref{eq: Phi descent 2 inexact};
 \item {\bf{Step 2:}} Solve sufficient conditions on $p$, $\rho$, $\beta$, $\zeta$, $\alpha$ for proving positiveness of $C_1-C_4$;
 \item {\bf{Step 3:}}  Summarize all the conditions of the convergence for the IPDC algorithm.
\end{itemize}
\begin{enumerate}[Step 1:]
    \item Obtain  \eqref{eq: Phi descent 2 inexact}:
    We can follow the same proof procedure of Lemma \ref{Fact descent potential}(a). Compared with \eqref{eq: Phi descent 2}, there are two additional terms in \eqref{eq: Phi descent 2 inexact} shown as follows.
         \begin{enumerate}[(1)]
         \item The first one is $-\frac{1}{2}  \|\yb^\rpo - \yb^r\|^2$, which is obtained by \eqref{eq: potential function L inexact} in Lemma \ref{fact: potential fun L inexact}.
         \item The second part $- C_4  \|\xb^\rpo - \xb^r\|^2$ consists of five terms. Note that the first two term is obtained by the bound in Lemma \ref{fact: potential func G inexact}, the third term is by applying lemma \ref{fact: perturb of y 2 inexact}(a) into \eqref{eq: combined term 1}${\mathrm{(iii)}}$, the fourth term is due to the use of Lemma \ref{fact: perturb of y inexact} in \eqref{eq: combined term 2}${\mathrm{(i)}}$, and the last term is resulted from the bound in Lemma \ref{fact: potential fun L inexact}.
        \end{enumerate}
    \item Prove the positiveness of $C_1$-$C_4$ in \eqref{eq: Phi descent 2 inexact}:
        \begin{enumerate}[(1)]
           \item Prove $C_1, C_3 > 0$: According to Lemma \ref{Fact descent potential}(a), we obtain that if $\alpha< \frac{\rho}{5}$ and $\beta$ satisfies \eqref{eqn: beta condition}, then $C_1$ and $ C_3 $ in \eqref{eq: bound for Phi 2} are both positive.
           \item Prove $C_2 >\frac{1}{2}$:
                As for $C_2$ in \eqref{eq: bound for Phi 2}, we firstly assume $\rho >\frac{1}{2}$. Different from the choice of $\delta$ in the proof of Lemma \ref{Fact descent potential}(a), we choose $\delta = \frac{{20} p\sigma_2^2 a_2 \lambda_{\max}^2}{\rho-\frac{1}{2}} >\frac{16 p\sigma_2^2 a_2 \lambda_{\max}^2}{\rho-\frac{1}{2} }$, i.e., $ \frac{4 p\sigma_2^2a_2\lambda_{\max}^2 }{\delta } < \frac{\rho-\frac{1}{2} }{4}$. Then, as long as  $ \alpha  < \frac{\rho -\frac{1}{2}}{8a_1 \lambda_{\max}^2}$ and $\rho > \frac{1}{2}$, we can have $C_2 -\frac{1}{2}> \frac{\rho}{2} - \frac{1}{4} > 0$ .
           \item Prove $C_4>0$:
                By inserting $\delta = \frac{{20} p\sigma_2^2 a_2 \lambda_{\max}^2}{\rho-1/2}$ into \eqref{eq: C4}, we have
                \begin{align}
                    C_4 =    \frac{1}{\zeta} - \frac{ \gamma^+ - \gamma^-}{2}  -2\alpha a_5 - \frac{(\rho-1/2)  a_6}{5a_2\lambda_{\max}^2 } -\frac{B_{\max}^2 \sigma_4^2}{2} .
                \end{align}

                To guarantee $C_4>0$, we have
                \begin{align}
                \alpha & < \frac{1}{2a_5} \bigg[
                            \frac{10 a_2 \lambda_{\max}^2 - (\gamma^+- \gamma^-)5 \zeta a_2 \lambda_{\max}^2 - 2(\rho -1/2) a_6 \zeta - 5B_{\max}^2 \sigma_4^2 \zeta a_2 \lambda_{\max}^2}
                            {10\zeta a_2 \lambda_{\max}^2}
                        \bigg]
                    \notag \\
                    &=\frac{10 a_2 \lambda_{\max}^2 - (\gamma^+- \gamma^-)5 \zeta a_2 \lambda_{\max}^2 - 2(\rho -1/2) a_6 \zeta - 5B_{\max}^2 \sigma_4^2 \zeta a_2 \lambda_{\max}^2}
                            {20\zeta a_2a_5 \lambda_{\max}^2}
                            \notag \\
                    & \triangleq   \frac{ D_1 }{  D_2  } .
                \end{align}

                Note that if $D_1>0$ and $\alpha <  \frac{D_1}{D_2}$, then we can complete the proof. To guarantee $D_1> 0$, it suffices to show the following three conditions
                \begin{subequations}
                    \begin{align}
                        & (\gamma^+- \gamma^-) \zeta a_2 \lambda_{\max}^2 < \frac{5a_2 \lambda_{\max}^2}{2} , \label{eq: C4 1}\\
                        & 2(\rho-\frac{1}{2}) a_6\zeta < \frac{5a_2 \lambda_{\max}^2}{2}, \label{eq: C4 2}\\
                        & B_{\max}^2 \sigma_4^2 \zeta a_2 \lambda_{\max}^2
                            < \frac{ a_2 \lambda_{\max}^2}{2}. \label{eq: C4 3}
                    \end{align}
                \end{subequations}

                Next, we prove these conditions one by one.
                \begin{enumerate}[1)]
                    \item To guarantee  \eqref{eq: C4 1} holds, i.e., $(\gamma^+- \gamma^-) \zeta a_2 \lambda_{\max}^2 < \frac{5a_2 \lambda_{\max}^2}{2}$ we require an upper bound on $\zeta$, i.e.,
                        \begin{align}\label{eq: C4 1 1}
                        \zeta < \frac{5}{2(\gamma^+- \gamma^-)}.
                        \end{align}
                    \item To guarantee \eqref{eq: C4 2} holds, i.e., $2(\rho-\frac{1}{2}) a_6\zeta < \frac{5a_2 \lambda_{\max}^2}{2}$, we obtain another bound on $\zeta$, i.e.,
                        \begin{align}\label{eq: C4 2 1}
                        \zeta < \frac{5a_2 \lambda_{\max}^2 \rho^2}{8(\rho - \frac{1}{2})\theta_1^2 B_{\max}^2 \sigma_4^2}.
                        \end{align}
                        Note that this bound is implicit, since $\sigma_4$ defined in Lemma \ref{fact: inexact gap for x} is a function of $\zeta$.

                         To overcome this, we assume $   \frac{1}{p+\gamma^-}< \zeta< \frac{2}{p+\gamma^+}$ and $p > \gamma^+ - 2 \gamma^-$. Then, we have $\sigma_4 =  \big(1 + \frac{3}{\zeta(p + \gamma^-)}\big) <4$. Therefore, we obtain an explicit sufficient condition of \eqref{eq: C4 2 1}, i.e.,
                        \begin{align}\label{eq: C4 2 2}
                             \frac{1}{p+\gamma^-}<  \zeta< \frac{5a_2 \lambda_{\max}^2 \rho^2}{128(\rho -\frac{1}{2})\theta_1^2 B_{\max}^2 }.
                        \end{align}

                        Since $\zeta$ has a lower bound $\frac{1}{p+\gamma^-}$, we need to guarantee that the solution set of $\zeta$ in \eqref{eq: C4 2 2} is non-empty. That is, the lower bound is required to be smaller than the upper bound. Inserting $a_2$ of \eqref{eq: a_2} into \eqref{eq: C4 2 2}, we have
                        \begin{align}
                            & \frac{1}{p+\gamma^-}
                              <  \frac{ 5\lambda_{\max}^2 \rho^2}{128(\rho - 1/2)\theta_1^2 B_{\max}^2 \sigma_4^2}
                            \bigg[\frac{\theta_1^4}{\rho^2} \bigg( \frac{1}{p+\gamma^-} + \rho^2 (p+\gamma^+)\bigg)^2 + 2\theta_1^2\bigg] \notag \\
                            \Leftarrow&~
                            5\lambda_{\max}^2 \theta_1^2 + 5\lambda_{\max}^2 \theta_1^2 \rho^4(p+\gamma^-)^2 (p+\gamma^+)^2  + 10\lambda_{\max}^2 \theta_1^2\rho^2(p+\gamma^-)  (p+\gamma^+)\notag\\
                            &~~~
                            + 10\lambda_{\max}^2  \rho^2(p+\gamma^-)^2
                            - 128( \rho-1/2) B_{\max}^2 (p+\gamma^-) >0 \notag\\
                            \Leftarrow &~10\lambda_{\max}^2  \rho^2(p+\gamma^-)^2
                            - 128( \rho-1/2)B_{\max}^2 (p+\gamma^-) >0 \notag\\
                            \Leftarrow &~ p> \frac{32 ( 2\rho-1 )B_{\max}^2}{5\lambda_{\max}^2 \rho^2} - \gamma^-,
                        \end{align}
                        which gives a sufficient condition for achieving \eqref{eq: C4 2 2}.

                        \item To guarantee \eqref{eq: C4 3} holds, i.e., $B_{\max}^2 \sigma_4^2 \zeta a_2 \lambda_{\max}^2
                            < \frac{ a_2 \lambda_{\max}^2}{2}$, we use the same operation on $\sigma_4$ in $\mathrm{(ii)}$, we get
                        \begin{align}\label{eq: C4 3 1}
                        \frac{1}{p+\gamma^-}<\zeta < \frac{1}{32B_{\max}^2}.
                        \end{align}
                        It is obvious that when $p > 32B_{\max}^2 - \gamma^-$, \eqref{eq: C4 3 1} is feasible.
                \end{enumerate}

                Note that for $\zeta < \frac{5}{2(\gamma^+- \gamma^-)}$ in \eqref{eq: C4 1 1} to hold, we also require  $\frac{1}{p+\gamma^-}< \frac{5}{2(\gamma^+- \gamma^-)}$,
                , which gives a condition of $ p> \frac{2\gamma^+ - 7\gamma^-}{5}. $
           \end{enumerate}

    \item
    Thus, combining all the results above above for proving positiveness of $C_1$-$C_4$, we have the following conditions on $\beta$, $\zeta$ and $\alpha$
        \begin{align}
            & \beta <  \bigg(\frac{1}{2} + \frac{20 p\sigma_2^2 a_2
    \lambda_{\max}^2}{\rho} + \frac{\rho(\sigma_1^2 + 2\sigma_2^2a_3)}{10
    p\sigma_2^2 a_2 \lambda_{\max}^2} \bigg)^{-1} ~~({\mathrm{Recall}}~\eqref{eqn: beta condition}),\\
            & \frac{1}{p+\gamma^-} < \zeta < \min \bigg\{ \frac{5}{2(\gamma^+- \gamma^-)},~\frac{2}{p+\gamma^+},~
                    \frac{5a_2 \lambda_{\max}^2 \rho^2}{64(2\rho - 1)\theta_1^2 B_{\max}^2 },~
                    \frac{1}{32B_{\max}^2}\bigg\} \\
                    &  \alpha <\min \bigg\{ \frac{\rho}{5},~ \frac{2\rho-1}{16a_1
            \lambda_{\max}^2}, ~ \frac{1}{8a_5\zeta}\bigg\},
        \end{align}
        provided that $\rho > \frac{1}{2}$ and $p$ satisfies
        \begin{align}
            p > \max\bigg\{0, \gamma^+ - 2 \gamma^-, \frac{32 ( 2\rho-1)B_{\max}^2}{5\lambda_{\max}^2 \rho^2} - \gamma^-, 32B_{\max}^2 - \gamma^-, \frac{2\gamma^+ - 7\gamma^-}{5}\bigg\}.
        \end{align}

    \end{enumerate}
    Therefore, given sufficiently large $p$ and $\rho >\frac{1}{2}$, for sufficiently small $\beta$, $\zeta$, and $\alpha$, the IPDC algorithm can have the same convergence properties as the PDC algorithm.

{\bf{ Proof of Lemma \ref{Fact descent potential inexact}}(b):}
The proof is similar to the proof of Lemma \ref{Fact descent potential}(b) (see Appendix \ref{appendix: lemma 8b}) hence is omitted.
\hfill $\blacksquare$

\vspace{0.2cm}
Similar to the proof of Theorem \ref{thm: conv and iteration comp}, we can also obtain that every limit point of the iteration sequence $\{(\zb^r, \yb^r)\}$ is a KKT solution of problem {\sf (P)}.
There exists a constant $\hat{C}>0$ and such that \begin{align}
        \|\zb^r-\xb(\zb^r)\|^2 \leq  \frac{\hat C}{t},
    \end{align}
where $
    {\hat C} \triangleq (\Phi^0-\underline{\Phi})\left(\frac{1}{\beta \sqrt{C_3}}+
    \frac{\sigma_2\sqrt{a_2}\lambda_{\max}}{\sqrt{C_2}} + \frac{\sigma_2\sqrt{a_6} +
      \sigma_4  }{\sqrt{C_4}} + \frac{\sigma_2\sqrt{a_4}}{\sqrt{\alpha}}\right)^2.$\hfill $\blacksquare$
\section{Additional Simulation Results} \label{Additional Numerical-LR}

\begin{figure}[!t]
\begin{minipage}[b]{1.0\linewidth}
  \centering
 \epsfig{figure=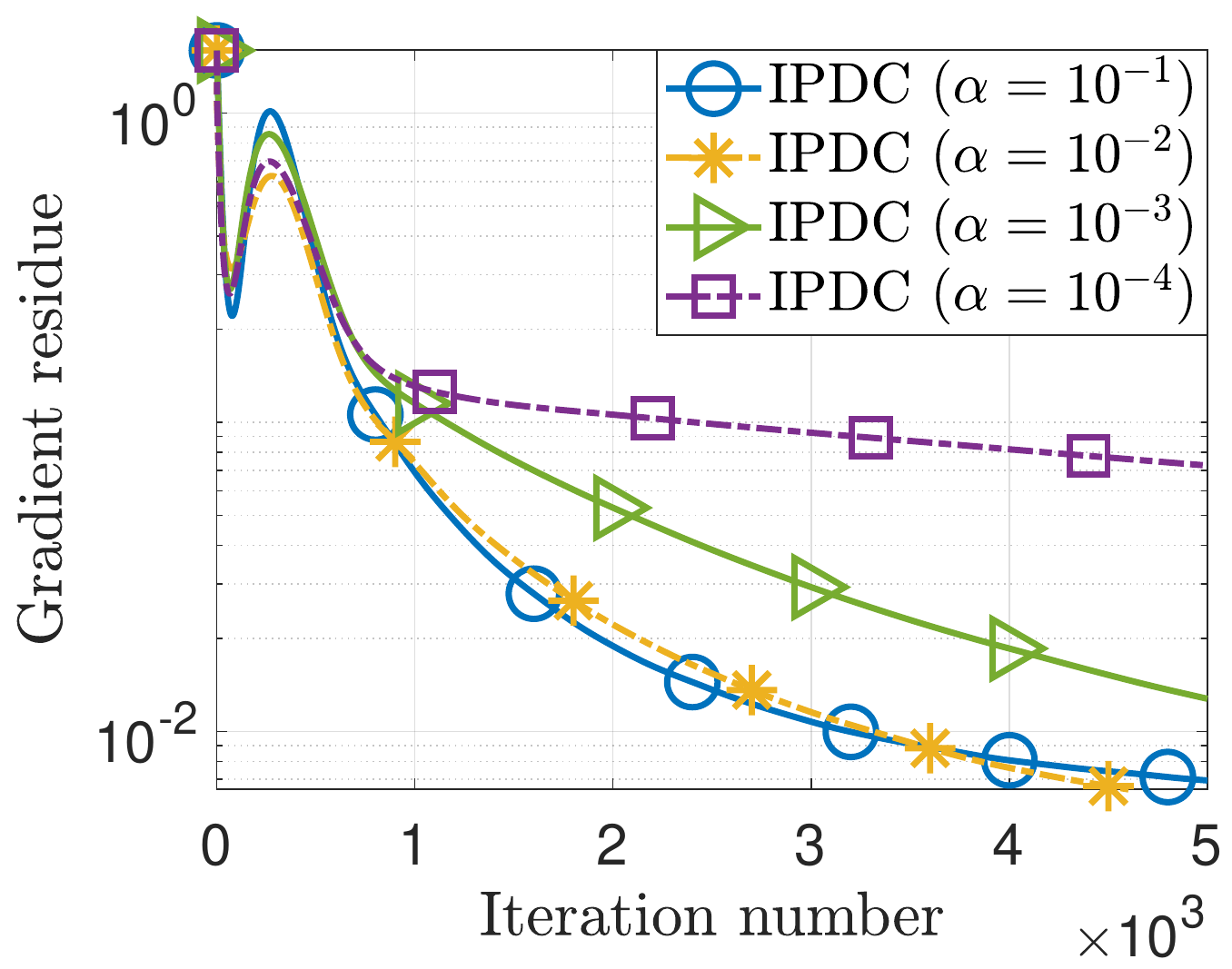,width=2.3in}
 \epsfig{figure=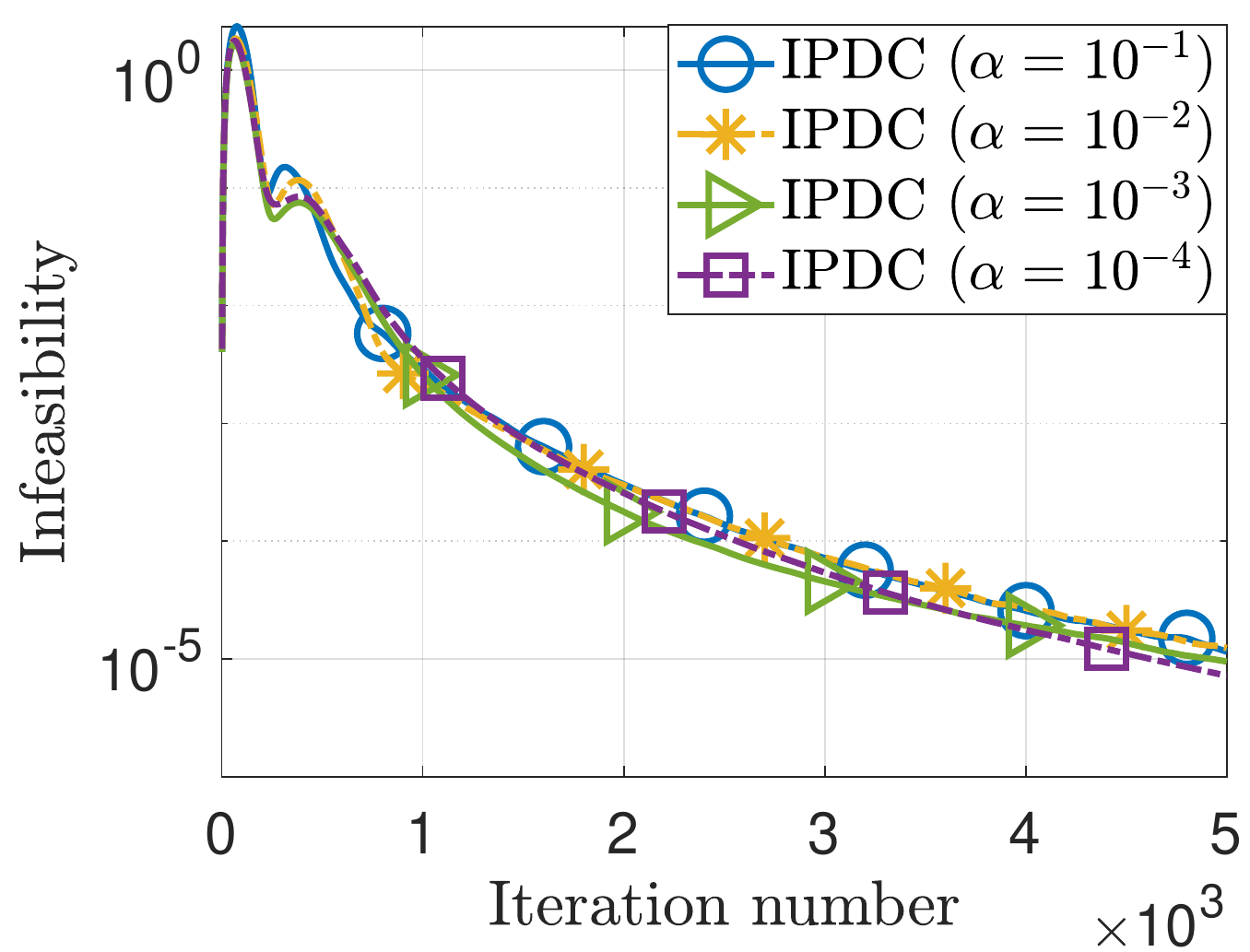,width=2.3in}
  \centerline{\small{(a) $\beta = 0.1$, $ p = 10$, $\rho = 1$, $\zeta = 0.1$, and various values of $\alpha$}}\medskip
\end{minipage}
\begin{minipage}[b]{1.0\linewidth}
  \centering
 \epsfig{figure=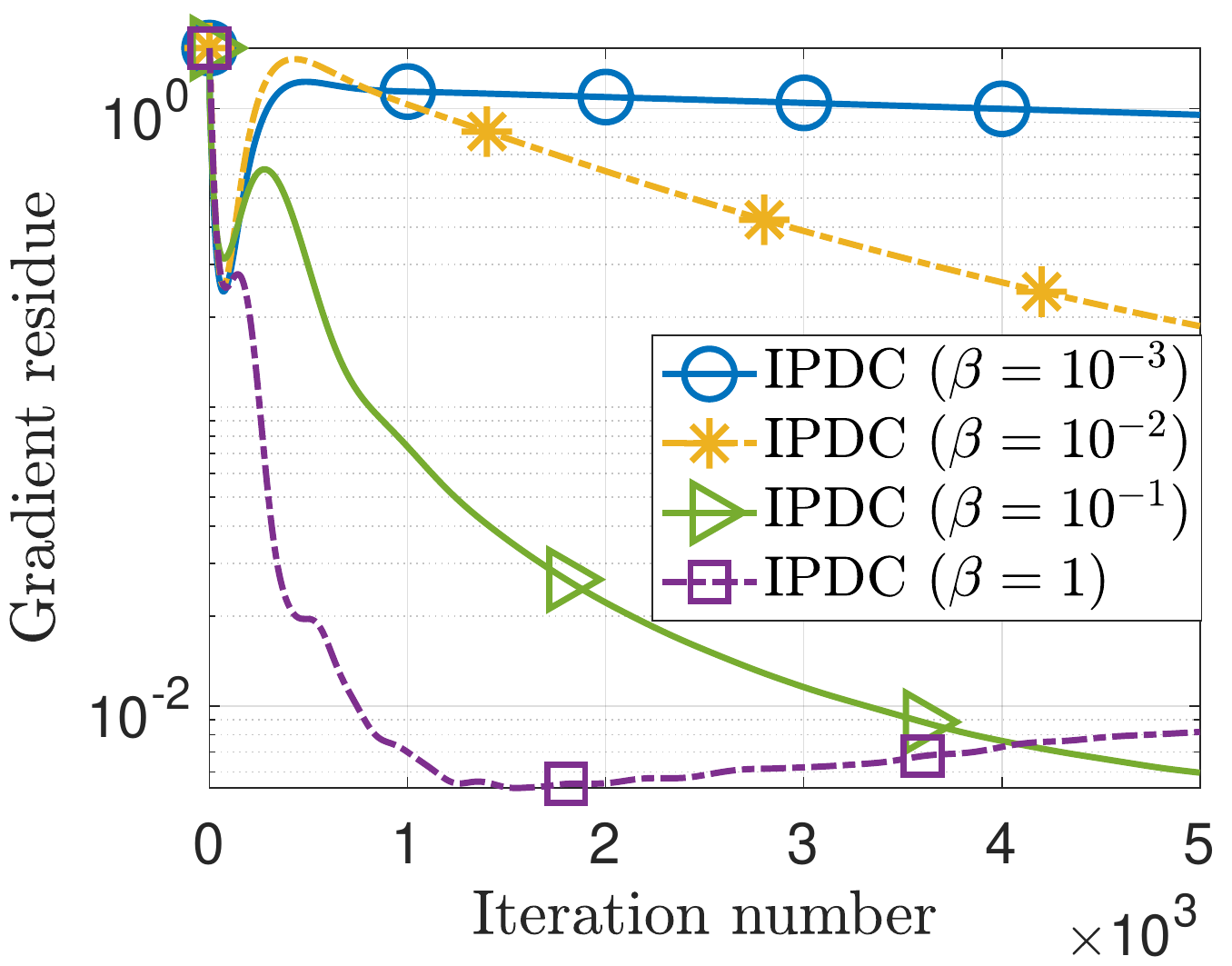,width=2.3in}
 \epsfig{figure=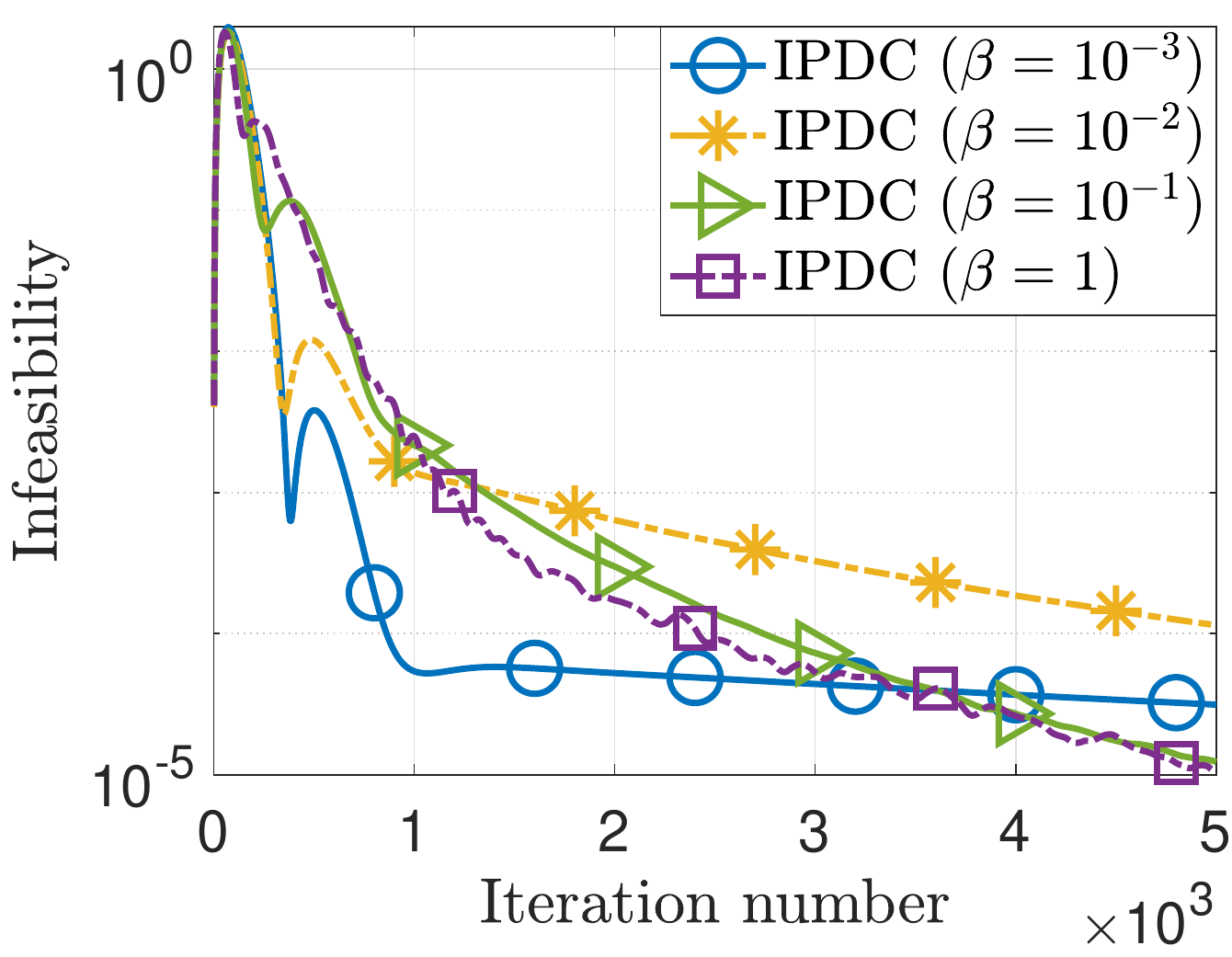,width=2.3in}
  \centerline{\small{(b) $\alpha = 0.01$, $\rho = 1$, $\zeta = 0.1$, $p = 10$, and various values of $\beta$}}\medskip
\end{minipage}
\caption{Convergence curves of the IPDC algorithm different step-sizes.} \label{fig: LR IPDC alpha beta}\vspace{-0.4cm}
\end{figure}

In this section, besides of the results displayed in the main text, we show the convergence results for the proposed IPDC algorithm with respect to the parameter $\alpha$, $\beta$, $\zeta$, $p$ and $\rho$.

\begin{figure}[!t]
\begin{minipage}[b]{1.0\linewidth}
  \centering
 \epsfig{figure=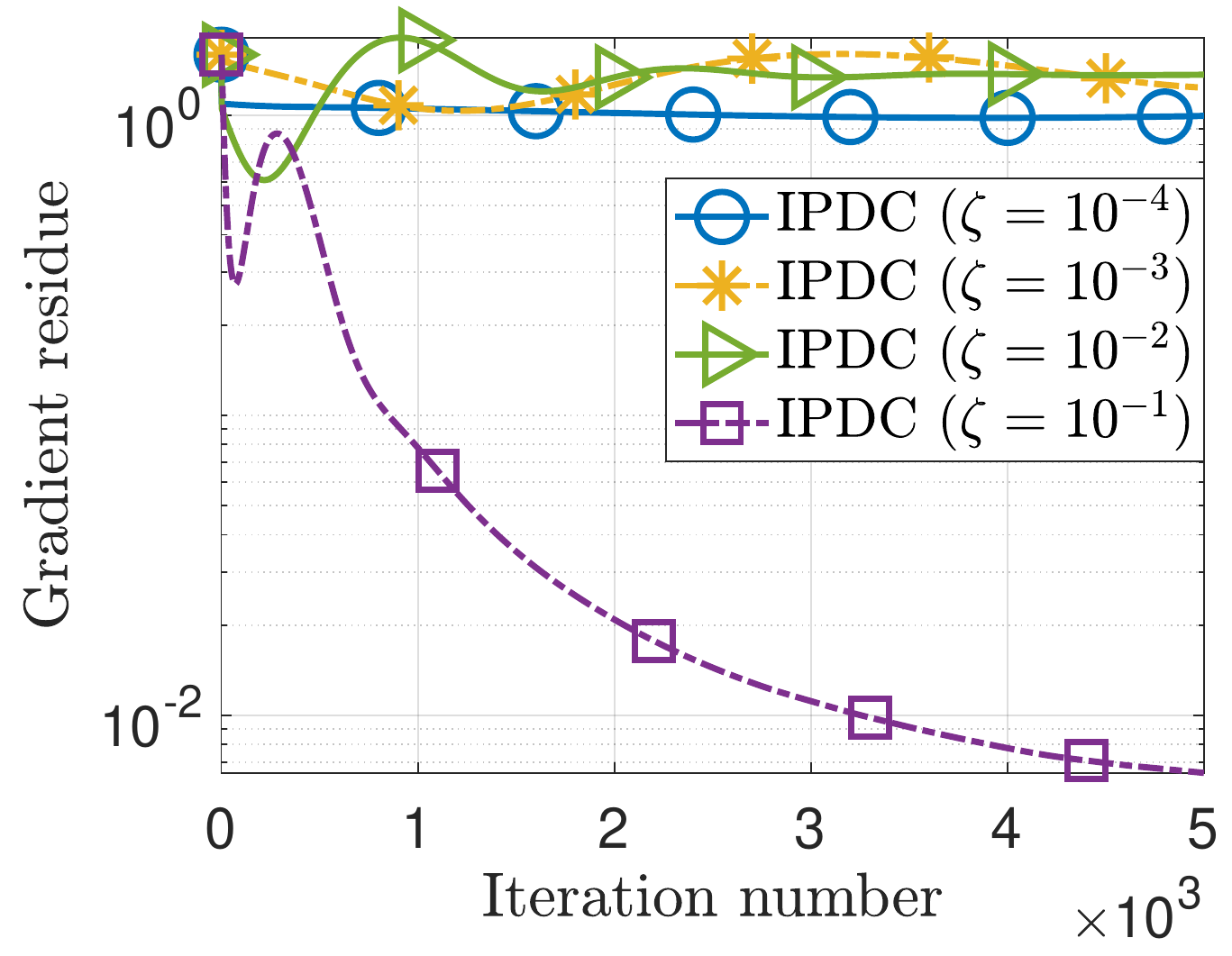,width=2.3in}
 \epsfig{figure=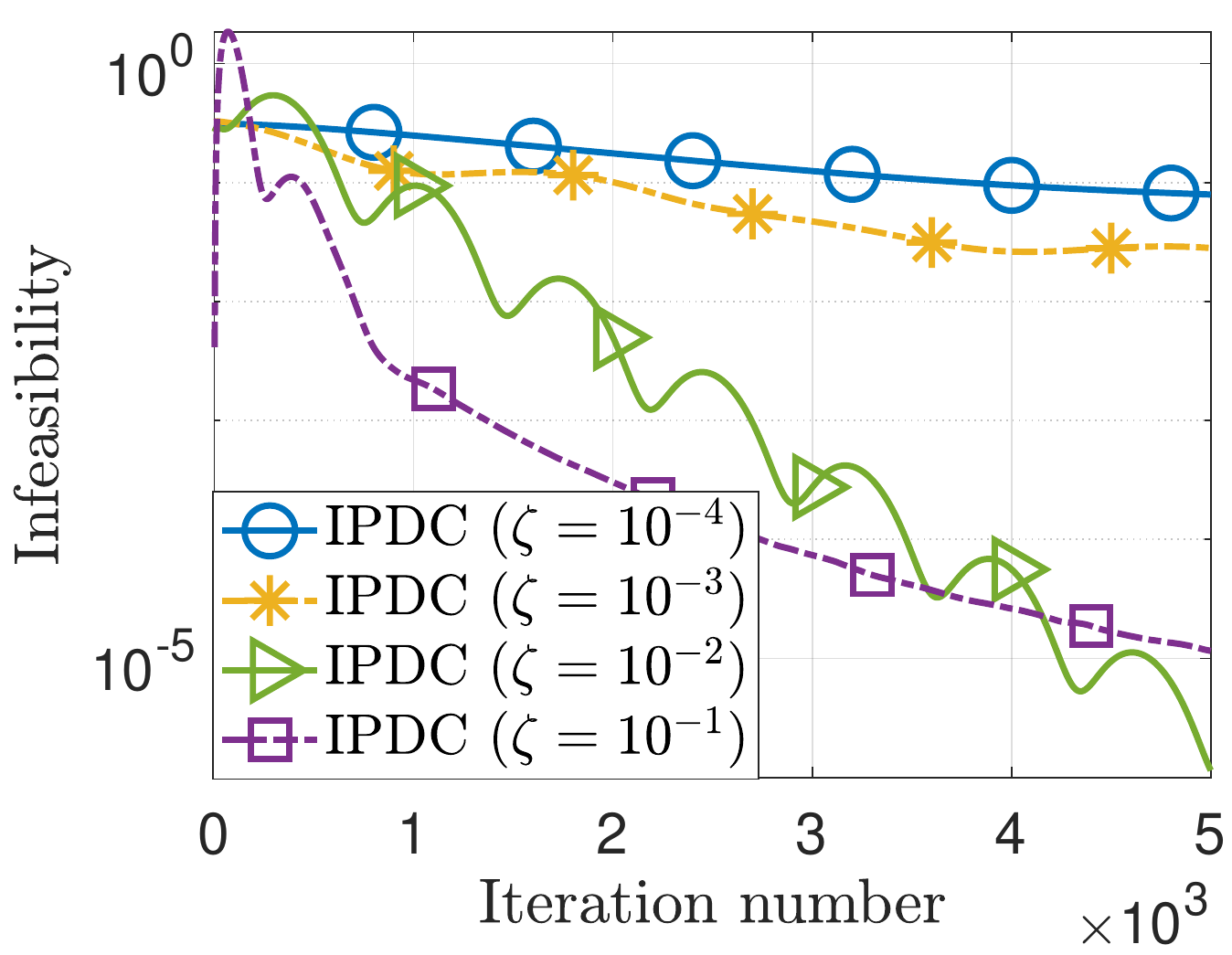,width=2.3in}
  \centerline{\small{(a) $\alpha = 0.01$, $ p= 10$, $\rho = 1$, $\beta = 0.1$, and various values of $\zeta$}}\medskip
\end{minipage}
\begin{minipage}[b]{1.0\linewidth}
  \centering
 \epsfig{figure=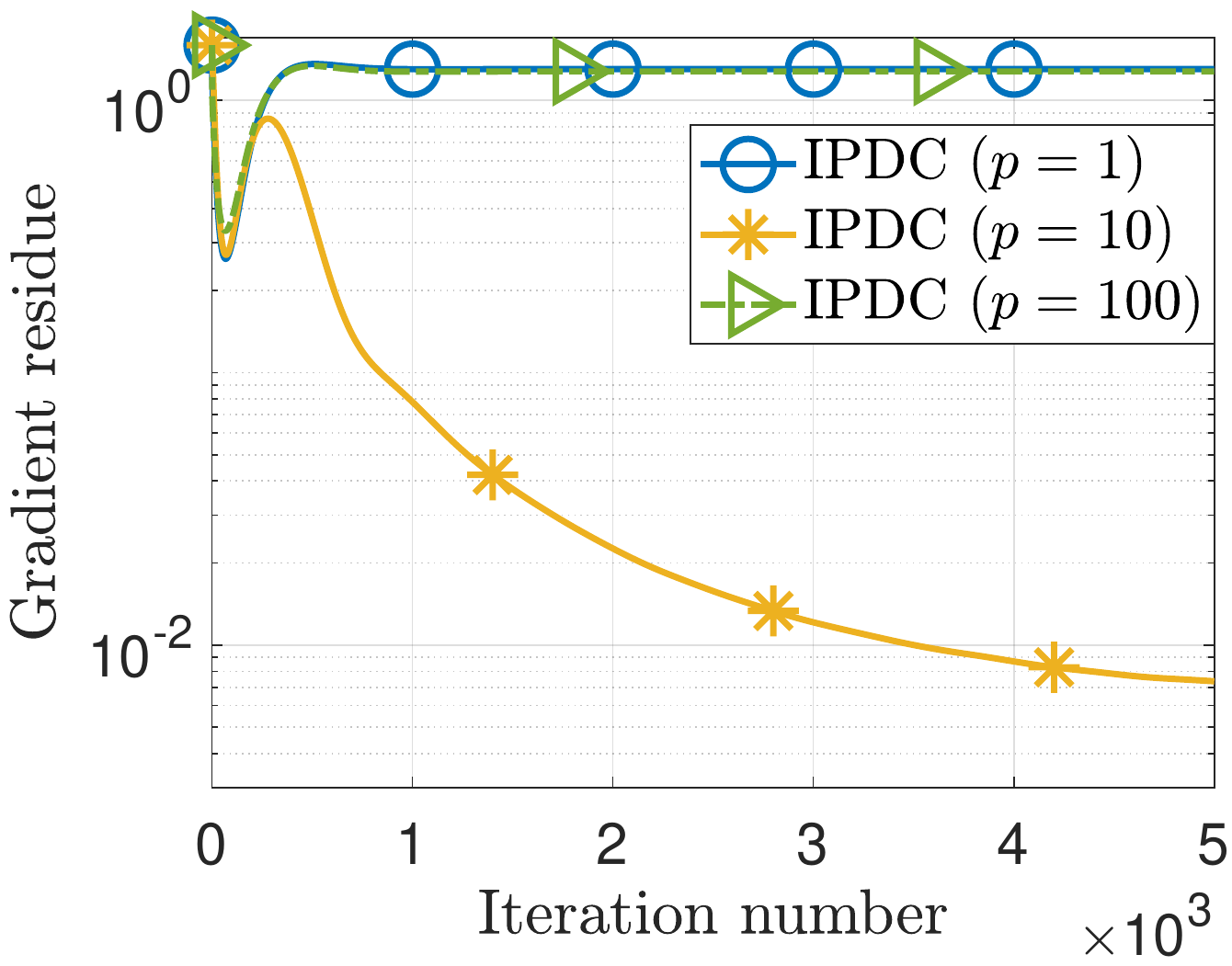,width=2.3in}
 \epsfig{figure=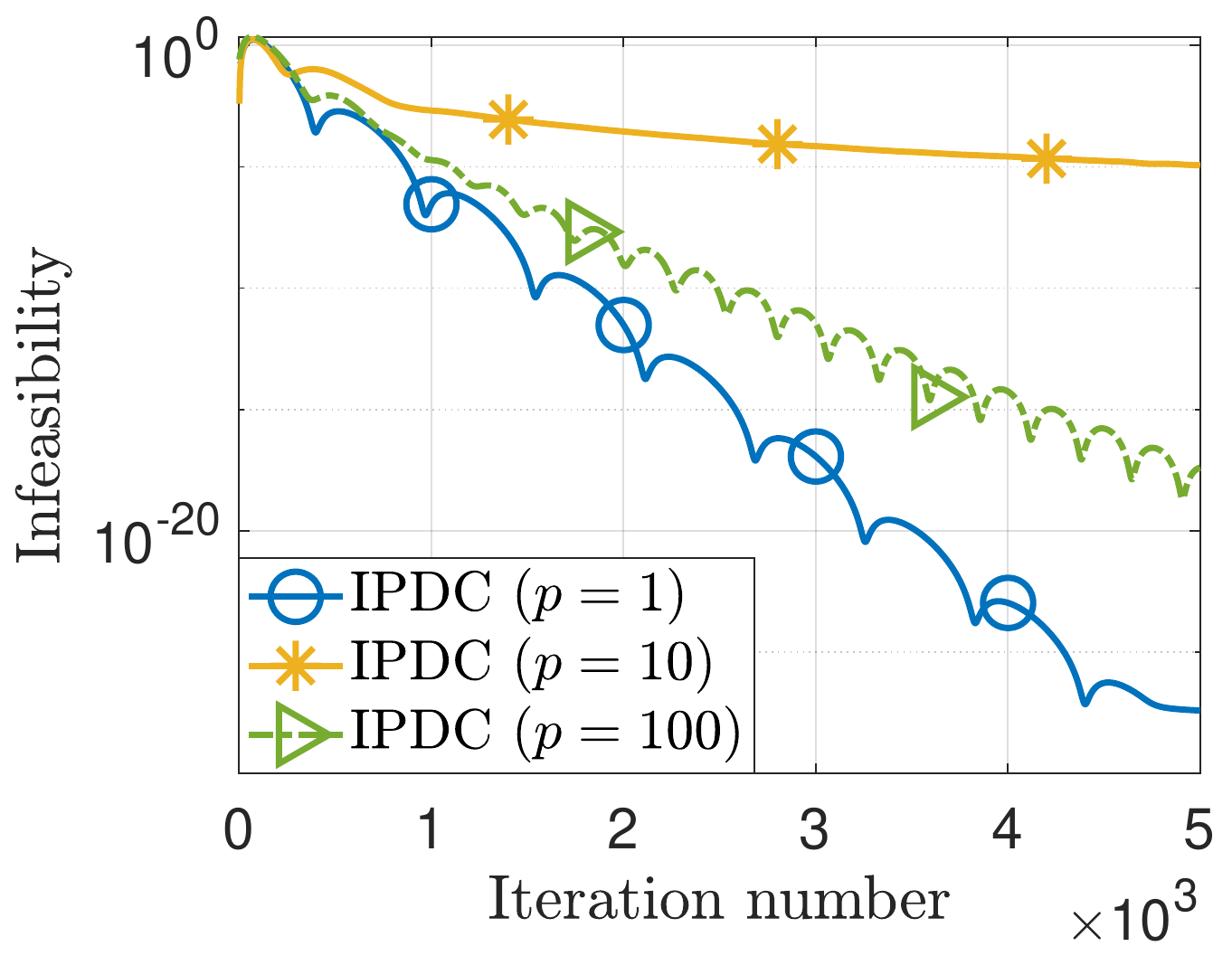,width=2.3in}
  \centerline{\small{(b) $\alpha = 0.01$, $\rho = 1$, $ \zeta = 0.1$, $\beta = 0.1$, and various values of $p$}}\medskip
\end{minipage}
\begin{minipage}[b]{1.0\linewidth}
  \centering
 \epsfig{figure=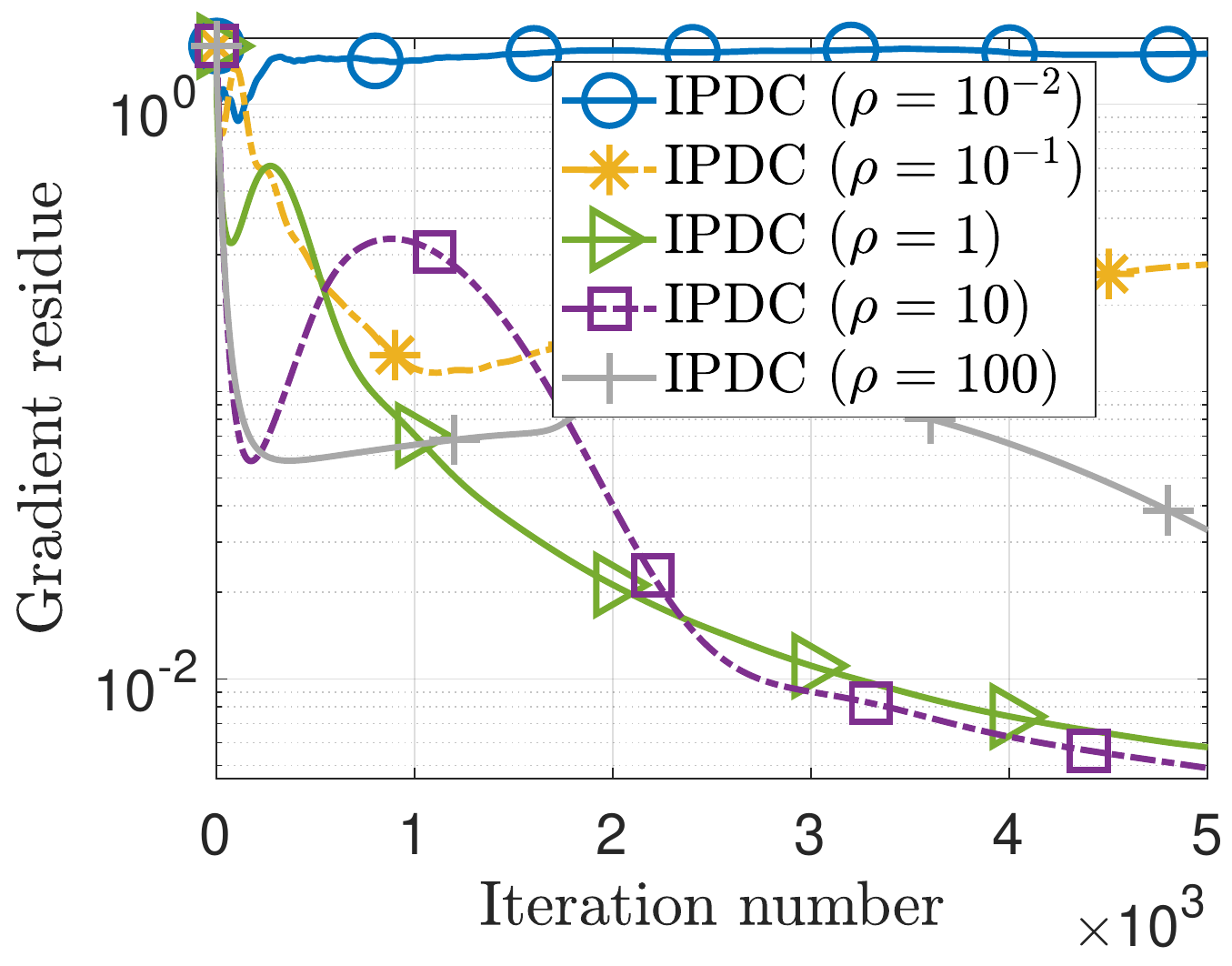,width=2.3in}
 \epsfig{figure=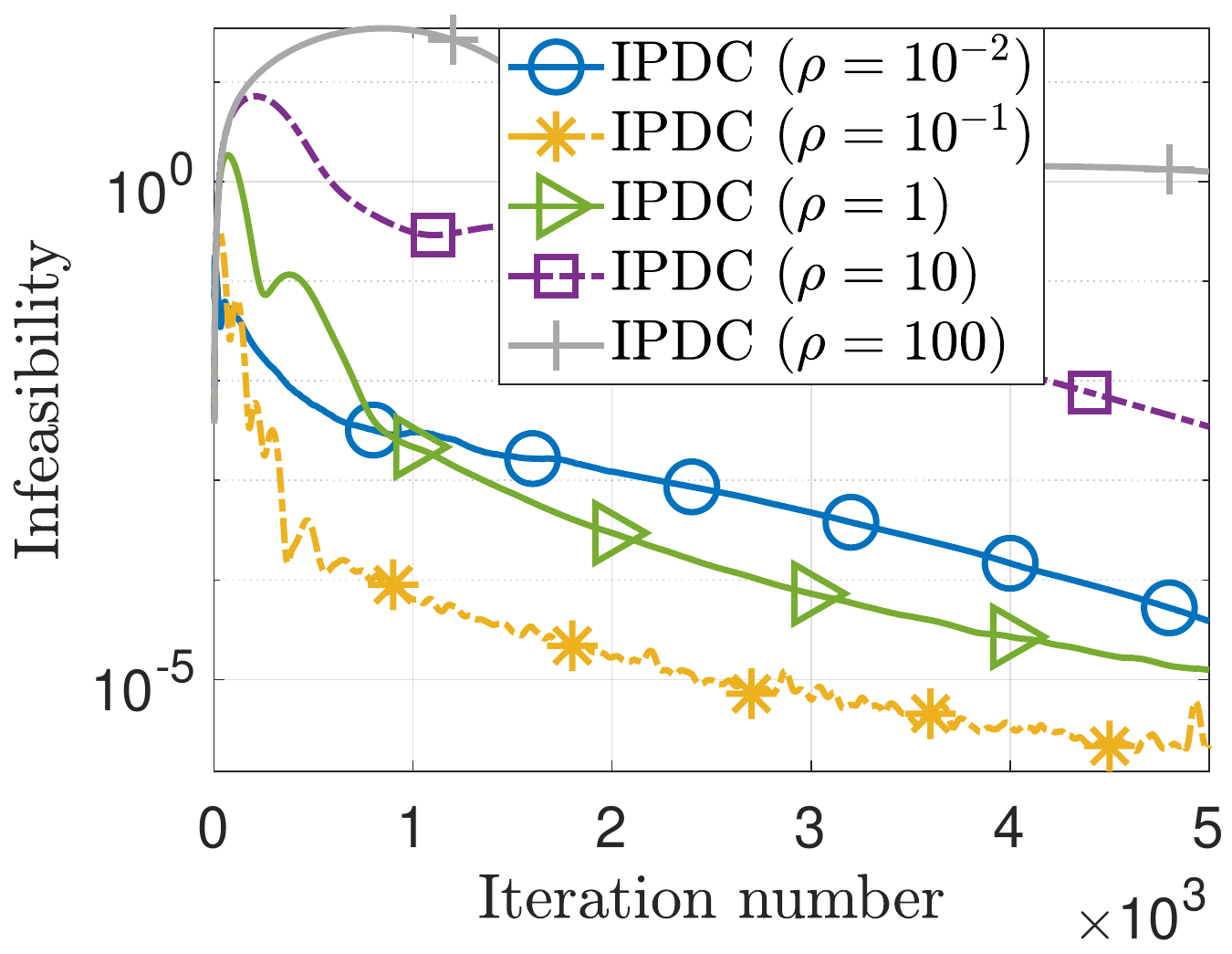,width=2.3in}
  \centerline{\small{(c) $\alpha = 0.01$, $ p= 10$, $\zeta = 0.1$, $ \beta = 0.1$, and various values of $\rho$}}\medskip
\end{minipage}
\caption{Convergence curves of the IPDC algorithm with respect to gradient residue and infeasibility.} \label{fig: LR IPDC p rho}\vspace{-0.4cm}
\end{figure}

In Fig.~\ref{fig: LR IPDC alpha beta}(a), we set $\beta = 0.1$, $ p = 10$, $\rho = 1$, $\zeta = 0.1$, and $\alpha \in \{10^{-4},  10^{-3}, 10^{-2}, 10^{-1} \}$. One can see that a larger $\alpha$ may speed up the convergence of gradient residue (left figure), while with distinct $\alpha$ the IPDC algorithm performs similar with respect to infeasibility (right figure).
In Fig.~\ref{fig: LR IPDC alpha beta}(b), we set $\alpha = 0.01$, $\rho = 1$, $\zeta = 0.1$, $p = 10$, and $\beta \in \{  10^{-3}, 10^{-2}, 10^{-1}, 1 \}$. we can see that $\beta$ influences the performance of the IPDC algorithm significantly. Interestingly, though Theorem 1
indicates that $\beta$ should be less than one and a small number, it is shown that
the algorithm behaves well with $\beta = 1$. This phenomenon is similar to that shown in Fig.~\ref{fig: LR PDC p rho}.

In Fig.~\ref{fig: LR IPDC p rho}(a), we set $\alpha = 0.01$, $ p= 10$, $\rho = 1$, $\beta = 0.1$, and $\zeta \in \{  10^{-4}, 10^{-3}, 10^{-2}, 10^{-1}\}$. One can observe that a larger step size $\zeta$ may speed up the convergence of gradient residue and infeasibility.
In Fig.~\ref{fig: LR IPDC p rho}(b), we set $\alpha = 0.01$, $\rho = 1$, $ \zeta = 0.1$, $\beta = 0.1$, and $p\in \{ 1, 10, 100\} $. One can observe that with a too small $p = 1$ or a too large $p=100$ the IPDC algorithm may not converge to a small gradient residue. Moreover, it is expected that a larger $p$ may slow down the convergence with respect to infeasibility, due to a more conservative convex approximation.
In Fig.~\ref{fig: LR IPDC p rho}(c), we set $\alpha = 0.01$, $ p= 10$, $\zeta = 0.1$, $ \beta = 0.1$, and $\rho  \in \{   10^{-2}, 10^{-1},  1, 10, 100\}$. One can see from these two figures that a larger $\rho$ can fasten the convergence of gradient residue, whereas a larger $\rho$ may slow down the speed satisfying the linear constraints.



\bibliographystyle{ieeetr}
{\footnotesize
\bibliography{distributed_opt}
}
\end{document}